\documentclass[12pt]{article}
\usepackage{graphics,graphicx}
\usepackage{stmaryrd}
\usepackage{amssymb,amsmath,amsopn,amsfonts}
\usepackage[dvips]{color}
\usepackage{colordvi,multicol}
\usepackage{epsf}
\newfam\msbfam
\font\tenmsb=msbm10 \textfont\msbfam=\tenmsb \font\sevenmsb=msbm7
\scriptfont\msbfam=\sevenmsb \font\fivemsb=msbm5
\scriptscriptfont\msbfam=\fivemsb

\def\th#1{\vspace{1mm}\noindent{\bf #1}\quad}

\def\proof{\vspace{1mm}\noindent{\it Proof}\quad}

\numberwithin{equation}{section}

\def\bc{\begin{center}}
\def\ec{\end{center}}
\def\no{\noindent}
\def\hang{\hangindent\parindent}
\def\textindent#1{\indent\llap{\qquad #1\ \ \enspace}\ignorespaces}
\def\ref{\par\hang\textindent}

\begin{document}
\title{ {\bf Three-dimensional Navier-Stokes equations driven by space-time white noise
\thanks{Research supported in part  by NSFC (No.11301026, No.11401019)  and DFG through IRTG 1132 and CRC 701, Key Lab of Random Complex Structures and
Data Science, Chinese Academy of Sciences (Grant No. 2008DP173182)}\\} }
\author{  {\bf Rongchan Zhu}$^{\mbox{a}}$, {\bf Xiangchan Zhu}$^{\mbox{b},}$\thanks{Corresponding author}
\date{}
\thanks{E-mail address:
zhurongchan@126.com(R. C. Zhu), zhuxc@bjtu.edu.cn(X. C. Zhu)}\\\\
$^{\mbox{a}}$Department of Mathematics, Beijing Institute of Technology, Beijing 100081,  China\\
$^{\mbox{b}}$School of Science, Beijing Jiaotong University, Beijing 100044, China  }

\maketitle

\noindent {\bf Abstract}

In this paper we prove existence and uniqueness of local solutions to the three dimensional (3D) Navier-Stokes (N-S) equation driven by space-time white noise  using two methods: first, the theory of regularity structures  introduced by Martin Hairer in \cite{Hai14} and second, the paracontrolled distribution proposed by Gubinelli, Imkeller, Perkowski in \cite{GIP13}. We also compare the two approaches.

\vspace{1mm}
\no{\footnotesize{\bf 2000 Mathematics Subject Classification AMS}:\hspace{2mm} 60H15, 82C28}
 \vspace{2mm}

\no{\footnotesize{\bf Keywords}:\hspace{2mm}  stochastic Navier-Stokes equation, regularity structure, paracontrolled distribution, space-time white noise, renormalisation}

\section{Introduction}
In this paper, we consider the three dimensional (3D) Navier-Stokes equation driven by space-time white noise:
 Recall that the Navier-Stokes equations describe the time
evolution of an incompressible fluid (see \cite{Te84}) and are given by
\begin{equation}\aligned\partial_t u+u\cdot \nabla u = &\nu \Delta u-\nabla p+\xi \\u(0)=&u_0, \quad div u = 0\endaligned\end{equation}
where $u( t,x)\in \mathbb{R}^3$ denotes the value of the velocity field at time $t$ and position
$x$, $p(t,x)$ denotes the pressure, and $\xi(t,x)$ is an external force field acting on
the fluid. We will consider the case when $x\in \mathbb{T}^3$, the three-dimensional torus.
Our mathematical model for the driving force $\xi$ is a Gaussian field which is
white in time and  space.

Random Navier-Stokes equations, especially the stochastic 2D Navier-Stokes equation driven by trace-class noise, have been studied in many articles (see e.g. \cite{FG95}, \cite{HM06}, \cite{De13}, \cite{RZZ14} and the reference therein). In the two dimensional case existence and uniqueness of  strong solutions have been obtained if the noisy forcing term is white in time and
colored in space. In the three dimensional case,  existence of martingale (=probabilistic weak) solutions, which form a Markov selection, have been constructed for  the stochastic 3D Navier-Stokes equation driven by trace-class noise in \cite{FR08}, \cite{DD03}, \cite{GRZ09}. Furthermore, the ergodicity has been obtained for every Markov selection of the martingale solutions if driven by non-degenerate trace-class noise (see \cite{FR08}).

This paper aims at giving a meaning to the equation (1.1) when $\xi$ is space-time white noise and at obtaining local (in time) solution.
Such a noise might not be relevant for the study of turbulence. However, in other cases, when a flow is subjected to an
external forcing with a very small time and space correlation length, a space-time white noise may be appropriate to model this situation. The main difficulty in this case is that $\xi$ and hence $u$ are so singular that the non-linear term is not well-defined.

In the two dimensional case, the Navier-Stokes equation driven by space-time white noise has been studied in \cite{DD02}, where a unique global solution in the (probabilistically) strong
sense has been obtained by using the Gaussian invariant measure for this equation. Thanks to the
incompressibility condition, we can write $u\cdot \nabla u=\frac{1}{2}\textrm{div} (u\otimes u)$.
The authors split the unknown into the solution to the
linear equation and  the solution to a modified version of the Navier-Stokes
equations:
$$\partial_t z=\nu\Delta z-\nabla \pi+\xi,\quad \textrm{div} z=0;$$
\begin{equation}\partial_t v=\nu\Delta v-\nabla q-\frac{1}{2}\textrm{div} [(v+z)\otimes (v+z)],\quad \textrm{div} v=0.\end{equation}
The first part $z$ is a  Gaussian process with non-smooth paths,  whereas the second part $v$ is smoother. The only term in   the nonlinear part, initially not well defined, is $z\otimes z$, which, however,  can be defined by using the Wick product. By a fixed point argument they obtain  existence and uniqueness of  local solutions in the two dimensional case. Then by using the Gaussian invariant measure for the 2D Navier-Stokes equation driven by space-time white noise,  existence and uniqueness of  (probabilistically) strong solutions starting from almost every initial value is obtained. (For the one-dimensional case we refer to \cite{DDT94}, \cite{RZZ14a}).

However, in the three dimensional case, the trick in the two dimensional case breaks down  since $v$ and  $z$ in (1.2) are so singular that not only $z\otimes z$ is not well-defined but also $v\otimes z$ and $v\otimes v$ have no meaning. Here $v$ is the solution to the nonlinear equation (1.2) and we cannot define these terms by using the Wick product.  As a result, we cannot make sense of  (1.2) and obtain  existence and uniqueness of  local solutions as in the two dimensional case. As a way out one might try to   iterate the above trick  as follows: we write $v=v_2+v_3$, where $v_2, v_3$ are the solutions to the following equations:
$$\partial_tv_2=\nu\Delta v_2-\nabla q_2-\frac{1}{2}\textrm{div} (z\otimes z),\quad \textrm{div} v_2=0,$$
\begin{equation}\partial_t v_3=\nu\Delta v_3-\nabla q_3-\frac{1}{2}\textrm{div} [(v_3+v_2)\otimes (v_3+v_2)]-\frac{1}{2}\textrm{div}((v_3+v_2)\otimes z)-\frac{1}{2}\textrm{div}(z\otimes (v_3+v_2)),\quad  \textrm{div} v_3=0.\end{equation}
Now we can make sense of the terms without $v_3$ in the right hand side of (1.3), hope $v_3$ becomes smoother such that the nonlinear terms including $v_3$ are well-defined and try to obtain a well-posed equation. However, this is not the case. For the unknown $v_3$ the nonlinear term on the right hand side of (1.3) including $v_3\otimes z$ is still not well-defined. Indeed, in this case $z\in \mathcal{C}^{-\frac{1}{2}-\kappa}$ for every $\kappa>0$. As a consequence, we cannot expect that  the regularity of $v_3$ is better than $\mathcal{C}^{\frac{1}{2}-\kappa}$ for every $\kappa>0$, which makes $v_3\otimes z$  not well-defined. No matter how many times we  modify this equation again as above, the equation  always contains the multiplication for the unknown and $z$, which is  not well-defined. Hence, this equation is ill-posed in the traditionally sense.

Thanks to the theory of regularity structures  introduced by Martin Hairer in \cite{Hai14} and the paracontrolled distribution proposed by Gubinelli, Imkeller and  Perkowski in \cite{GIP13} we can solve this problem and obtain  existence and uniqueness of  local solutions to the stochastic three dimensional Navier-Stokes equations driven by space-time white noise. Recently, these two approaches have been successful in giving a meaning to a lot of ill-posed stochastic PDEs like the Kardar-Parisi-Zhang (KPZ) equation (\cite{KPZ86}, \cite{BG97}, \cite{Hai13}), the dynamical $\Phi_3^4$ model (\cite{Hai14}, \cite{CC13}) and so on.  From a “philosophical” perspective, the theory of regularity structures and  the paracontrolled distribution  are inspired by the theory
of controlled rough paths \cite{Lyo98}, \cite{Gub04}, \cite{Hai11}. The main difference is that the regularity structure theory considers the problem locally, while the paracontrolled distribution method is a global approach  using Fourier analysis. For a comparison of these two methods we refer to Remark 3.13.

The key idea of the theory of regularity structures is as follows:  we perform an abstract Talyor expansion on  both sides of the equation. Originally Talyor expansions are only for functions. Here  the right objects, e.g. regularity structure that could
possibly take the place of Taylor polynomials, can be constructed. The regularity structure can  be endowed with a model $\iota\xi$, which is a concrete way of associating every element in  the abstract regularity structure to the actual Taylor polynomial at every point. Multiplication, differentiation, the state space of solutions,  and the convolution with singular kernels can be defined on this regularity structure, which is the major difficulty when trying to give a meaning to such singular stochastic partial differential equations as above. On the regularity   structure, a fixed point argument can be applied to obtain  local existence and uniqueness of the solution $\Phi$ to the equation lifted onto the regularity structure. Furthermore, we can go back to the real world with the help of another central tool of the theory, namely the reconstruction operator $\mathcal{R}$. If $\xi$ is a smooth process, $\mathcal{R}\Phi$ coincides with the classic solution to the equation. Now we have the following maps
$$\xi\mapsto\iota\xi\mapsto\Phi\mapsto\mathcal{R}\Phi,$$
and one is led to the following question: Given a sequence $\xi_\varepsilon$ of regularisations of the space-time white noise $\xi$, can we obtain the solution associated with $\xi$ by taking the limit of $\mathcal{R}\Phi_\varepsilon$, as $\varepsilon$ goes to $0$, where $\Phi_\varepsilon$ is the solution associated to $\xi_\varepsilon$. However, the answer to this question is no. Indeed, while the last two maps are continuous with respect to suitable topologies,   the above sequence $\iota\xi_\varepsilon$ of
canonical models fails to converge. It may, however,  still be possible to renormalize the model $\iota\xi_\varepsilon$ into some converging
model $\hat{\iota}\xi_\varepsilon$, which in turn can be related to a specific renormalised equation.

With these considerations in mind, let us go back to the 3D Navier-Stokes equations driven by space-time white noise.  We  apply Martin Hairer's regularity structure theory  to solve it. First, as in the two dimensional case we  write the nonlinear term $u\cdot \nabla u=\frac{1}{2}\textrm{div} (u\otimes u)$ and construct the associated regularity structure (Theorem 2.8). As in \cite{Hai14} we  construct different admissible models to denote different realizations of the equations corresponding to different noises. Then for any suitable models, we  obtain local existence and uniqueness of solutions by a fixed point argument. Finally, we renormalize the models associated with the approximations as mentioned above such that the  solution  to the equations associated with these renormalised models  converge to the solution to the 3D Navier-Stokes equation driven by space-time white noise in probability, locally in time.

The theory of paracontrolled distributions combines the idea of Gubinelli's controlled rough path \cite{Gub04} and Bony's paraproduct \cite{Bon81}, which is defined as  follows: Let $\Delta_jf$ be the jth Littlewood-Paley block of a distribution $f$ and define
$$\pi_<(f,g)=\pi_>(g,f)=\sum_{j\geq-1}\sum_{i<j-1}\Delta_if\Delta_jg, \quad\pi_0(f,g)=\sum_{|i-j|\leq1}\Delta_if\Delta_jg.$$ Formally $fg=\pi_<(f,g)+\pi_0(f,g)+\pi_>(f,g)$. Observing that, if $f$ is regular, $\pi_<(f,g)$ behaves like $g$ and is the only term in  Bony's paraproduct not increasing the regularity, the authors in \cite{GIP13}  consider a paracontrolled ansatz of  type
$$u=\pi_<(u',g)+u^\sharp,$$ where $\pi_<(u',g)$ represents the "bad-term" in the solution, $g$ is a functional of the Gaussian field and $u^\sharp$ is regular enough to allow the required multiplication. Then to make sense of the product  $uf$ we only need to define $gf$ by using a commutator estimate (Lemma 3.3).

In the second part of this paper we apply the paracontrolled distribution method to the 3D Navier-Stokes equations driven by space-time white noise. First we split the equation into four equations and consider the approximation equations. Here as in the theory of regularity structures, we still approximate $\xi$ by smooth functions $\xi_\varepsilon$ and obtain the  approximation equation associated with $\xi_\varepsilon$. By using the paracontrolled ansatz we obtain uniform estimates for the  approximation equations and moreover we also get the local Lipschitz continuity of solutions with respect to initial values and some extra terms $\mathbb{Z}(\xi_\varepsilon)$, which are independent of the solutions. These extra terms $\mathbb{Z}(\xi_\varepsilon)$ play a similar role as the models associated with the "distributional-like" elements in the abstract regularity structures. If $\mathbb{Z}(\xi_\varepsilon)$ converges to some $\mathbb{Z}$ in some suitable space, then the solution $u_\varepsilon$ associated with $\mathbb{Z}(\xi_\varepsilon)$ will converge to the desired solution. However, as in the theory of regularity structures,  we have to do suitable renormalisations for these terms such that they  converge in  suitable spaces. Here, inspired by \cite{Hai14}, we prove Lemma 3.10, which makes the calculations for the renormalisation  easier. Moreover  taking the limit of the solutions to the approximation equations we obtain local existence and uniqueness of the solutions. Indeed, by choosing a suitable solution space we can also give a meaning to the original equation (see Remark 3.9).

The main result of this article is the following theorem.
\vskip.10in

\th{Theorem 1.1} Let  $u_0\in \mathcal{C}^{\eta}$ for $\eta\in (-1,\alpha+2]$ with $\alpha\in(-\frac{13}{5},-\frac{5}{2})$. Let $\xi=(\xi^1,\xi^2,\xi^3)$, with $\xi^i, i=1,2,3$ being independent  white noises on $\mathbb{R}\times \mathbb{T}^3$, which we extend periodically to $\mathbb{R}^4$. Let $\rho:\mathbb{R}^4\rightarrow\mathbb{R}$ be a smooth compactly supported function with Lebesgue integral equal to $1$, and  symmetric with respect to space variable,  set $\rho_\varepsilon(t,x)=\varepsilon^{-5}\rho(\frac{t}{\varepsilon^2},\frac{x}{\varepsilon})$ and define $\xi_\varepsilon^i=\rho_\varepsilon*\xi^i$.  Consider the maximal  solution $u_\varepsilon$  to the following equation $$\partial_tu_{\varepsilon}^{i}=\Delta u_{\varepsilon}^i+\sum_{i_1=1}^3P^{ii_1}\xi_\varepsilon^{i_1}-\frac{1}{2}\sum_{i_1=1}^3P^{ii_1}(\sum_{j=1}^3D_j(u_{\varepsilon}^{i_1}u_{\varepsilon}^{j})),\quad u_\varepsilon(0)=Pu_0.$$
Then there exists $u\in C([0,\tau); \mathcal{C}^{\eta})$ and a sequence of random time $\tau_L$ converging to the explosion time $\tau$ of $u$ such that $$\sup_{t\in[0,\tau_L]}\|u^\varepsilon-u\|_{\eta}\rightarrow^P0.$$
\vskip.10in

\th{Remark 1.2}i) From Theorem 1.1 we know that although some diverging terms appear in the intermediate stages of the analysis, no renormalisation is actually necessary in (1.1).

ii)The results obtained by using paracontrolled distribution method are expressed a little bit differently (see Theorem 3.12).

\vskip.10in

This paper is organized as follows. In Section 2, we use the regularity structure theory to obtain local existence and uniqueness of  solutions
 to the 3D Navier-Stokes equations  driven by space-time white noise. In Section 3, we apply the paracontrolled distribution method to deduce local existence and uniqueness of  solutions. In Remark 3.13 we compare the two approaches.

\section{N-S equation by regularity structure theory}

\subsection{Preliminary on regularity structure theory}
In this subsection we recall some preliminaries for the  theory of regularity structures from \cite{Hai14} and \cite{HM15}. From this section we fixed a scaling $\mathfrak{s}=(\mathfrak{s}_0,1,...,1)$ of $\mathbb{R}^{d+1}$. We call $|\mathfrak{s}|=\mathfrak{s}_0+d$ scaling dimension.  We define the  associate  metric on $\mathbb{R}^{d+1}$  by
$$\|z-z'\|_\mathfrak{s}:=\sum_{i=0}^{d}|z_i-z'_i|^{1/\mathfrak{s}_i}.$$ For $k=(k_0,...,k_d)$ we define $|k|_\mathfrak{s}=\sum_{i=0}^d\mathfrak{s}_ik_i$.
\vskip.10in
\th{Definition 2.1} A regularity structure $\mathfrak{T}=(A,T,G)$ consists of the following elements:

(i) An index set $A\subset \mathbb{R}$ such that $0\in A$, $A$ is bounded from below and  locally finite.

(ii) A model space $T$, which is a graded vector space $T=\oplus_{\alpha\in A}T_\alpha$, with each $T_\alpha$ a Banach space. Furthermore, $T_0$ is one-dimensional and has a basis vector $\mathbf{1}$. Given $\tau\in T$ we write $\|\tau\|_\alpha$ for the norm of its component in $T_\alpha$.

(iii) A structure group $G$ of (continuous) linear operators acting on $T$ such that for every $\Gamma\in G$, every $\alpha\in A$ and every $\tau_\alpha\in T_\alpha$ one has
$$\Gamma\tau_\alpha-\tau_\alpha\in T_{<\alpha}:=\bigoplus_{\beta<\alpha}T_\beta.$$
Furthermore, $\Gamma\mathbf{1}=\mathbf{1}$ for every $\Gamma\in G$.
\vskip.10in
Now we have the regularity structure $\bar{T} =\bigoplus_{n\in\mathbb{N}}\bar{T}_n$ given by all polynomials in $d+1$ indeterminates, let us call them $X_0,..., X_d$, which denote the time and space directions respectively.  Denote $X^k=X_0^{k_0}\cdot\cdot\cdot X_d^{k_d}$ with $k$ a multi-index. In this case,   $A = \mathbb{N}$ and $\bar{T}_n$ denote the
space of monomials that are homogeneous of degree $n$.  The structure group can be defined by $\Gamma_hX^k=(X-h)^k$, $h\in\mathbb{R}^{d+1}$.

Given a smooth compactly supported test function $\varphi$, $x, y\in\mathbb{R}^d$, $\lambda>0$, we define
$$\varphi_x^\lambda(y)=\lambda^{-d}\varphi(\frac{y-x}{\lambda}).$$

Denote by $\mathcal{B}_r$ the set of smooth test functions $\varphi:\mathbb{R}^{d}\mapsto\mathbb{R}$ that are supported in the centred  ball of radius $1$ and such that their derivatives of order up to $r$ are uniformly bounded by $1$. We denote by $\mathcal{S}'$ the space of all distributions on $\mathbb{R}^{d}$. Now we give the definition of a model, which is a concrete way of associating every element in  the abstract regularity structure to the actual Taylor polynomial at every point.

For the Navier-Stokes equation we need to consider heat kernel composed with the Leray projection, which is not smooth on $\mathbb{R}^{d+1}\backslash\{0\}$. So we cannot apply [16, Lemma 5.5] directly. Instead we use  the inhomogeneous modelled distribution introduced in \cite{HM15}.
\vskip.10in

\th{Definition 2.2} Given a regularity structure $\mathfrak{T}$, an inhomogeneous model $(\Pi,\Gamma, \Sigma)$ consists of the following three elements:
\begin{itemize}
  \item  A collection of maps $\Gamma^t:\mathbb{R}^d\times \mathbb{R}^d\rightarrow \mathcal{G}$ parametrized by $t\in\mathbb{R}$, such that
  $$\Gamma^t_{xx}=1,\quad \Gamma^t_{xy}\Gamma^t_{yz}=\Gamma^t_{xz},$$
  for any $x,y,z\in\mathbb{R}^d$ and $t\in \mathbb{R}$, and the action of $\Gamma_{xy}^t$ on polynomials is given as above with $h=(0,y-x)$.
  \item A collection of maps $\Sigma_x:\mathbb{R}\times \mathbb{R}\rightarrow\mathcal{G}$, parametrized by $x\in\mathbb{R}^d$, such that, for any $x\in\mathbb{R}^d$ and $s,r,t\in\mathbb{R}$, one has
      $$\Sigma_x^{tt}=1, \quad \Sigma_x^{sr}\Sigma_x^{rt}=\Sigma_x^{st}, \quad \Sigma_x^{st}\Gamma^t_{xy}=\Gamma^s_{xy}\Sigma^{st}_y,$$
      and the action of $\Sigma_x^{st}$ on polynomials is given as above with $h=(t-s,0)$.
  \item A collection of linear maps $\Pi^t_x:T\rightarrow \mathcal{S}'$, such that
  $$\Pi^t_y=\Pi^t_x\Gamma^t_{xy}, \quad (\Pi^t_xX^{(0,\bar{k})})(y)=(y-x)^{\bar{k}},\quad  (\Pi^t_xX^{(k_0,\bar{k})})(y)=0,$$
  for all $x,y\in\mathbb{R}^d, t\in\mathbb{R}, \bar{k}\in\mathbb{N}^d, k_0\in\mathbb{N}$ such that $k_0>0$.
\end{itemize}
Moreover, for any $\gamma>0$ and every $T>0$, there is a constant $C$ for which the analytic bounds
$$|\langle \Pi^t_x\tau, \varphi^\lambda_x\rangle|\leq C\|\tau\|_l\lambda^l,\quad \|\Gamma_{xy}^t\tau\|_m\leq C\|\tau\|_l|x-y|^{l-m},  $$
$$\|\Sigma_x^{st}\tau\|_m\leq C\|\tau\|_l|t-s|^{(l-m)/\mathfrak{s}_0},$$
holds uniformly over all $\tau\in T_l$,  $l\in A$ with $l<\gamma$, all $m\in A$ such that $m<l$, and all test functions $\varphi\in \mathcal{B}_r$ with $r>-\inf A$, and all $t,s,\in[-T,T]$ and $x,y\in\mathbb{R}^d$ such that $|t-s|\leq 1$ and $|x-y|\leq 1$.

\vskip.10in
For a model $Z=(\Pi,\Gamma,\Sigma)$ we denote by $\|\Pi\|_{\gamma;T}, \|\Gamma\|_{\gamma;T}$ and $\|\Sigma\|_{\gamma;T}$ the smallest constants $C$ such that the bounds on $\Pi, \Gamma$ and $\Sigma$ in the above analytic bounds hold. Furthermore, we define
$$\interleave Z\interleave_{\gamma;T}:= \|\Pi\|_{\gamma;T}+\|\Gamma\|_{\gamma;T}+\|\Sigma\|_{\gamma;T}.$$
If $\bar{Z}=(\bar{\Pi}, \bar{\Gamma},\bar{\Sigma})$ is another model we define
$$\interleave Z;\bar{Z}\interleave_{\gamma;T}:=\|\Pi-\bar{\Pi}\|_{\gamma;T}+\|\Gamma-\bar{\Gamma}\|_{\gamma;T}+\|\Sigma-\bar{\Sigma}\|_{\gamma;T},$$
This
gives a natural topology for the space of all models for a given regularity structure. In the following we consider the models are periodic in space, which allows us to require the bounds to hold globally.

Now we have the following definition for the spaces of distributions $\mathcal{C}^\alpha$, $\alpha<0$, which is an extension of the definition of H\"{o}lder space to include $\alpha<0$.
\vskip.10in

\th{Definition 2.3} For $\alpha<0$, $\mathcal{C}^\alpha$ consists of $\eta\in\mathcal{S}'$, belonging to the dual space of the space $C^r_0$, i.e. the space of compactly supported $C^r$ functions, with $r>-\alpha+1$ and such that
$$\|\eta\|_{\alpha}:=\sup_{x\in\mathbb{R}^d}\sup_{\varphi\in\mathcal{B}_r}\sup_{\lambda\leq 1}\lambda^{-\alpha}|\eta(\varphi_x^\lambda)|<\infty.$$

\vskip.10in

On a bounded domain, $\mathcal{C}^\alpha$ coincides with the Besov space $B^\alpha_{\infty,\infty}$ defined in Section 3.
\vskip.10in
We also have the following definition of spaces of inhomogeneous modelled distributions, which are the  H\"{o}lder spaces on the regularity structure.

\vskip.10in
\th{Definition 2.4} Given a model $Z=(\Pi,\Gamma,\Sigma)$ for a regularity structure $\mathfrak{T}$   as above. Then for any $\gamma>0$ and $\eta\in\mathbb{R}$, the space $\mathcal{D}^{\gamma,\eta}$  consists of all functions $H:(0,T]\times\mathbb{R}^{d}\rightarrow \bigoplus_{\alpha<\gamma}T_\alpha$ such that
$$\interleave H\interleave_{\gamma,\eta;T}:=\|H\|_{\gamma,\eta;T}+\sup_{\tiny\aligned
s\neq t\in&(0,T]\\
|t-s|\leq& |t,s|_0^{\mathfrak{s}_0}
\endaligned }\sup_{x\in\mathbb{R}^d}
\sup_{l<\gamma}\frac{\|H_t(x)-\Sigma_x^{ts}H_s(x)\|_l}{|t-s|^{(\gamma-l)/\mathfrak{s}_0}|t,s|_0^{\eta-\gamma}}<\infty,$$
with
$$\| H\|_{\gamma,\eta;T}:=\sup_{t\in(0,T]}\sup_{x\in\mathbb{R}^d}\sup_{l<\gamma}|t|_0^{(l-\eta)\vee 0}\|H_t(x)\|_l+\sup_{t\in(0,T]}\sup_{\tiny\aligned
x\neq y\in&\mathbb{R}^d\\
|x-y|\leq&1
\endaligned }
\sup_{l<\gamma}\frac{\|H_t(x)-\Gamma_{xy}^{t}H_t(y)\|_l}{|x-y|^{\gamma-l}|t|_0^{\eta-\gamma}},$$
Here  we wrote $\|\tau\|_l$ for the norm of the component of $\tau$ in $T_l$ and $|t|_0:=|t|^{\frac{1}{\mathfrak{s}_0}}\wedge1$ and $|t,s|_0:=|t|_0\wedge|s|_0$.
\vskip.10in

For $H\in \mathcal{D}^{\gamma,\eta}$ and $\bar{H}\in \bar{\mathcal{D}}^{\gamma,\eta}$ (denoting by $\bar{\mathcal{D}}^{\gamma,\eta}$ the space built over another model $(\bar{\Pi},\bar{\Gamma},\bar{\Sigma})$), we also set
$$\aligned\| H;\bar{H}\|_{\gamma,\eta;T}:=&\sup_{t\in(0,T]}\sup_{x\in\mathbb{R}^d}\sup_{l<\gamma}|t|_0^{(l-\eta)\vee 0}\|H_t(x)-\bar{H}_t(x)\|_l\\&+\sup_{t\in(0,T]}\sup_{\tiny\aligned
x\neq y\in&\mathbb{R}^d\\
|x-y|\leq&1
\endaligned }
\sup_{l<\gamma}\frac{\|H_t(x)-\Gamma_{xy}^{t}H_t(y)-\bar{H}_t(x)+\bar{\Gamma}_{xy}^{t}\bar{H}_t(y)\|_l}{|x-y|^{\gamma-l}|t|_0^{\eta-\gamma}},\endaligned$$
$$\interleave H;\bar{H}\interleave_{\gamma,\eta;T}:=\|H;\bar{H}\|_{\gamma,\eta;T}+\sup_{\tiny\aligned
s\neq t\in&(0,T]\\
|t-s|\leq& |t,s|_0
\endaligned }\sup_{x\in\mathbb{R}^d}
\sup_{l<\gamma}\frac{\|H_t(x)-\Sigma_x^{ts}H_s(x)-\bar{H}_t(x)+\bar{\Sigma}^{ts}_x\bar{H}_s(x)\|_l}{|t-s|^{(\gamma-l)/\mathfrak{s}_0}|t,s|_0^{\eta-\gamma}},$$
which gives a natural
distance between elements $H\in \mathcal{D}^{\gamma,\eta}$ and $\bar{H}\in \bar{\mathcal{D}}^{\gamma,\eta}$.
\vskip.10in
Given a regularity structure, we say that a subspace $V\subset T$ is a sector of regularity $\alpha$ if it is invariant under the action of the structure group $G$ and it can be written as $V=\oplus_{\beta\in A}V_\beta$ with $V_\beta\subset T_\beta$, and $V_\beta=\{0\}$ for $\beta<\alpha$. We will use $\mathcal{D}^{\gamma,\eta}(V)$ to denote all functions in $\mathcal{D}^{\gamma,\eta}$ taking values in $V$.

On the regularity structure  a product $\star$ is a bilinear map on $T$ satisfying that for every $a\in T_\alpha$ and $b\in T_\beta$ one has $a\star b\in T_{\alpha+\beta}$ and $\mathbf{1}\star a=a\star \mathbf{1}=a$ for every $a\in T$. The product induces the pointwise product between modelled distribution under some conditions. For more details we refer to [16, Section 4].

The reconstruction theorem, which defines the so-called reconstruction operator, is one of the most fundamental result in
 the regularity structures theory.

\vskip.10in
\th{Theorem 2.5} (cf. [18, Theorem 2.11]) Given a model $Z=(\Pi,\Gamma,\Sigma)$ for a regularity structure $\mathfrak{T}$ with $\alpha:=\min A$ . Then for every $\eta\in\mathbb{R}, \gamma>0$ and $T>0$, there is a unique family of linear operators $\mathcal{R}_t:\mathcal{D}^{\gamma,\eta}\rightarrow \mathcal{C}^\alpha(\mathbb{R}^d)$, parametrised by $t\in(0,T]$, such that the bound
$$|\langle \mathcal{R}_tH_t-\Pi^t_xH_t(x),\varphi^\lambda_x\rangle|\lesssim \lambda^{\gamma}|t|_0^{\eta-\gamma}\|H\|_{\gamma,\eta;T}\|\Pi\|_{\gamma;T},$$
holds uniformly in $H\in\mathcal{D}^{\gamma,\eta}, t\in(0,T], x\in\mathbb{R}^d, \lambda\in (0,1]$ and $\varphi\in \mathcal{B}_r$ with $r>-\alpha+1$.

\vskip.10in
In order to define the integration against a space-time singular kernel $K$, Martin Hairer in \cite{Hai14} introduced an abstract integration map $\mathcal{I}:T\rightarrow T$ to provide an "abstract" representation of $\mathcal{K}$ operating at the level of the regularity structure. In the regularity structure theory  $\mathcal{I}$ is a linear map from $T$ to $T$ such that $\mathcal{I}T_\alpha\subset T_{\alpha+\beta}$ and $\mathcal{I}\bar{T}=0$ and for every $\Gamma\in G, \tau\in T$ one has $\Gamma \mathcal{I}\tau-\mathcal{I}\Gamma\tau\in\bar{T}$.

Furthermore, we say that $K$ is a  $\beta$-regularising kernel if one can write $K=\sum_{n\geq0}K_n$ where each  $K_n:\mathbb{R}^{d+1}\rightarrow\mathbb{R}$ is smooth and compactly supported in a ball of radius $2^{-n}$ around the origin. Moreover, we assume that for every multi-index $k$, one has a constant $C$ such that
$$\sup_z|D^kK_n(z)|\leq C2^{n(d+1-\beta+|k|_\mathfrak{s})},$$
holds uniformly in $n$. Finally, we assume that $\int K_n(z)E(z)dz=0$ for every polynomial $E$ of degree at most $r$ for some sufficiently large value of $r$.

We will write $K_t(x)=K(z)$, for $z=(t,x)$.
We say that a model $Z=(\Pi,\Gamma,\Sigma)$ realises $K$ for $\mathcal{I}$ if, for every $\alpha\in A$, every $\tau\in T_\alpha$ and every $x\in\mathbb{R}^d$, one has \begin{equation}\Pi_x^t(\mathcal{I}\tau+\mathcal{J}_{t,x}\tau)(y)=\int_{\mathbb{R}}\langle \Pi^s_x\Sigma^{st}_x\tau, K_{t-s}(y-\cdot)\rangle ds,\end{equation}
where $$\mathcal{J}_{t,x}\tau=\sum_{|k|_{\mathfrak{s}}<\alpha+\beta}\frac{X^k}{k!}\int_{\mathbb{R}}\langle \Pi^s_x\Sigma^{st}_x\tau, D^kK_{t-s}(x-\cdot)\rangle ds,$$
where $k\in\mathbb{N}^{d+1}$ and the derivative $D^k$ is in time-space. Moreover, we require that
\begin{equation}\Gamma^t_{xy}(\mathcal{I}+\mathcal{J}_{t,y})=(\mathcal{I}+\mathcal{J}_{t,x})\Gamma^t_{xy}, \quad \Sigma_x^{st}(\mathcal{I}+\mathcal{J}_{t,x})=(\mathcal{I}+\mathcal{J}_{s,x})\Sigma_x^{st},\end{equation}
for all $s,t\in\mathbb{R}$, and $x,y\in\mathbb{R}^d$.

Now we introduce the following operator acting on modelled distribution $H\in\mathcal{D}^{\gamma,\eta}$ with $\gamma+\beta>0$:
$$(\mathcal{K}_\gamma H)_t(x):=\mathcal{I}H_t(x)+\mathcal{J}_{t,x}H_t(x)+(\mathcal{N}_\gamma H)_t(x).$$
Here
$$(\mathcal{N}_\gamma H)_t(x):=\sum_{|k|_{\mathfrak{s}}<\gamma+\beta}\frac{X^k}{k!}\int_{\mathbb{R}}\langle \mathcal{R}_sH_s-\Pi^s_x\Sigma^{st}_xH_t(x), D^kK_{t-s}(x-\cdot)\rangle ds,$$
where $k\in\mathbb{N}^{d+1}$ and the derivative $D^k$ is in time-space.

Then we have the following results from [18, Theorem 2.21].
\vskip.10in
\th{Theorem 2.6} Let $\mathfrak{T}=(A,T,G)$ be a regularity structure with the minimal homogeneity $\alpha$. Let $K$ be a $\beta$-regularising kernel for some $\beta>0$, let $\mathcal{I}$ be an abstract integration map and let $Z=(\Pi,\Gamma,\Sigma)$ be a model realising $K$ for $\mathcal{I}$. Let $\gamma>0$, $\eta< \gamma<\eta+\mathfrak{s}_0, \eta>-\mathfrak{s}_0$. Then for $\gamma+\beta, \eta+\beta\notin\mathbb{N}$, $\mathcal{K}_\gamma$ maps $\mathcal{D}^{\gamma,\eta}$ into $\mathcal{D}^{\bar{\gamma},\bar{\eta}}$ with $\bar{\gamma}=\gamma+\beta$ and $\bar{\eta}=(\eta\wedge\alpha)+\beta$, and for any $H\in\mathcal{D}^{\gamma,\eta}$ the following bound holds
$$ \interleave \mathcal{K}_\gamma H\interleave_{\bar{\gamma},\bar{\eta};T}\lesssim \interleave H\interleave_{\gamma,\eta;T}\|\Pi\|_{\gamma;T}\|\Sigma\|_{\gamma;T}(1+\|\Gamma\|_{\bar{\gamma};T}+\|\Sigma\|_{\bar{\gamma};T}). $$
Furthermore, for every $t\in[0,T]$ one has
$$\mathcal{R}_t(\mathcal{K}_\gamma H )_t(x)=\int_0^t \langle \mathcal{R}_sH_s,K_{t-s}(x-\cdot)\rangle ds.$$
Let $\bar{Z}=(\bar{\Pi}, \bar{\Gamma},\bar{\Sigma})$ be another model realising $K$ for $\mathcal{I}$, which satisfies the same assumptions, and let $\bar{\mathcal{K}}_\gamma$ be defined as above for this model. Then one has
$$ \interleave \mathcal{K}_\gamma H; \bar{\mathcal{K}}_\gamma \bar{H}\interleave_{\bar{\gamma},\bar{\eta};T}\lesssim \interleave H;\bar{H}\interleave_{\gamma,\eta;T}+\interleave Z;\bar{Z}\interleave_{\bar{\gamma};T}, $$
for all $H\in \mathcal{D}^{\gamma,\eta}$ and $\bar{H}\in\bar{\mathcal{D}}^{\gamma,\eta}$. Here, the proportionality constant depends on $\interleave H\interleave_{\gamma,\eta;T}$, $\interleave \bar{H}\interleave_{\gamma,\eta;T}$ and the norms on the models $Z$ and $\bar{Z}$.

\vskip.10in
In order to deal with the Leray Projection, we have to consider convolution with the singular kernel for space variable. As in [16] we introduce an abstract integration map $\mathcal{I}_0:T\rightarrow T$ to provide an "abstract" representation of $\mathcal{P}$ operating at the level of the regularity structure. In the regularity structure theory  $\mathcal{I}_0$ is a linear map from $T$ to $T$ such that $\mathcal{I}_0T_\alpha\subset T_{\alpha}$ and $\mathcal{I}_0\bar{T}=0$ and for every $\Gamma\in G, \tau\in T$ one has $\Gamma \mathcal{I}_0\tau-\mathcal{I}_0\Gamma\tau\in\bar{T}$.

Furthermore, we say that $P$ is a  $0$-regularising kernel on $\mathbb{R}^d$ if one can write $P=\sum_{n\geq0}P_n$, where each  $P_n:\mathbb{R}^{d}\rightarrow\mathbb{R}$ is smooth and compactly supported in a ball of radius $2^{-n}$ around the origin. Furthermore, we assume that for every multi-index $k$, one has a constant $C$ such that
$$\sup_x|D^kP_n(x)|\leq C2^{n(d+|k|)},$$
holds uniformly in $n$. Finally, we assume that $\int P_n(x)E(x)dz=0$ for every polynomial $E$ of degree at most $r$ for some sufficiently large value of $r$.

We say that a model $Z=(\Pi,\Gamma,\Sigma)$ realises $P$ for $\mathcal{I}_0$ if, for every $\alpha\in A$, every $\tau\in T_\alpha$ and every $x\in\mathbb{R}^d$, one has \begin{equation}\Pi_x^t(\mathcal{I}_0\tau+\mathcal{J}^0_{t,x}\tau)(y)=\langle \Pi^t_x\tau, P(y-\cdot)\rangle,\end{equation}
where $$\mathcal{J}_{t,x}^0\tau=\sum_{|k|<\alpha}\frac{X^k}{k!}\langle \Pi^t_x\tau, D^kP(x-\cdot)\rangle ,$$
where $k\in\mathbb{N}^{d}$ and the derivative $D^k$ is in space. Moreover, we require that
\begin{equation}\Gamma^t_{xy}(\mathcal{I}_0+\mathcal{J}_{t,y}^0)=(\mathcal{I}_0+\mathcal{J}_{t,x}^0)\Gamma^t_{xy}, \quad \Sigma_x^{st}(\mathcal{I}_0+\mathcal{J}_{t,x}^0)=(\mathcal{I}_0+\mathcal{J}_{s,x}^0)\Sigma_x^{st},\end{equation}
for all $s,t\in\mathbb{R}$, and $x,y\in\mathbb{R}^d$.

Now we introduce the following operator acting on modelled distribution $H\in\mathcal{D}^{\gamma,\eta}$ with $\gamma+\beta>0$:
$$(\mathcal{P}_{\gamma,t} H_t)(x):=\mathcal{I}_0H_t(x)+\mathcal{J}_{t,x}^0H_t(x)+(\mathcal{N}^0_{\gamma,t} H_t)(x).$$
Here
$$(\mathcal{N}^0_{\gamma,t} H_t)(x):=\sum_{|k|<\gamma}\frac{X^k}{k!}\langle \mathcal{R}_tH_t-\Pi^t_xH_t(x), D^kP(x-\cdot)\rangle ,$$
where $k\in\mathbb{N}^{d}$ and the derivative $D^k$ is in space.

\th{Theorem 2.7} Let $\mathfrak{T}=(A,T,G)$ be a regularity structure with the minimal homogeneity $\alpha$. Let $P$ be a $0$-regularising kernel on $\mathbb{R}^d$, let $\mathcal{I}_0$ be an abstract integration map and let $Z=(\Pi,\Gamma,\Sigma)$ be a model realising $P$ for $\mathcal{I}_0$. Let $\gamma>0$, $\eta< \gamma, \eta>-\mathfrak{s}_0, r>\gamma-\alpha$. Then for $\gamma, \eta \notin \mathbb{N}$, $\mathcal{P}_\gamma$ maps $\mathcal{D}^{\gamma,\eta}$ into $\mathcal{D}^{{\gamma},{\eta}}$ , and for any $H\in\mathcal{D}^{\gamma,\eta}$ the following bound holds
$$ \interleave \mathcal{P}_\gamma H\interleave_{{\gamma},{\eta};T}\lesssim \interleave H\interleave_{\gamma,\eta;T}\|\Pi\|_{\gamma;T}\|\Sigma\|_{\gamma;T}(1+\|\Gamma\|_{{\gamma};T}+\|\Sigma\|_{{\gamma};T}). $$
Furthermore, for every $t\in[0,T]$ one has
$$\mathcal{R}_t(\mathcal{P}_{\gamma,t} H_t)(x)= \langle \mathcal{R}_tH_t,P(x-\cdot)\rangle .$$
Let $\bar{Z}=(\bar{\Pi}, \bar{\Gamma},\bar{\Sigma})$ be another model realising $P$ for $\mathcal{I}_0$, which satisfies the same assumptions, and let $\bar{\mathcal{P}}_\gamma$ be defined as above for this model.
Then one has
$$ \interleave \mathcal{P}_{\gamma} H; \bar{\mathcal{P}}_{\gamma} \bar{H}\interleave_{{\gamma},{\eta};T}\lesssim \interleave H;\bar{H}\interleave_{\gamma,\eta;T}+\interleave Z;\bar{Z}\interleave_{\gamma;T}, $$
for all $H\in \mathcal{D}^{\gamma,\eta}$ and $\bar{H}\in\bar{\mathcal{D}}^{\gamma,\eta}$. Here, the proportionality constant depends on $\interleave H\interleave_{\gamma,\eta;T}$, $\interleave \bar{H}\interleave_{\gamma,\eta;T}$ and the norms on the models $Z$ and $\bar{Z}$.

\proof The required bounds on the components of $(\mathcal{P}_{\gamma,t} H_t)(x)$ and $(\mathcal{P}_{\gamma,t} H_t)(y)-\Gamma^t_{yx}(\mathcal{P}_{\gamma,t} H_t)(x)$ as well as on the components of $(\mathcal{P}_{\gamma,t} H_t)(x)-\Sigma^{ts}_{x}(\mathcal{P}_{\gamma,s} H_s)(x)$ with non-integer homogeneities, can be obtained in exactly the same way as in [16, Prop. 6.16]. In the following we estimate the elements of $(\mathcal{P}_{\gamma,t} H_t)(x)-\Sigma^{ts}_{x}(\mathcal{P}_{\gamma,s} H_s)(x)$ with integer homogeneities: We have the identity
$$\aligned((\mathcal{P}_{\gamma,t} H_t)(x)-(\Sigma^{ts}_{x}\mathcal{P}_{\gamma,s} H_s)(x))_k=&((\mathcal{N}^0_{\gamma,t} H_t)(x))_k-((\Sigma^{ts}_{x}\mathcal{N}^0_{\gamma,s} H_s)(x))_k\\&+(\mathcal{J}_{t,x}^0(H_t(x)-\Sigma^{ts}_xH_s(x)))_k. \endaligned$$

We decompose $\mathcal{J}^0$ as $\mathcal{J}^0=\sum_{n\geq0}\mathcal{J}^{0,(n)}_{t,x}$ and $\mathcal{N}_\gamma^0$ as $\mathcal{N}_\gamma^0=\sum_{n\geq0}\mathcal{N}^{0,(n)}_{\gamma}$, where the nth term in each sum is obtained by replacing $P$ by $P_n$ in the expressions for $\mathcal{J}^0$ and $\mathcal{N}_\gamma^0$ respectively.
$$((\mathcal{N}^{0,(n)}_{\gamma,t} H_t)(x))_k=\frac{1}{k!}\langle \mathcal{R}_tH_t-\Pi^t_xH_t(x), D^kP_n(x-\cdot)\rangle.$$
$$((\Sigma^{ts}_{x}\mathcal{N}^{0,(n)}_{\gamma,s} H_s)(x))_k=\frac{1}{k!}\langle \mathcal{R}_sH_s-\Pi^s_xH_s(x), D^kP_n(x-\cdot)\rangle.$$
$$(\mathcal{J}_{t,x}^{0,(n)}(H_t(x)-\Sigma^{ts}_xH_s(x)))_k=\frac{1}{k!}\sum_{|k|<\zeta<\gamma}\langle \Pi^t_x\mathcal{Q}_\zeta (H_t(x)-\Sigma^{ts}_xH_s(x)), D^kP_n(x-\cdot)\rangle.$$
We first consider the case $2^{-n}\leq |t-s|^{\frac{1}{\mathfrak{s}_0}}\leq \frac{1}{2}|t,s|_0$:  by Theorem 2.5 we have for $|k|\leq \gamma$
$$|((\mathcal{N}^{0,(n)}_{\gamma,t} H_t)(x))_k|\lesssim |t|_0^{\eta-\gamma}2^{n(-\gamma+|k|)}\lesssim \sum_{\delta<0}|t,s|_0^{\eta-\gamma}|t-s|^{\frac{\gamma-|k|+\delta}{\mathfrak{s}_0}}2^{n\delta},$$
where we used the fact that $|t,s|_0\thicksim |t|_0$. The same bound also holds for $((\Sigma^{ts}_{x}\mathcal{N}^{0,(n)}_{\gamma,s} H_s)(x))_k$. Moreover, we obtain that
$$|(\mathcal{J}_{t,x}^{0,(n)}(H_t(x)-\Sigma^{ts}_xH_s(x)))_k|\lesssim \sum_{|k|<\zeta<\gamma}|t,s|_0^{\eta-\gamma}|t-s|^{\frac{\gamma-\zeta}{\mathfrak{s}_0}}2^{n(-\zeta+|k|)}
\lesssim\sum_{\delta<0}|t,s|_0^{\eta-\gamma}|t-s|^{\frac{\gamma-|k|+\delta}{\mathfrak{s}_0}}2^{n\delta}.$$
Regarding the corresponding term arising in $(\mathcal{P}_{\gamma,t} H_t)(x)-(\bar{\mathcal{P}}_{\gamma,t} \bar{H}_t)(x)$, we can use similar arguments as in the proof of [16, Theorem 5.12].

Now we consider the case $ |t-s|^{\frac{1}{\mathfrak{s}_0}}\leq 2^{-n}$. For this case we define
$$(\tilde{\Pi}_{(t,x)}\tau)(s,y):=(\Pi^s_x\Sigma^{st}_x\tau)(y), \quad \tilde{\Gamma}_{(t,x),(s,y)}:=\Gamma^t_{xy}\Sigma^{ts}_y=\Sigma^{ts}_x\Gamma^s_{xy}$$
for $\tau\in T$. By [18, Remark 2.7] we know that  the pair $(\tilde{\Pi},\tilde{\Gamma})$ is a model in the original sense of [16, Def. 2.17].
By [18, Remark 2.13]  $˜\tilde{R}H(t,\cdot):=R_tH_t(\cdot)$ is the reconstruction
operator for the model $(\tilde{\Pi},\tilde{\Gamma})$. One has the following identity:

$$\aligned&((\mathcal{P}_{\gamma,t} H_t)(x)-(\Sigma^{ts}_{x}\mathcal{P}_{\gamma,s} H_s)(x))_k\\=&\frac{1}{k!}\langle \mathcal{R}_tH_t-\mathcal{R}_sH_s-\Pi^t_x\Sigma^{ts}_xH_s(x)+\Pi^s_xH_s(x), D^kP_n(x-\cdot)\rangle\\&-\frac{1}{k!}\sum_{\zeta\leq |k|}\langle \Pi^t_x\mathcal{Q}_\zeta (H_t(x)-\Sigma^{ts}_xH_s(x)), D^kP_n(x-\cdot)\rangle:= T_1^k+T_2^k. \endaligned$$
By the same argument as above we have
$$|T_2^k|\lesssim \sum_{\zeta< |k|}|t,s|_0^{\eta-\gamma}|t-s|^{\frac{\gamma-\zeta}{\mathfrak{s}_0}}2^{n(-\zeta+|k|)}
\lesssim\sum_{\delta>0}|t,s|_0^{\eta-\gamma}|t-s|^{\frac{\gamma-|k|+\delta}{\mathfrak{s}_0}}2^{\delta n},$$
where the sum runs over a finite number of exponents.
In the following we consider $T_1^k$:
$$\aligned|T_1^k|=&\frac{1}{k!}|\langle \tilde{\mathcal{R}}H(t,\cdot)-\tilde{\mathcal{R}}H(s,\cdot)-\tilde{\Pi}_{(s,x)}H_s(x)(t,\cdot)+\tilde{\Pi}_{(s,x)}H_s(x)(s,\cdot), D^kP_n(x-\cdot)\rangle|\\=&\lim_{m\rightarrow\infty}\frac{1}{k!}|\sum_{(s_0,y)\in\Lambda_m^{\mathfrak{s}}}\langle \tilde{\Pi}_{(s_0,y)}H_{s_0}(y)-\tilde{\Pi}_{(s,x)}H_s(x),\varphi^{m,\mathfrak{s}}_{(s_0,y)}\rangle\langle \varphi^{m,\mathfrak{s}}_{(s_0,y)}(t)-\varphi^{m,\mathfrak{s}}_{(s_0,y)}(s), D^kP_n(x-\cdot)\rangle|,\endaligned$$
where $\varphi^{m,\mathfrak{s}}_{(s_0,y)}$ is the basis introduced in [16, Section 3], $\Lambda_m^\mathfrak{s}=\{\sum_{j=0}^3 2^{-m\mathfrak{s}_j}k_je_j:k_j\in\mathbb{Z}\},$
with $e_j$ denoting the jth element of the canonical basis of $\mathbb{R}^4$. Here in the second equality we used [16, Theorem 3.23] and the proof of [16, Theorem 3.10].
By the definition of the model we have
$$\aligned|\langle \tilde{\Pi}_{(s_0,y)}H_{s_0}(y)-\tilde{\Pi}_{(s,x)}H_s(x),\varphi^{m,\mathfrak{s}}_{(s_0,y)}\rangle|=&|\langle \tilde{\Pi}_{(s_0,y)}H_{s_0}(y)-\tilde{\Pi}_{(s_0,y)}\tilde{\Gamma}_{(s_0,y),(s,x)}H_s(x),\varphi^{m,\mathfrak{s}}_{(s_0,y)}\rangle|
\\\lesssim&|s_0,s|_0^{\eta-\gamma}\sum_{l<\gamma}
\|(s,x)-(s_0,y)\|_{\mathfrak{s}}^{\gamma-l}2^{-\frac{m|\mathfrak{s}|}{2}-lm}.\endaligned$$
For $\langle \varphi^{m,\mathfrak{s}}_{(s_0,y)}(t)-\varphi^{m,\mathfrak{s}}_{(s_0,y)}(s), D^kP_n(x-\cdot)\rangle$ we choose $m$ large enough such that $2^{-m}\leq |t-s|^{\frac{1}{\mathfrak{s}_0}}\leq 2^{-n}$. In this case by a similar calculation as in the proof of [16, Theorem 3.10] we know that
$$|\langle \varphi^{m,\mathfrak{s}}_{(s_0,y)}(t)-\varphi^{m,\mathfrak{s}}_{(s_0,y)}(s), D^kP_n(x-\cdot)\rangle|\lesssim 2^{n|k|}2^{-\frac{3m}{2}-rm}2^{n(3+r)}2^{\frac{\mathfrak{s}_0}{2}m}.$$
Furthermore, $|\langle \varphi^{m,\mathfrak{s}}_{(s_0,y)}(t)-\varphi^{m,\mathfrak{s}}_{(s_0,y)}(s), D^kP_n(x-\cdot)\rangle|=0$
unless $|x-y|\lesssim 2^{-n}$ and $|s-s_0|^{\frac{1}{\mathfrak{s}_0}}\wedge|t-s_0|^{\frac{1}{\mathfrak{s}_0}}\lesssim 2^{-m}$.

Hence we obtain that
$$\aligned|T_1^k|\lesssim&\lim_{m\rightarrow\infty}2^{3m}2^{-3n}|t,s|_0^{\eta-\gamma}\sum_{l<\gamma}
2^{-n(\gamma-l)}2^{-\frac{m|\mathfrak{s}|}{2}-lm}2^{n|k|}2^{-\frac{3m}{2}-rm}2^{n(3+r)}2^{\frac{\mathfrak{s}_0m}{2}}
\\\lesssim&\lim_{m\rightarrow\infty}|t,s|_0^{\eta-\gamma}\sum_{l<\gamma}
2^{-m(\gamma-l)}2^{-lm}2^{n|k|}2^{(\gamma-l-r)(m-n)}\\\lesssim&\sum_{\delta>0}|t,s|_0^{\eta-\gamma}|t-s|^{\frac{\gamma-|k|+\delta}{\mathfrak{s}_0}}2^{\delta n},\endaligned$$
where the sum runs over a finite number of exponents and in the first inequality we used the fact that  $|s_0,s|_0\thicksim |s,t|_0$ and the factor $2^{3m}2^{-3n}$ counts the number of non-zero terms appearing in the sum over $(s_0,y)$ and in the third inequality we used $r>\gamma-\alpha$. Taking summation over $n$ the required bound follows.
Regarding the corresponding term arising in $(\mathcal{P}_{\gamma,t} H_t)(x)-(\bar{\mathcal{P}}_{\gamma,t} \bar{H}_t)(x)$, we can use similar arguments as in the proof of [16, Theorem 5.12]. $\hfill\Box$

\subsection{N-S equation}

In this subsection we apply the regularity structure theory to the 3D Navier-Stokes equations on $\mathbb{T}^3$ driven by space-time white noise.
In this case we have the scaling $\mathfrak{s}=(2,1,1,1)$, so that the scaling dimension of space-time is $5$.
Since  the heat kernel $G$ is smooth on $\mathbb{R}^{4}\backslash\{0\}$ and has the scaling property $G(\frac{t}{\delta^2},\frac{x}{\delta})=\delta^3G(t,x)$ for $\delta>0$, by [16, Lemma 5.5] it can be decomposed into $K+R$ where $K$ is a $2$-regularising kernel and $R\in \mathcal{C}^\infty$

We know that the kernel $P^{ij}, i,j=1,2,3,$ for the Leray projection  is smooth on $\mathbb{R}^{3}\backslash\{0\}$ and has the scaling property $P^{ij}(\frac{x}{\delta})=\delta^3P^{ij}(x)$ for $\delta>0$, by [16, Lemma 5.5] it can be decomposed into $\bar{P}^{ij}+R_0^{ij}, i,j=1,2,3,$ where $\bar{P}^{ij}$ is a $0$-regularising kernel on $\mathbb{R}^3$ and $R_0^{ij}\in \mathcal{C}^\infty$.
 Define $$K^{ij}:=K*\bar{P}^{ij}.$$
By \cite{KT01} we have  $K^{ij}$ is of order $-3$, i.e. $|D^kK^{ij}(z)|\leq C\|z\|_\mathfrak{s}^{-3-|k|_\mathfrak{s}}$ for every $z$ with $\|z\|_\mathfrak{s}\leq 1$ and every multi-index $k$. We also use $D_lK^{ij}$, $l=1,2,3$, to represent the derivative of $K$ with respect to the $l$-th space variable and $D_lK^{ij}$ is of order $-4$.

Consider the regularity structure generated by the stochastic N-S equation with $\beta=2, -\frac{13}{5}<\alpha<-\frac{5}{2}$. In the regularity structure we use symbol the $\Xi_i$ to replace the driving noise $\xi^i$. We introduce the integration map $\mathcal{I}$  associated with $K$ and the integration map $\mathcal{I}_0^{ij}$ associated with $\bar{P}^{ij}$, which helps us to define $\mathcal{K}_\gamma$ and $\bar{\mathcal{P}}_{\gamma}^{ij}$. We also need
the integration maps $\mathcal{I}^{ii_1}_{0,k}, i,i_1=1,2,3, \mathcal{I}_k$ for a multiindex $k$, which  represents integration against $D^kP^{ii_1}, i,i_1=1,2,3, D^k K$ respectively.
We recall the following notations from \cite{Hai14}: defining a set $\mathcal{F}$ by postulating that $\{\mathbf{1},\Xi_i,X_j\}\subset \mathcal{F}$ and whenever $\tau,\bar{\tau}\in\mathcal{F}$, we have
$\tau\bar{\tau}\in\mathcal{F}$ and $\mathcal{I}^{ij}_{0,k}(\tau), \mathcal{I}_k(\tau)\in\mathcal{F}$; defining $\mathcal{F}_+$ as the set of all elements $\tau\in\mathcal{F}$ such that either $\tau=\mathbf{1}$ or $|\tau|_\mathfrak{s}>0$ and such that, whenever $\tau$ can be written as $\tau=\tau_1\tau_2$ we have either $\tau_i=\mathbf{1}$ or $|\tau_i|_\mathfrak{s}>0$; $\mathcal{H}, \mathcal{H}_+$ denote the sets of finite linear combinations of all elements in $\mathcal{F}, \mathcal{F}_+$, respectively.
Here for each $\tau\in \mathcal{F}$ a weight $|\tau|_{\mathfrak{s}}$  is obtained by setting
$|\mathbf{1}|_\mathfrak{s}=0$,
$$|\tau\bar{\tau}|_\mathfrak{s}=|\tau|_\mathfrak{s}+|\bar{\tau}|_\mathfrak{s},$$
for any two formal expressions $\tau$ and $\bar{\tau}$ in $\mathcal{F}$ such that
$$|\Xi_i|_\mathfrak{s}=\alpha,\quad |X_i|_\mathfrak{s}=\mathfrak{s}_i,\quad|\mathcal{I}_{k}(\tau)|_\mathfrak{s}=|\tau|_\mathfrak{s}+2-|k|_\mathfrak{s}, \quad |\mathcal{I}^{ii_1}_{0,k}(\tau)|_\mathfrak{s}=|\tau|_\mathfrak{s}-|k|_\mathfrak{s}.$$
To apply the regularity structure theory we write the equation as follows: for $i=1,2,3$
\begin{equation}\aligned\partial_t v_1^i=&\nu \Delta v_1^{i}+\sum_{i_1=1}^3P^{ii_1}\xi^{i_1},\quad \textrm{div} v_1=0,\\
\partial_t v^i=&\nu \Delta v^{i}-\sum_{i_1,j=1}^3P^{ii_1}\frac{1}{2}D_j [(v^{i_1}+v_1^{i_1}) (v^j+v_1^j)],\quad \textrm{div} v=0.\endaligned\end{equation}
Then $v_1+v$ is the solution to the 3D Navier-Stokes equations driven by space-time white noise.
Now we consider the second equation in (2.5). Define for $i,j,i_1=1,2,3$,
$$\mathcal{I}^{ij}:=\mathcal{I}_0^{ij}\mathcal{I},\quad \mathcal{I}^{ij}_{i_1}:=\mathcal{I}_0^{ij}\mathcal{I}_{i_1}$$
$$\mathfrak{M}_F^{ij}=\{1,\mathcal{I}^{ii_1}(\Xi_{i_1}),\mathcal{I}^{jj_1}(\Xi_{j_1}), \mathcal{I}^{ii_1}(\Xi_{i_1})\mathcal{I}^{jj_1}(\Xi_{j_1}), U_i, U_j, U_iU_j, \mathcal{I}^{ii_1}(\Xi_{i_1})U_j, U_i\mathcal{I}^{jj_1}(\Xi_{j_1}), i_1,j_1=1,2,3\}.$$
Then we build subsets $\{\mathcal{P}_n^i\}_{n\geq0}$, $\{\mathcal{U}_n\}_{n\geq0}$ and $\{\mathcal{W}_n\}_{n\geq0}$ by the following algorithm: For $i,j=1,2,3$, set $\mathcal{W}_0^{ij}=\mathcal{P}_0^i=\mathcal{U}_0=\varnothing$ and
$$\mathcal{W}_n^{ij}=\mathcal{W}_{n-1}^{ij}\cup\bigcup_{\mathcal{Q}\in \mathfrak{M}_F^{ij}}\mathcal{Q}(\mathcal{P}_{n-1}^i,\mathcal{P}_{n-1}^j),$$
$$\mathcal{P}_{n}^i=\{X^k\}\cup\{\mathcal{I}^{ii_1}_{i_2}(\tau):\tau\in\mathcal{W}_{n-1}^{i_1i_2}, i_1,i_2=1,2,3\},$$
$$\mathcal{U}_{n}=\{  \mathcal{I}_{i_2}(\tau):\tau\in\mathcal{W}_{n-1}^{i_1i_2},  i_1,i_2=1,2,3\},$$
and
$$\mathcal{F}_F:=\bigcup_{n\geq0}(\bigcup_{i,j=1}^3\mathcal{W}_n^{ij}\cup\mathcal{U}_{n} ), \quad \mathcal{F}_F^{ij}:=\bigcup_{n\geq0}\mathcal{W}_n^{ij}, i,j=1,2,3.$$
Then $\mathcal{F}_F$ contains the elements required to describe both the solution and the terms  in the equation (2.5).
We denote by $\mathcal{H}_F, \mathcal{H}_F^{ij}, i,j=1,2,3,$ the set of finite linear combinations of elements in $\mathcal{F}_F$, $\mathcal{F}_F^{ij}$, respectively.
\vskip.10in

\th{Remark} Here we construct $\mathcal{F}_F$  in a slightly different way from \cite{Hai14}. From (2.5) we observe that the integration map $\mathcal{I}^{ii_1}_j$ only acts on the elements belonging to $\mathcal{W}_n^{i_1j}$. The regularity structure  does not contain the elements belonging to $\mathcal{I}^{ii_1}_j(\mathcal{W}_n^{i_2j_1})$ for $(i_1,j)\neq (i_2,j_1)$ and $(i_1,j)\neq (j_1,i_2)$, which is enough for us to describe the solution and the equations.
\vskip.10in

Now we follow \cite{Hai14} to construct the structure group $G$. Define a linear projection operator $P_+:\mathcal{H}\rightarrow\mathcal{H}_+$ by imposing that
$$P_+\tau=\tau,\quad \tau\in\mathcal{F}_+,\quad P_+\tau=0,\quad \tau\in\mathcal{F}\setminus\mathcal{F}_+,$$
and two linear maps $\Delta:\mathcal{H}\rightarrow\mathcal{H}\otimes \mathcal{H}_+$ and $\Delta^+:\mathcal{H}_+\rightarrow\mathcal{H}_+\otimes \mathcal{H}_+$
by
$$\Delta\mathbf{1}=\mathbf{1}\otimes \mathbf{1},\quad \Delta^+\mathbf{1}=\mathbf{1}\otimes \mathbf{1},$$
$$\Delta X_i=X_i\otimes \mathbf{1}+\mathbf{1}\otimes X_i, \quad \Delta^+ X_i=X_i\otimes \mathbf{1}+\mathbf{1}\otimes X_i,$$
$$\Delta \Xi^i=\Xi^i\otimes \mathbf{1},$$
and recursively by $$\Delta(\tau\bar{\tau})=(\Delta\tau)(\Delta\bar{\tau})$$
$$\Delta(\mathcal{I}_k\tau)=(\mathcal{I}_{k}\otimes I)\Delta\tau+\sum_{l,m}\frac{X^l}{l!}\otimes \frac{X^m}{m!}(P_+\mathcal{I}_{k+l+m}\tau),$$
$$\Delta^+(\tau\bar{\tau})=(\Delta^+\tau)(\Delta^+\bar{\tau})$$
$$\Delta^+(\mathcal{I}_k\tau)=(I\otimes \mathcal{I}_{k}\tau)+\sum_{l}(P_+\mathcal{I}_{k+l}\otimes \frac{(-X)^l}{l!})\Delta\tau.$$
The above equalities still hold if $\mathcal{I}$ is replaced by $\mathcal{I}^{ij}_0$.

 By using the theory of regularity structures (see [16, Section 8]) we can  define a structure group $G_F$ of linear operators acting on  $\mathcal{H}_F$ satisfying Definition 2.1 as follows: For group-like elements $g\in \mathcal{H}_+^*$, the dual of $\mathcal{H}_+$,  $\Gamma_g:\mathcal{H}\rightarrow\mathcal{H}, \Gamma_g\tau=(I\otimes g)\Delta\tau$.  By [16, Theorem 8.24] we construct the following regularity structure.
\vskip.10in
\th{Theorem 2.8}  Let $T=\mathcal{H}_F$ with $T_\gamma=\langle \{\tau\in\mathcal{F}_F:|\tau|_\mathfrak{s}=\gamma\}\rangle$, $A=\{|\tau|_\mathfrak{s}:\tau\in\mathcal{F}_F\}$ and let $G_F$ be as above. Then $\mathfrak{T}_F=(A,\mathcal{H}_F, G_F)$ defines a regularity structure $\mathfrak{T}$. Furthermore,  $\mathcal{I}$ is an  abstract integration map of order $2$. For every $i,i_1=1,2,3$, $\mathcal{I}_0^{ii_1}$ is an abstract integration map of order $0$.

\proof In our case, the nonlinearity is locally subcritical. (i) (ii) in Definition 2.1 can be checked easily. (iii) in Definition 2.1 and the last results for $\mathcal{I}$ and $\mathcal{I}_0^{ii_1}$ follow from the definitions of $\Delta$ and $\Gamma_g$. $\hfill\Box$
\vskip.10in

We also endow $\mathfrak{T}_F$
with a natural commutative product $\star$ by setting $\tau\star\tau' = \tau\tau'$ for all basis vectors $\tau, \tau'$.

Now we come to construct suitable models associated with the regularity structure above.
Given any continuous approximation $\xi_\varepsilon$ to the driving noise $\xi$, we set for $s,t\in \mathbb{R}, x,y\in\mathbb{R}^3$
$$(\Pi_x^{(\varepsilon,t)}\mathcal{I}^{ii_1}(\Xi_{i_1}))(y)=K^{ii_1}*\xi^{i_1}_\varepsilon(t,y),$$
and
recursively define
$$(\Pi_x^{(\varepsilon,t)}\tau\bar{\tau})(y)=(\Pi_x^{(\varepsilon,t)}\tau)(y)(\Pi_x^{(\varepsilon,t)}\bar{\tau})(y),$$
\begin{equation}(\Sigma_x^{(\varepsilon,st)}\tau\bar{\tau})=(\Sigma_x^{(\varepsilon,st)}\tau)(\Sigma_x^{(\varepsilon,st)}\bar{\tau}),\quad (\Gamma_{xy}^{(\varepsilon,t)}\tau\bar{\tau})=(\Gamma_{xy}^{(\varepsilon,t)}\tau)(\Gamma_{xy}^{(\varepsilon,t)}\bar{\tau}).\end{equation}
For $\mathcal{I}\tau$ we define the actions of the maps  $(\Pi^{(\varepsilon)},\Gamma^{(\varepsilon)}, \Sigma^{(\varepsilon)})$ by (2.1) and (2.2). For $\mathcal{I}_0^{ij}$ we define the actions of the maps $(\Pi^{(\varepsilon)},\Gamma^{(\varepsilon)}, \Sigma^{(\varepsilon)})$ by (2.3) and (2.4) with $\mathcal{I}_0$ and $P$ replaced by $\mathcal{I}_0^{ij}$ and $\bar{P}^{ij}$. By this we can extend $(\Pi^{(\varepsilon)},\Gamma^{(\varepsilon)}, \Sigma^{(\varepsilon)})$ to the whole $\mathcal{H}_F$.

\vskip.10in
\th{Proposition 2.9}  $(\Pi^{(\varepsilon)},\Gamma^{(\varepsilon)}, \Sigma^{(\varepsilon)})$ is a model for the regularity structure $\mathfrak{T}_F$ constructed in Theorem 2.8.

\proof As in the proof of Theorem 2.7 we introduce the following model
$$(\tilde{\Pi}_{(t,x)}^{(\varepsilon)}\tau)(s,y)=(\Pi^{(\varepsilon,s)}_x\Sigma^{(\varepsilon,st)}_x\tau)(y), \quad \tilde{\Gamma}^{(\varepsilon)}_{(t,x),(s,y)}=\Gamma^{(\varepsilon,t)}_{xy}\Sigma^{{(\varepsilon,ts)}}_y=\Sigma^{{(\varepsilon,ts)}}_x\Gamma^{(\varepsilon,s)}_{xy},$$
which is a model in the original sense of [16, Def 2.17].
We can easily check that
$(\tilde{\Pi}_{(t,x)}^{(\varepsilon)}\mathcal{I}\tau)(s,y)$ and $\tilde{\Gamma}^{(\varepsilon)}_{(t,x),(s,y)}\mathcal{I}\tau $ coincide with the canonical model acting on $\mathcal{I}\tau$, which by [16, Prop. 8.27] implies the analytic and algebraic relations of the model $(\tilde{\Pi}^{(\varepsilon)},\tilde{\Gamma}^{(\varepsilon)})$ for the element $\mathcal{I}\tau$.
Since
$$(\Pi^{(\varepsilon,t)}_x\tau)(y)=(\tilde{\Pi}_{(t,x)}^{(\varepsilon)}\tau)(t,y), \quad \Gamma^{(\varepsilon,t)}_{xy}=\tilde{\Gamma}^{(\varepsilon)}_{(t,x),(t,y)}, \quad \Sigma^{{(\varepsilon,ts)}}_x=\tilde{\Gamma}^{(\varepsilon)}_{(t,x),(s,x)},$$
the analytic and algebraic relations also holds for the model $(\Pi^{(\varepsilon)},\Gamma^{(\varepsilon)}, \Sigma^{(\varepsilon)})$ acting on the element $\mathcal{I}\tau$. In the following we consider $\mathcal{I}_0^{ij}\tau$. The algebraic relations and the analytic bounds for $\Pi^{(\varepsilon,t)}_{x}\mathcal{I}_0^{ij}$ and $\Gamma^{(\varepsilon,t)}_{xy}\mathcal{I}_0^{ij}$ can be checked easily by similar arguments as in the proof of [16, Proposition 8.27]. In the following we will prove the bound on $\Sigma^{(\varepsilon)}\mathcal{I}_0^{ij}$. For $\tau\in T_l$, $k<l$ such that $k\notin \mathbb{N}$, (2.3) yields
$$\|\Sigma^{(\varepsilon,ts)}_x\mathcal{I}_0^{ij}\tau\|_k=\|\mathcal{I}_0^{ij}(\Sigma^{(\varepsilon,ts)}_x\tau)\|_k\leq \|\Sigma^{(\varepsilon,ts)}_x\tau\|_k\lesssim|s-t|^{\frac{l-k}{2}}.$$
 For $k\in\mathbb{N}^3$, by (2.4) we have the identity
$$(\Sigma^{(\varepsilon,ts)}_x\mathcal{I}_0\tau)_k=((\mathcal{J}_{t,x}^{0}\Sigma^{(\varepsilon,ts)}_x-\Sigma^{(\varepsilon,ts)}_x\mathcal{J}_{s,x}^{0})\tau)_k.$$
Here and in the following we omit superscript for $\mathcal{J}^0, \mathcal{I}_0$ and $\bar{P}$ for simplicity. We decompose $\mathcal{J}^0$ as $\mathcal{J}^0=\sum_{n\geq0}\mathcal{J}^{0,(n)}_{t,x}$, where the nth term in each sum is obtained by replacing $\bar{P}$ by $\bar{P}_n$ in the expressions for $\mathcal{J}^0$. Moreover, for $\tau\in T_l$
$$(\mathcal{J}^{0,(n)}_{t,x}\Sigma^{(\varepsilon,ts)}_x\tau)_k=\frac{1}{k!}\sum_{|k|<\zeta< \gamma}\langle \Pi^{(\varepsilon,t)}_x\mathcal{Q}_\zeta\Sigma_x^{(\varepsilon,ts)}\tau, D^k\bar{P}_n(x-\cdot)\rangle,$$
$$(\Sigma^{(\varepsilon,ts)}_x\mathcal{J}^{0,(n)}_{s,x}\tau)_k=\frac{1}{k!}\langle \Pi^{(\varepsilon,s)}_x\tau, D^k\bar{P}_n(x-\cdot)\rangle.$$

We first consider the case $2^{-n}\leq |t-s|^{\frac{1}{2}}$: by Definition 2.2 we have
$$|(\mathcal{J}^{0,(n)}_{t,x}\Sigma^{(\varepsilon,ts)}_x\tau)_k|\lesssim \sum_{|k|<\zeta< l}2^{n|k|}2^{-n\zeta}|s-t|^{\frac{l-\zeta}{2}}\lesssim \sum_{\delta<0}|t-s|^{\frac{l-|k|+\delta}{2}}2^{\delta n},$$
and
$$|(\Sigma^{(\varepsilon,ts)}_x\mathcal{J}^{0,(n)}_{s,x}\tau)_k|\lesssim 2^{n|k|}2^{-nl}\lesssim \sum_{\delta<0}|t-s|^{\frac{l-|k|+\delta}{2}}2^{\delta n}.$$
For the case that $ |t-s|^{\frac{1}{2}}\leq 2^{-n}$ we have
$$\aligned&(\mathcal{J}^{0,(n)}_{t,x}\Sigma^{(\varepsilon,ts)}_x\tau)_k-(\Sigma^{(\varepsilon,ts)}_x\mathcal{J}^{0,(n)}_{s,x}\tau)_k
\\=&-\frac{1}{k!}\sum_{\zeta\leq |k| }\langle \Pi^{(\varepsilon,t)}_x\mathcal{Q}_\zeta\Sigma^{(\varepsilon,ts)}\tau, D^k\bar{P}_n(x-\cdot)\rangle+\frac{1}{k!}\langle\Pi^{(\varepsilon,t)}_x\Sigma^{(\varepsilon,ts)}_x\tau- \Pi^{(\varepsilon,s)}_x\tau, D^k\bar{P}_n(x-\cdot)\rangle\\=&:T_1^k+T_2^k.\endaligned$$
For $T_1^k$ we have
$$|T_1^k|\lesssim \sum_{\zeta< |k|}2^{n|k|}2^{-n\zeta}|s-t|^{\frac{l-\zeta}{2}}\lesssim\sum_{\delta>0}|t-s|^{\frac{l-|k|+\delta}{2}}2^{\delta n},$$
where the sum runs over a finite number of exponents.
In the following we consider $T_2^k$:
$$\aligned|T_2^k|=&\frac{1}{k!}|\langle \tilde{\Pi}_{(s,x)}\tau(t,\cdot)-\tilde{\Pi}_{(s,x)}\tau(s,\cdot), D^k\bar{P}_n(x-\cdot)\rangle|\\=&\lim_{m\rightarrow\infty}\frac{1}{k!}|\sum_{(s_0,y)\in\Lambda_m^{\mathfrak{s}}}\langle \tilde{\Pi}_{(s,x)}\tau,\varphi^{m,\mathfrak{s}}_{(s_0,y)}\rangle\langle \varphi^{m,\mathfrak{s}}_{(s_0,y)}(t)-\varphi^{m,\mathfrak{s}}_{(s_0,y)}(s), D^kP_n(x-\cdot)\rangle|\\\lesssim&\lim_{m\rightarrow\infty}2^{3m}2^{-3n}
2^{-\frac{m|\mathfrak{s}|}{2}-l m}2^{n|k|}2^{-\frac{3m}{2}-rm}2^{n(3+r)}2^m
\\\lesssim&\lim_{m\rightarrow\infty}2^{-l m}2^{n|k|}2^{-r(m-n)}\\\lesssim&\sum_{\delta>0}|t-s|^{\frac{l-|k|+\delta}{2}}2^{\delta n},\endaligned$$
where the sum runs over a finite number of exponents and in the first inequality we used similar arguments as in the proof of Theorem 2.7 and the factor $2^{3m}2^{-3n}$ counts the number of non-zero terms appearing in the sum over $(s_0,y)$ and in the last inequality we choose $m$ large enough such that $2^{-m}<|t-s|^{\frac{1}{2}}$.
Taking the sum over $n$ we obtain the desired bounds for $\|\Sigma^{(\varepsilon,ts)}_x\mathcal{I}_0\tau\|_k$ and the result follows.$\hfill\Box$

\vskip.10in
\th{Definition 2.10} A model $(\Pi,\Gamma,\Sigma)$ for $\mathfrak{T}$ is admissible if it satisfies  (2.6) and furthermore realizes $K, \bar{P}^{ij}, i,j=1,2,3,$ for $\mathcal{I}, \mathcal{I}_0^{ij}$ respectively. We denote by $\mathcal{M}_F$ the set of admissible models.
\vskip.10in

Set
$$\aligned\mathcal{F}_0=&\{\mathbf{1},  \mathcal{I}^{ii_1}(\Xi_{i_1}), \mathcal{I}^{ii_1}(\Xi_{i_1})\mathcal{I}^{jj_1}(\Xi_{j_1}), \mathcal{I}^{ii_1}_j(\mathcal{I}^{i_1i_2}(\Xi_{i_2})),\mathcal{I}^{ii_1}_j(\mathcal{I}^{i_1i_2}(\Xi_{i_2})\mathcal{I}^{jj_1}(\Xi_{j_1})), \mathcal{I}^{ii_1}_j(\mathcal{I}^{jj_1}(\Xi_{j_1})),\\&\mathcal{I}^{ii_1}_k(\mathcal{I}^{i_1i_2}(\Xi_{i_2}))\mathcal{I}^{jj_1}(\Xi_{j_1}),
\mathcal{I}^{ii_1}_k(\mathcal{I}^{kk_1}(\Xi_{k_1}))\mathcal{I}^{jj_1}(\Xi_{j_1}),
\mathcal{I}^{ii_1}_k(\mathcal{I}^{i_1i_2}(\Xi_{i_2})\mathcal{I}^{kk_1}(\Xi_{k_1}))\mathcal{I}^{jj_1}(\Xi_{j_1}),\\&
\mathcal{I}^{ii_1}_k(\mathcal{I}^{i_1i_2}(\Xi_{i_2})\mathcal{I}^{kk_1}(\Xi_{k_1}))
\mathcal{I}^{jj_1}_l(\mathcal{I}^{j_1j_2}(\Xi_{j_2})\mathcal{I}^{ll_1}(\Xi_{l_1})),
\mathcal{I}^{ii_1}_l(\mathcal{I}^{i_1i_2}_k(\mathcal{I}^{i_2i_3}(\Xi_{i_3})\mathcal{I}^{kk_1}(\Xi_{k_1}))
\mathcal{I}^{ll_1}(\Xi_{l_1}))\mathcal{I}^{jj_1}(\Xi_{j_1}),\\&
\mathcal{I}^{ii_1}_l(\mathcal{I}^{ll_1}_k(\mathcal{I}^{l_1l_2}(\Xi_{l_2})\mathcal{I}^{kk_1}(\Xi_{k_1}))\mathcal{I}^{i_1i_2}(\Xi_{i_2}))\mathcal{I}^{jj_1}
(\Xi_{j_1}), \mathcal{I}_j(\mathcal{I}^{ii_1}(\Xi_{i_1})\mathcal{I}^{jj_1}(\Xi_{j_1})),
 \\&i,j,k,l,i_1,i_2,i_3,j_1,j_2,k_1,l_1,l_2=1,2,3 \},\endaligned$$
 and
$$\mathcal{F}_*=\{\mathcal{I}^{ik}(\Xi_k), \mathcal{I}^{ii_1}_k(\mathcal{I}^{i_1i_2}(\Xi_{i_2})\mathcal{I}^{kk_1}(\Xi_{k_1}))\mathcal{I}^{jj_1}(\Xi_{j_1}), i,k, i_1,i_2,j,j_1,k_1=1,2,3\}.$$
To make our paper more readable we use the tree notation from \cite{Hai14} to explain the complicated elements in $\mathcal{F}_0$. However, unlike as in the $\Phi^4_3$ case, the solution to the  stochastic N-S equation is vector valued and there are a lot of superscripts and subscripts for the elements in $\mathcal{F}_0$, which  will not be noticeable in the tree notation. The tree notation only helps us to make the complicated notation  clearer.

For $\Xi$ we simply draw a dot. The integration map $\mathcal{I}^{ij}$ is then represented by a downfacing line while the integration map  $\mathcal{I}_0\mathcal{I}_j$ is then represented by a downfacing dotted line. The integration map $\mathcal{I}_j$ is represented by $\includegraphics[height=0.3cm]{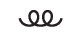}$. The multiplication of symbols is obtained by joining them at the root.

$$\aligned\mathcal{F}_0=\{\mathbf{1}, &\includegraphics[height=0.5cm]{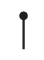}, \includegraphics[height=0.5cm]{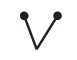},\includegraphics[height=0.7cm]{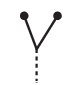}, \includegraphics[height=0.7cm]{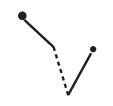},\includegraphics[height=0.7cm]{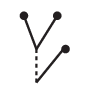},
\includegraphics[height=0.7cm]{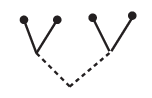},\includegraphics[height=0.7cm]{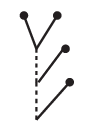}, \includegraphics[height=0.7cm]{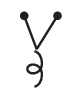}\},\endaligned$$
$$\aligned\mathcal{F}_*=&\{ \includegraphics[height=0.5cm]{01.eps}, \includegraphics[height=0.7cm]{06.eps}\}.\endaligned$$
We choose $\alpha\in(-\frac{13}{5},-\frac{5}{2})$ and the reason for $\alpha>-\frac{13}{5}$ is that this is precisely the value of $\alpha$ at which the homogeneity of the term $\mathcal{I}_j(\tau)\mathcal{I}(\Xi)$ vanishes for $\tau=\includegraphics[height=0.7cm]{08.eps}$.

Then $\mathcal{F}_0\subset \mathcal{F}_F$ contains every $\tau\in\mathcal{F}_F$ with $|\tau|_\mathfrak{s}\leq0$ and for every $\tau\in\mathcal{F}_0$, $\Delta\tau\in\langle\mathcal{F}_0\rangle\otimes \langle \rm{Alg}(\mathcal{F}_*)\rangle$. Here $\langle\mathcal{F}_0\rangle$ denotes the linear span of $\mathcal{F}_0$ and $\rm{Alg}(\mathcal{F}_*)$ denotes the set of all elements in $\mathcal{F}_+$ of the form $X^k\prod_{i,i_1,i_2}\mathcal{I}_{l_i}^{i_1i_2}\tau_i$ for some multiindices $k$ and $l_i$ such that $|\mathcal{I}^{ii_1}_{l_i}\tau_i|_\mathfrak{s}>0$ and $\tau_i\in\mathcal{F}_*$.

As mentioned in the introduction, we should do renormalisations for the model $(\Pi^\varepsilon,\Gamma^\varepsilon,\Sigma^\varepsilon)$ built from $\xi_\varepsilon$ such that it converges as $\xi_\varepsilon\rightarrow\xi$ in a suitable sense. In the theory of regularity structure, this has been transferred to find a sequence of $M_\varepsilon$ belonging to the renormalisation group $\mathfrak{R}_0$ defined in [16, Definition 8.43] such that $M_\varepsilon(\Pi^\varepsilon,\Gamma^\varepsilon,\Sigma^\varepsilon)$ converges to a finite limit. In the following we use the notations and definitions in [16, Section 8.3] and follow Hairer's idea to define $M$. We also use the tree notation as above to make it  clearer.

For  constants $C^1_{ii_1jj_1}, C^2_{ii_1i_2jj_1j_2kk_1ll_1}, C^3_{ii_1i_2i_3kk_1ll_1jj_1}, C^4_{ii_1i_2kk_1ll_1l_2jj_1},  i,j,k,l, i_1,i_2,i_3,j_1,$ $k_1,l_1,l_2=1,2,3$, we define a linear map $M$ on $\langle\mathcal{F}_0\rangle$ by
$$\aligned &M(\mathcal{I}^{ii_1}(\Xi_{i_1})\mathcal{I}^{jj_1}(\Xi_{j_1}))=\mathcal{I}^{ii_1}(\Xi_{i_1})\mathcal{I}^{jj_1}(\Xi_{j_1})-C^1_{ii_1jj_1}\mathbf{1},\\
&M\includegraphics[height=0.5cm]{02.eps}=\includegraphics[height=0.5cm]{02.eps}-C^1_{ii_1jj_1}\mathbf{1},\\
 &M(\mathcal{I}^{ii_1}_k(\mathcal{I}^{i_1i_2}(\Xi_{i_2})\mathcal{I}^{kk_1}(\Xi_{k_1}))
\mathcal{I}^{jj_1}_l(\mathcal{I}^{j_1j_2}(\Xi_{j_2})\mathcal{I}^{ll_1}(\Xi_{l_1})))\\=&\mathcal{I}^{ii_1}_k(\mathcal{I}^{i_1i_2}(\Xi_{i_2})\mathcal{I}^{kk_1}(\Xi_{k_1}))
\mathcal{I}^{jj_1}_l(\mathcal{I}^{j_1j_2}(\Xi_{j_2})\mathcal{I}^{ll_1}(\Xi_{l_1}))-C^2_{ii_1i_2jj_1j_2kk_1ll_1}\mathbf{1},
\\&M\includegraphics[height=0.7cm]{07.eps}=\includegraphics[height=0.7cm]{07.eps}-C^2_{ii_1i_2jj_1j_2kk_1ll_1}\mathbf{1},\\ &M(\mathcal{I}^{ii_1}_l(\mathcal{I}^{i_1i_2}_k(\mathcal{I}^{i_2i_3}(\Xi_{i_3})\mathcal{I}^{kk_1}(\Xi_{k_1}))
\mathcal{I}^{ll_1}(\Xi_{l_1}))\mathcal{I}^{jj_1}(\Xi_{j_1}))\\=&\mathcal{I}^{ii_1}_l(\mathcal{I}^{i_1i_2}_k(\mathcal{I}^{i_2i_3}(\Xi_{i_3})\mathcal{I}^{kk_1}(\Xi_{k_1}))
\mathcal{I}^{ll_1}(\Xi_{l_1}))\mathcal{I}^{jj_1}(\Xi_{j_1})-C^3_{ii_1i_2i_3kk_1ll_1jj_1}\mathbf{1},
\\&M\includegraphics[height=0.7cm]{08.eps}=\includegraphics[height=0.7cm]{08.eps}-C^3_{ii_1i_2i_3kk_1ll_1jj_1}\mathbf{1},\endaligned$$
\begin{equation}\aligned &M(\mathcal{I}^{ii_1}_l(\mathcal{I}^{ll_1}_k(\mathcal{I}^{l_1l_2}(\Xi_{l_2})\mathcal{I}^{kk_1}(\Xi_{k_1}))\mathcal{I}^{i_1i_2}(\Xi_{i_2}))\mathcal{I}^{jj_1}(\Xi_{j_1}))
\\=&\mathcal{I}^{ii_1}_l(\mathcal{I}^{ll_1}_k(\mathcal{I}^{l_1l_2}(\Xi_{l_2})\mathcal{I}^{kk_1}(\Xi_{k_1}))
\mathcal{I}^{i_1i_2}(\Xi_{i_2}))\mathcal{I}^{jj_1}(\Xi_{j_1})-C^4_{ii_1i_2kk_1ll_1l_2jj_1}\mathbf{1},\endaligned\end{equation}
as well as $M(\tau)=\tau$ for the remaining basis vectors in $\mathcal{F}_0$. Here we omit the tree notation for the last one since it is the same as the one including $C^3$.
We claim that for any $\tau\in\mathcal{F}_0$,
\begin{equation}\Delta^M\tau=(M\tau)\otimes \textbf{1}.\end{equation}
Since $\tau$ satisfies $M\tau=\tau-C\mathbf{1}$ for any $\tau\in\mathcal{F}_0$, it is easy to check that (2.8) holds. Here for the definitions of $\Delta^M, \mathcal{A}, \hat{M}, \hat{\Delta}^M$ we refer to [16, Section 8.3].

For $\tau=\includegraphics[height=0.7cm]{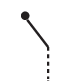}$,  we have
$$\Delta^+\includegraphics[height=0.7cm]{03.eps}=\includegraphics[height=0.7cm]{03.eps}
\otimes \mathbf{1}+\mathbf{1}\otimes \includegraphics[height=0.7cm]{03.eps}. $$
$$(\mathcal{A}\hat{M}\mathcal{A}\otimes \hat{M})\Delta^+\includegraphics[height=0.7cm]{03.eps}=\includegraphics[height=0.7cm]{03.eps}\otimes \mathbf{1}+\mathbf{1}\otimes \includegraphics[height=0.7cm]{03.eps},$$
It follows that
$$\hat{\Delta}^M\includegraphics[height=0.7cm]{03.eps}=\includegraphics[height=0.7cm]{03.eps}\otimes \mathbf{1}.$$
For $\tau=\mathcal{I}^{ii_1}_l(\tau_1),$ where $\tau_1=\includegraphics[height=0.7cm]{06.eps}$, $i,i_1=1,2,3$,  we have
$$\Delta^+\mathcal{I}^{ii_1}_l(\tau_1)=\mathcal{I}^{ii_1}_l(\tau_1)\otimes \mathbf{1}+\mathbf{1}\otimes \mathcal{I}^{ii_1}_l(\tau_1). $$
$$(\mathcal{A}\hat{M}\mathcal{A}\otimes \hat{M})\Delta^+\mathcal{I}^{ii_1}_l(\tau_1)=\mathcal{I}^{ii_1}_l(\tau_1)\otimes \mathbf{1}+\mathbf{1}\otimes \mathcal{I}^{ii_1}_l(\tau_1),$$
which implies that
$$\hat{\Delta}^M\mathcal{I}^{ii_1}_l(\tau_1)=\mathcal{I}^{ii_1}_l(\tau_1)\otimes \mathbf{1}.$$
As a consequence of this expression,  $M$ belongs to the renormalisation group $\mathfrak{R}_0$  defined in [16, Definition 8.43]. Then similar as in [16, Theorem 8.46] we can define $(\Pi^M,\Gamma^M,\Sigma^M)$ and it is an admissible model for $\mathfrak{T}_F$ on $\langle\mathcal{F}_0\rangle$. Furthermore, it extends uniquely to an admissible model for all of $\mathfrak{T}_F$.
By (2.8) we also have
$$\Pi_x^M\tau=\Pi_xM\tau.$$

Now we  lift the equation onto the abstract regularity structure.
First, we define for any $\alpha_0<0$ and compact set $\mathfrak{R}$ the norm
$$|\xi|_{\alpha_0;\mathfrak{R}}=\sup_{s\in\mathbb{R}}\|\xi1_{t\geq s}\|_{\alpha_0;\mathfrak{R}},$$
and we denote by $\bar{\mathcal{C}}_\mathfrak{s}^{\alpha_0}$ the intersections of the completions of smooth functions under $|\cdot|_{\alpha_0;\mathfrak{R}}$ for all compact sets $\mathfrak{R}$.

Since $\alpha<-\frac{5}{2}$,  Theorem 2.5 does not apply to $\mathbf{R}^+\Xi_i$ directly, where $\mathbf{R}^+:\mathbb{R}\times \mathbb{R}^d\rightarrow\mathbb{R}$ is given by $\mathbf{R}^+(t,x)=1$ for $t>0$ and $\mathbf{R}^+(t,x)=0$ otherwise.  To define the reconstruction operator  for $\mathbf{R}^+\Xi_i$ by hand,  we need the following results, which can be proved by similar arguments as in [16, Proposition 9.5] and using Lemma 3.6 below.

\vskip.10in
\th{Proposition 2.11} Let $\xi=(\xi^1,\xi^2,\xi^3)$, with $\xi^i, i=1,2,3$ being independent  white noises on $\mathbb{R}\times \mathbb{T}^3$, which we extend periodically to $\mathbb{R}^4$. Let $\rho:\mathbb{R}^4\rightarrow\mathbb{R}$ be a smooth compactly supported function with Lebesgue integral equal to $1$, set $\rho_\varepsilon(t,x)=\varepsilon^{-5}\rho(\frac{t}{\varepsilon^2},\frac{x}{\varepsilon})$ and define $\xi_\varepsilon^i=\rho_\varepsilon*\xi^i$. Then for every $i,i_1=1,2,3$,  $K^{ii_1}*\xi^{i_1}\in C(\mathbb{R},\mathcal{C}^{\alpha+2}(\mathbb{R}^3)) $ almost surely and  for every $i,i_0,i_1,j,j_1=1,2,3$,  $D_jK^{i_0i}*(K^{ii_1}*\xi^{i_1}\diamond K^{jj_1}*\xi^{j_1}):=\lim_{\varepsilon\rightarrow0}D_jK^{i_0i}*(K^{ii_1}*\xi_\varepsilon^{i_1} K^{jj_1}*\xi_\varepsilon^{j_1})$ in $C(\mathbb{R},\mathcal{C}^{2\alpha+5}(\mathbb{R}^3)) $ almost surely. Moreover,  for every compact set $\mathfrak{R}\subset \mathbb{R}^4$ and every $0<\theta<-\alpha-\frac{5}{2}$
we have
$$E|\xi^i-\xi_\varepsilon^i|_{\alpha;\mathfrak{R}}\lesssim \varepsilon^\theta.$$
 Finally for every $0<\kappa<-\alpha-\frac{5}{2}$, we have the bound
$$E\sup_{t\in[0,1]}\|K^{ii_1}*\xi^{i_1}(t,\cdot)- K^{ii_1}*\xi^{i_1}_\varepsilon(t,\cdot)\|_{\alpha+2}\lesssim\varepsilon^\kappa.$$
$$E\sup_{t\in[0,1]}\|D_jK^{i_0i}*(K^{ii_1}*\xi^{i_1}\diamond K^{jj_1}*\xi^{j_1})- D_jK^{i_0i}*(K^{ii_1}*\xi_\varepsilon^{i_1} K^{jj_1}*\xi_\varepsilon^{j_1})\|_{2\alpha+5}\lesssim\varepsilon^\kappa.$$

Now we reformulate the fixed point map as
\begin{equation}\aligned v_1^i=&\sum_{i_1=1}^3(\bar{\mathcal{P}}_{\gamma}^{ii_1}+R_{0,\gamma}^{ii_1}\mathcal{R})(\mathcal{K}_{\bar{\gamma}}+R_\gamma\mathcal{R})\emph{\textbf{R}}^+\Xi_{i_1},
\\u^i=&-\frac{1}{2}\sum_{i_1,j=1}^3(\bar{\mathcal{P}}_{\gamma}^{ii_1}+R_{0,\gamma}^{ii_1}\mathcal{R})((\mathcal{D}_j\mathcal{K})_{\bar{\gamma}}+(D_jR)_\gamma\mathcal{R})\emph{\textbf{R}}^+ (u^{i_1}\star u^j)+v^i_1+\sum_{i_1=1}^3\mathcal{G}^{ii_1}u_0^{i_1}
.\endaligned\end{equation}
Here  for $j=1,2,3,$ $\mathcal{K}_{\bar{\gamma}}$ and $(\mathcal{D}_j\mathcal{K})_{\bar{\gamma}}$ are the continuous linear operators obtained by Theorem 2.6 associated with the kernel $K$ and $D_jK$ respectively, while for $i,i_1=1,2,3,$ $\bar{\mathcal{P}}_{{\gamma}}^{ii_1}$ is the continuous linear operators obtained by Theorem 2.7 associated with the kernel $\bar{P}^{ij}$, for $f:(0,T]\rightarrow\mathcal{C}^\alpha$
 $$(R_\gamma f)_t(x):=\sum_{|k|_\mathfrak{s}<\gamma}\frac{X^k}{k!}\int \langle D_1^kR_{t-s}(x-\cdot),f_s\rangle ds,$$
$$(R_{0,\gamma}^{ii_1}f)_t(x)=\sum_{|k|_\mathfrak{s}<\gamma}\frac{X^k}{k!}\langle D_1^k(R_0^{ii_1})(x-\cdot),f_t\rangle d\bar{z},$$
 $$\mathcal{G}u_0(z)=\sum_{|k|_\mathfrak{s}<\gamma}\frac{X^k}{k!} D^k(P*Gu_0)(z),$$
where $\gamma,\bar{\gamma}$ will be chosen below. We also use that $\int K(x-y) D_jf(y)dy=\int D_jK(x-y)f(y)dy$  and  define $\mathcal{R}\emph{\textbf{R}}^+\Xi$ as the distribution $\xi\textbf{1}_{t\geq0}$.

We consider the second equation in (2.9): Define
$$V^i:=\oplus_{i_1,j=1}^3\mathcal{I}^{ii_1}_j(\mathcal{H}_F^{i_1j})\oplus \textrm{span}\{\mathcal{I}^{ii_1}(\Xi_{i_1}), i_1=1,2,3\} \oplus\bar{T},$$
$$V=V^1\times V^2\times V^3.$$

For $\gamma>0, \eta\in\mathbb{R}$ we also define $$\mathcal{D}^{\gamma,\eta}(V):=\mathcal{D}^{\gamma,\eta}(V^1)\times \mathcal{D}^{\gamma,\eta}(V^2)\times \mathcal{D}^{\gamma,\eta}(V^3).$$
$$(\mathcal{D}^{\gamma,\eta})^3:=\mathcal{D}^{\gamma,\eta}\times\mathcal{D}^{\gamma,\eta}\times\mathcal{D}^{\gamma,\eta}.$$
\vskip.10in
\th{Lemma 2.12} For $\gamma>|\alpha+2|$ and $-1<\eta\leq\alpha+2$, the map $u\mapsto u^iu^j$ is locally Lipschitz continuous from $\mathcal{D}^{\gamma,\eta}(V)$ into
$\mathcal{D}^{\gamma+\alpha+2,2\eta}$.

\proof This is a consequence of [16, Proposition 6.12, Proposition 6.15].$\hfill\Box$
\vskip.10in
Now for $\gamma,\eta$ as in Lemma 2.12 and  $u_0^{i_1}\in \mathcal{C}^\eta(\mathbb{R}^3), i_1=1,2,3$, periodic, we have $P^{ii_1}u_0^{i_1}\in \mathcal{C}^\eta(\mathbb{R}^3), i, i_1=1,2,3$ (see Lemma 3.6), which by [16, Lemma 7.5] implies that $\mathcal{G}^{ii_1}u_0^{i_1}\in \mathcal{D}^{\gamma,\eta}, i,i_1=1,2,3$. By Proposition 2.11 and [16, Remark 6.17] we also have that $v_1^i\in \mathcal{D}^{\gamma,\eta}$  for $i=1,2,3$. Now we can apply a fixed point argument in $(\mathcal{D}^{\gamma,\eta})^3$ to obtain existence and uniqueness of local solutions to (2.9).
\vskip.10in
\th{Proposition 2.13} Let $\mathfrak{T}_F$ be the regularity structure from Theorem 2.8 associated to the stochastic N-S equation driven by space-time white noise with $\alpha\in (-\frac{13}{5},-\frac{5}{2})$. Let $\eta\in (-1,\alpha+2],|\alpha+2|<\gamma<\eta-\alpha$, $u_0\in \mathcal{C}^\eta(\mathbb{R}^3)$, periodic and let $Z=(\Pi,\Gamma,\Sigma)\in\mathcal{M}_F$ be an admissible model for $\mathfrak{T}_F$  with the additional properties that for $i,i_0,i_1,j,j_1=1,2,3$, $\xi^i:=\mathcal{R}\Xi^i$ belongs to $\bar{\mathcal{C}}_\mathfrak{s}^\alpha$ and that $K^{ii_1}*\xi^{i_1}\in C(\mathbb{R},\mathcal{C}^{\eta})$, $D_jK^{i_0i}*(K^{ii_1}*\xi^{i_1}\diamond K^{jj_1}*\xi^{j_1})\in C(\mathbb{R},\mathcal{C}^{2\alpha+5}) $. Then there exists a maximal solution $\mathcal{S}^L\in(\mathcal{D}^{\gamma,\eta})^3$ to the equation (2.9).

\proof Consider the second equation in (2.9). We have that $u$ takes values in a sector of regularity $\zeta=\alpha+2$ and $u^iu^j, i,j=1,2,3,$ takes value in a sector of regularity $\bar{\zeta}=2\alpha+4$ satisfying $\zeta<\bar{\zeta}+1$. For $\eta$ and $\gamma$ we have $\bar{\eta}=2\eta$ and $\gamma>\bar{\gamma}=\gamma+\alpha+2>0$, $\bar{\gamma}<\eta+2$ and  $\bar{\gamma}+1>\gamma$. By Lemma 2.12 for $i,j=1,2,3$, $u^iu^j$ is locally Lipschitz continuous from $\mathcal{D}^{\gamma,\eta}(V)$ to $\mathcal{D}^{\bar{\gamma},\bar{\eta}}$.
  Then $\eta<(\bar{\eta}\wedge\bar{\zeta})+1$ and $(\bar{\eta}\wedge\bar{\zeta})+2>0$ are satisfied by our assumptions. We consider a fixed model.
  Denote by $M_F^i(u)$ the right hand side of the second equation in (2.9). By  Theorems 2.6, 2.7, [16, Theorem 7.1, Lemma 7.3] and local Lipschitz continuity of $u\mapsto u^iu^j$ we obtain that there exists $\kappa>0$ such that for every $R>0$
  $$\aligned\sum_{i=1}^3\interleave M_F^i(u)-M_F^i(\bar{u})\interleave_{\gamma,\eta;T}\lesssim & T^\kappa\sum_{i,j=1}^3\interleave u^iu^j-\bar{u}^i\bar{u}^j
  \interleave_{\bar{\gamma},\bar{\eta};T}\\\lesssim& T^\kappa\sum_{i=1}^3\interleave u^i-
  \bar{u}^i\interleave_{\gamma,\eta;T},\endaligned$$
  uniformly over $T\in[0,1]$ and over all $u,\bar{u}$ such that $\interleave u^i\interleave_{\gamma,\eta;T}+
  \interleave\bar{u}^i\interleave_{\gamma,\eta;T}\leq R$.
Then we  obtain local existence and uniqueness of the solutions by similar arguments as in the proof of [16, Theorem 7.8]. Here we consider  vector valued  solutions and the corresponding norm is the sum of the norm for each component.
To extend this local map up to the first time where $\sum_{i=1}^3\|(\mathcal{R}u^i)(t,\cdot)\|_\eta$ blows up,  we  write $u=v_1+v_2+v_3$ with $v_1$ in (2.9) and
$$\aligned
v_2^i=&\mathcal{T}^{ii_1}_j(v_1^{i_1}\star v_1^j),
\\v_3^i=&\mathcal{T}^{ii_1}_j[ (v_3^{i_1}+v_2^{i_1})\star (v_3^j+v_2^j)\\&+(v_3^{i_1}+v_2^{i_1})\star v_1^j+v_1^{i_1}\star (v_3^j+v_2^j)]+\sum_{i_1=1}^3\mathcal{G}^{ii_1}u_0^{i_1}
,\endaligned$$
where $\mathcal{T}^{ii_1}_j=-\frac{1}{2}\sum_{i_1,j=1}^3(\bar{\mathcal{P}}_{\gamma}^{ii_1}
+R_{0,\gamma}^{ii_1}\mathcal{R})((\mathcal{D}_j\mathcal{K})_{\bar{\gamma}}+(D_jR)_\gamma\mathcal{R})$.
In this case $v_3^i$ takes values in a function-like sector of regularity $3\alpha+8$ and we can use similar arguments as in the proof of [16, Proposition 7.11] to conclude the results.  $\hfill\Box$
\vskip.10in
\th{Remark 2.14} Here the lower bound for $\eta$ is $-1$, which seems to be optimal  by the theory of regularity structures. The reason for this is as follows: the nonlinear term always contains $v\star v$ and thus  $\bar{\eta}\leq 2\eta$ which should be larger than $-2$ required by [16, Theorem 7.8]. As a result, $\eta>-1$.
\vskip.10in

Set $O:=[-1,2]\times \mathbb{R}^3$. Given a  model $Z=(\Pi,\Gamma,\Sigma)$ for $\mathfrak{T}_F$, a periodic initial condition $u_0\in (\mathcal{C}^\eta)^3$, and some cut-off value $L>0$, we denote by $u=\mathcal{S}^L(u_0,Z)\in(\mathcal{D}^{\gamma,\eta})^3$ and $T=T^L(u_0,Z)\in\mathbb{R}_+\cup\{+\infty\}$ the (unique) modelled distribution and time such that (2.9) holds on $[0,T]$, such that $\|(\mathcal{R}u)(t,\cdot)\|_\eta<L$ for $t<T$, and such that $\|(\mathcal{R}u)(t,\cdot)\|_\eta\geq L$ for $t\geq T$. Then by [16, Corollary 7.12] we obtain the following result.
\vskip.10in

\th{Proposition 2.15} Let $L>0$ be fixed. In the setting of Proposition 2.13, for every $\varepsilon>0$ and $C>0$ there exists $\delta>0$ such that setting $T=1\wedge T^L(u_0,Z)\wedge T^L(\bar{u}_0,\bar{Z})$ we have
$$\|\mathcal{S}^L(u_0,Z)-\mathcal{S}^L(\bar{u}_0,\bar{Z})\|_{\gamma,\eta;T}\leq \varepsilon,$$
for all $u_0,\bar{u}_0, Z,\bar{Z}$ provided that $\interleave Z\interleave_{\gamma;T}\leq C,  \interleave\bar{Z}\interleave_{\gamma;T}\leq C, \|u_0\|_\eta\leq L/2,
\|\bar{u}_0\|_\eta\leq L/2, \|u_0-\bar{u}_0\|_\eta\leq \delta,$ and $\interleave Z;\bar{Z}\interleave_{\gamma;T}\leq \delta$ and
$$|\xi|_{\alpha;O}+|\bar{\xi}|_{\alpha;O}\leq C,\quad\sum_{i,i_1=1}^3\sup_{t\in[0,1]}\big[\|(K^{ii_1}*\xi^{i_1})(t,\cdot)\|_{\eta}+\|(K^{ii_1}*\bar{\xi}^{i_1})(t,\cdot)\|_{\eta}\big]\leq C,$$
$$\sup_{t\in[0,1]}\sum_{i,i_0,i_1,j,j_1=1}^3\big[\|D_jK^{i_0i}*(K^{ii_1}*\xi^{i_1}\diamond K^{jj_1}*\xi^{j_1})\|_\eta+ \|D_jK^{i_0i}*(K^{ii_1}*\bar{\xi}^{i_1} K^{jj_1}*\bar{\xi}^{j_1})\|_{\eta}\big]\leq C.$$
as well as
$$|\xi-\bar{\xi}|_{\alpha;O}\leq \delta,\quad\sum_{i,i_1=1}^3\sup_{t\in[0,1]}\|(K^{ii_1}*\xi^{i_1})(t,\cdot)-(K^{ii_1}*\bar{\xi}^{i_1})(t,\cdot)\|_{\eta}\leq \delta,$$
$$\sup_{t\in[0,1]}\sum_{i,i_0,i_1,j,j_1=1}^3\|D_j K^{i_0i}*(K^{ii_1}*\xi^{i_1}\diamond K^{jj_1}*\xi^{j_1})- D_j K^{i_0i}*(K^{ii_1}*\bar{\xi}^{i_1} K^{jj_1}*\bar{\xi}^{j_1})\|_{\eta}\leq \delta.$$
where $\bar{\xi}^i=\bar{\mathcal{R}}\Xi^i$  and $\bar{\mathcal{R}}$ is the reconstruction operator associated to $\bar{Z}$.
\vskip.10in

As in [16, Section 9] we now identify solutions corresponding to a model that has been renormalised by $M$ with classical solutions to a modified equation.
\vskip.10in
\th{Proposition 2.16} Given a continuous periodic vector $\xi_\varepsilon=(\xi_\varepsilon^1,\xi_\varepsilon^2,\xi_\varepsilon^3)$, denote by $Z_\varepsilon=(\Pi^{(\varepsilon)},\Gamma^{(\varepsilon)},\Sigma^{(\varepsilon)})$ the associated canonical model realising
$\mathfrak{T}_F$ given in Proposition 2.9. Let $M$ be the renormalisation map defined in (2.7). Then for every $L>0$ and periodic $u_0\in C^\eta(\mathbb{R}^3;\mathbb{R}^3)$, $u_\varepsilon=\mathcal{R}S^L(u_0,Z_\varepsilon)$ satisfies the following equation on $[0,T^L(u_0,Z_\varepsilon)]$ in the mild sense:
$$\partial_tu_\varepsilon=\Delta u_\varepsilon-\frac{1}{2}P\textrm{\textrm{div}}(u_\varepsilon\otimes u_\varepsilon)+P\xi_\varepsilon,\quad \textrm{div} u_\varepsilon=0,\quad  u_\varepsilon(0)=Pu_0.$$
Furthermore,
$u_\varepsilon^M=\mathcal{R}S^L(u_0,MZ_\varepsilon)$ also satisfies the same equation on $[0,T^L(u_0,MZ_\varepsilon)]$ in the mild sense.

\proof We follow a similar argument as in the proof of [16, Proposition 9.4].

 For $i=1,2,3$,  the solution $u^i$ to the abstract fixed point map can be expanded as
$$\aligned u^i=&\sum_{i_1=1}^3\mathcal{I}^{ii_1}(\Xi_{i_1})-\frac{1}{2}\sum_{j,i_1,i_2,j_1=1}^3\mathcal{I}^{ii_1}_j(\mathcal{I}^{i_1i_2}(\Xi_{i_2})\mathcal{I}^{jj_1}(\Xi_{j_1}))
+\varphi^i\mathbf{1}-\frac{1}{2}\sum_{j,i_1,j_1=1}^3\mathcal{I}^{ii_1}_j(\mathcal{I}^{jj_1}(\Xi_{j_1}))\varphi^{i_1}
\\&-\frac{1}{2}\sum_{j,i_1,i_2=1}^3\mathcal{I}^{ii_1}_j(\mathcal{I}^{i_1i_2}(\Xi_{i_2}))\varphi^j
+\frac{1}{4}\sum_{i_1,i_2,i_3,j,j_1,k,k_1=1}^3\mathcal{I}^{ii_1}_k(\mathcal{I}^{i_1i_2}_j(\mathcal{I}^{i_2i_3}(\Xi_{i_3})\mathcal{I}^{jj_1}(\Xi_{j_1}))\mathcal{I}^{kk_1}(\Xi_{k_1}))
\\&+\frac{1}{4}\sum_{i_1,i_2,j,j_1,k,k_1,k_2=1}^3\mathcal{I}^{ii_1}_k(\mathcal{I}^{i_1i_2}(\Xi_{i_2})\mathcal{I}^{kk_1}_j
(\mathcal{I}^{k_1k_2}(\Xi_{k_2})\mathcal{I}^{jj_1}(\Xi_{j_1})))+\rho_u.\endaligned$$
i.e.
$$\aligned u^i=&\includegraphics[height=0.7cm]{01.eps}-\frac{1}{2}\includegraphics[height=0.7cm]{04.eps}
+\varphi^i\mathbf{1}-\frac{1}{2}\sum_{i_1=1}^3\includegraphics[height=0.7cm]{03.eps}\varphi^{i_1}\\&-\frac{1}{2}\sum_{j=1}^3\includegraphics[height=0.7cm]{03.eps}\varphi^j
+\frac{1}{4}\includegraphics[height=0.7cm]{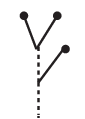}+\frac{1}{4}\includegraphics[height=0.7cm]{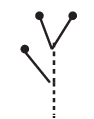}+\rho_u.\endaligned$$
Here every component of $\rho_u$ has homogeneity strictly greater than $3\alpha+8$.
Then we have
$$\aligned u^iu^j=&\frac{1}{4}\sum_{i_1,i_2,j_1,j_2,k,k_1,l,l_1=1}^3\mathcal{I}^{ii_1}_k(\mathcal{I}^{i_1i_2}(\Xi_{i_2})\mathcal{I}^{kk_1}(\Xi_{k_1}))
\mathcal{I}^{jj_1}_l(\mathcal{I}^{j_1j_2}(\Xi_{j_2})\mathcal{I}^{ll_1}(\Xi_{l_1}))\\&-
\frac{1}{2}\sum_{i_1,i_2,k,k_1=1}^3\mathcal{I}^{ii_1}_k(\mathcal{I}^{i_1i_2}(\Xi_{i_2})\mathcal{I}^{kk_1}(\Xi_{k_1}))\varphi^j-
\frac{1}{2}\varphi^i\sum_{j_1,j_2,k,k_1=1}^3\mathcal{I}^{jj_1}_k(\mathcal{I}^{j_1j_2}(\Xi_{j_2})\mathcal{I}^{kk_1}(\Xi_{k_1}))\\ &+\varphi^i\varphi^j-
\frac{1}{2}\sum_{i_1,i_2,j_1,k,k_1=1}^3\mathcal{I}^{ii_1}_k(\mathcal{I}^{i_1i_2}(\Xi_{i_2})\mathcal{I}^{kk_1}(\Xi_{k_1}))\mathcal{I}^{jj_1}(\Xi_{j_1})
+\varphi^i\sum_{j_1=1}^3\mathcal{I}^{jj_1}(\Xi_{j_1})
\\&-\frac{1}{2}\sum_{i_1,j_1,k,k_1=1}^3\mathcal{I}^{ii_1}_k(\mathcal{I}^{kk_1}(\Xi_{k_1}))\varphi^{i_1}\mathcal{I}^{jj_1}(\Xi_{j_1})
-\frac{1}{2}\sum_{i_1,i_2,j_1,k=1}^3\mathcal{I}^{ii_1}_k(\mathcal{I}^{i_1i_2}(\Xi_{i_2}))\varphi^k\mathcal{I}^{jj_1}(\Xi_{j_1})
\\&+\frac{1}{4}\sum_{i_1,i_2,i_3,l,l_1,k,k_1,j_1=1}^3
\mathcal{I}^{ii_1}_k(\mathcal{I}^{i_1i_2}_l(\mathcal{I}^{i_2i_3}(\Xi_{i_3})
\mathcal{I}^{ll_1}(\Xi_{l_1}))\mathcal{I}^{kk_1}(\Xi_{k_1}))
\mathcal{I}^{jj_1}(\Xi_{j_1}) \\&+\frac{1}{4}\sum_{i_1,i_2,k,k_1,k_2,l,l_1,j_1=1}^3\mathcal{I}^{ii_1}_k(\mathcal{I}^{i_1i_2}(\Xi_{i_2})\mathcal{I}^{kk_1}_l(\mathcal{I}^{k_1k_2}(\Xi_{k_2})
\mathcal{I}^{ll_1}(\Xi_{l_1})))
\mathcal{I}^{jj_1}(\Xi_{j_1})
\\ &-\frac{1}{2}\sum_{i_1,j_1,j_2,k,k_1=1}^3\mathcal{I}^{jj_1}_k(\mathcal{I}^{j_1j_2}(\Xi_{j_2})\mathcal{I}^{kk_1}(\Xi_{k_1}))\mathcal{I}^{ii_1}(\Xi_{i_1})
+\sum_{i_1=1}^3\varphi^j\mathcal{I}^{ii_1}(\Xi_{i_1})
\\ &-\frac{1}{2}\sum_{i_1,j_1,k,k_1=1}^3\mathcal{I}^{jj_1}_k(\mathcal{I}^{kk_1}(\Xi_{k_1}))\varphi^{j_1}\mathcal{I}^{ii_1}(\Xi_{i_1})
-\frac{1}{2}\sum_{i_1,j_1,j_2,k=1}^3\mathcal{I}^{jj_1}_k(\mathcal{I}^{j_1j_2}(\Xi_{j_2}))\varphi^k\mathcal{I}^{ii_1}(\Xi_{i_1})
\endaligned$$
$$\aligned&+\frac{1}{4}\sum_{i_1,j_1,j_2,j_3,l,l_1,k,k_1=1}^3\mathcal{I}^{jj_1}_k(\mathcal{I}^{j_1j_2}_l(\mathcal{I}^{j_2j_3}(\Xi_{j_3})
\mathcal{I}^{ll_1}(\Xi_{l_1}))\mathcal{I}^{kk_1}(\Xi_{k_1}))
\mathcal{I}^{ii_1}(\Xi_{i_1})
\\&+\frac{1}{4}\sum_{i_1,j_1,j_2,l,l_1,k,k_1,k_2=1}^3\mathcal{I}^{jj_1}_k(\mathcal{I}^{j_1j_2}(\Xi_{j_2})\mathcal{I}^{kk_1}_l(\mathcal{I}^{k_1k_2}
(\Xi_{k_2})\mathcal{I}^{ll_1}(\Xi_{l_1})))
\mathcal{I}^{ii_1}(\Xi_{i_1})
\\&+\sum_{i_1,j_1=1}^3\mathcal{I}^{ii_1}(\Xi_{i_1})\mathcal{I}^{jj_1}(\Xi_{j_1})+\rho_F,\endaligned$$
i.e. $$\aligned u^iu^j=&\frac{1}{4}\includegraphics[height=0.7cm]{07.eps}-
\frac{1}{2}\includegraphics[height=0.7cm]{04.eps}\varphi^j-
\frac{1}{2}\varphi^i\includegraphics[height=0.7cm]{04.eps}\\&+\varphi^i\varphi^j-
\frac{1}{2}\includegraphics[height=0.7cm]{06.eps}
+\varphi^i\includegraphics[height=0.5cm]{01.eps}
-\frac{1}{2}\sum_{i_1=1}^3\includegraphics[height=0.7cm]{05.eps}\varphi^{i_1}
-\frac{1}{2}\sum_{k=1}^3\includegraphics[height=0.7cm]{05.eps}\varphi^k
\\&+\frac{1}{4}\includegraphics[height=0.7cm]{08.eps}+\frac{1}{4}\includegraphics[height=0.7cm]{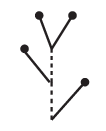}
-\frac{1}{2}\includegraphics[height=0.7cm]{06.eps}
+\varphi^j\includegraphics[height=0.5cm]{01.eps}
-\frac{1}{2}\sum_{j_1=1}^3\includegraphics[height=0.7cm]{05.eps}\varphi^{j_1}
-\frac{1}{2}\sum_{k=1}^3\includegraphics[height=0.7cm]{05.eps}\varphi^{k}
\\&+\frac{1}{4}\includegraphics[height=0.7cm]{08.eps}+\frac{1}{4}\includegraphics[height=0.7cm]{09.eps}
+\includegraphics[height=0.5cm]{02.eps}+\rho_F,\endaligned$$
where $\rho_F$ has strictly positive homogeneity. Moreover, we have
$$\mathcal{R}u^i=-\frac{1}{2}\sum_{i_1,i_2,j,j_1=1}^3D_jK^{ii_1}*(K^{i_1i_2}*\xi^{i_2}_\varepsilon\cdot K^{jj_1}*\xi^{j_1}_\varepsilon)+\varphi^i+\sum_{i_1=1}^3K^{ii_1}*\xi^{i_1}_\varepsilon,$$
where $\mathcal{R}$ is the reconstruction operator associated with $Z_\varepsilon$.
Since $\Delta^M\tau=M\tau\otimes 1$, one has the identity $(\Pi_z^{M,(\varepsilon)}\tau)(z)=(\Pi_z^{(\varepsilon)}M\tau)(z)$. It follows that for the reconstruction operator $\mathcal{R}^M$ associated with $MZ_\varepsilon$
$$\aligned\mathcal{R}^M(u^iu^j)=&\mathcal{R}u^i\mathcal{R}u^j
-\frac{1}{4}\sum_{i_1,i_2,j_1,j_2,k,k_1,l,l_1=1}^3C^2_{ii_1i_2jj_1j_2kk_1ll_1}-\sum_{i_1,j_1=1}^3C^1_{ii_1jj_1}\\&-\frac{1}{4}\sum_{i_1,i_2,i_3,k,k_1,l,l_1,j_1=1}^3C^3_{ii_1i_2i_3ll_1kk_1jj_1}
-\frac{1}{4}\sum_{i_1,i_2,k,k_1,k_2,l,l_1,j_1=1}^3C^4_{ii_1i_2ll_1kk_1k_2jj_1}\\&-\frac{1}{4}
\sum_{i_1,k,k_1,l,l_1,j_1,j_2,j_3=1}^3C^3_{jj_1j_2j_3ll_1kk_1ii_1}
-\frac{1}{4}\sum_{i_1,k,k_1,k_2,l,l_1,j_1,j_2=1}^3C^4_{jj_1j_2ll_1kk_1k_2ii_1},\endaligned$$
which together with the fact that $\int_0^t\int D_jG(t-s,x-y)dyds=0$ implies  the results.$\hfill\Box$
\vskip.10in
Now we follow [16, Section 10] to show that if $\xi_\varepsilon\rightarrow\xi$ with $Z_\varepsilon$ denoting the corresponding model, then one can find a sequence $M_\varepsilon\in\mathfrak{R}_0$ such that $M_\varepsilon Z_\varepsilon\rightarrow \hat{Z}$.
\vskip.10in
\th{Theorem 2.17} Let $\mathfrak{T}_F$ be the regularity structure associated to the stochastic N-S equation driven by space-time white noise for $\beta=2, \alpha\in(-\frac{13}{5},-\frac{5}{2})$, let $\xi_\varepsilon=\rho_\varepsilon*\xi$ be as in Proposition 2.11, $\rho_\varepsilon$ symmetric in the sense that $\rho_\varepsilon(t,x)=\rho_\varepsilon(t,-x)$, and let $Z_\varepsilon$ be the associated canonical model and $M_\varepsilon$ be a sequence of renormalisation linear maps defined in (2.7) corresponding to  $C^{1,\varepsilon}, C^{2,\varepsilon}, C^{3,\varepsilon}, C^{4,\varepsilon}$, which will be defined in the proof. Set $\hat{Z}_\varepsilon=M_\varepsilon Z_\varepsilon$.  Then, there exists a random model $\hat{Z}$ independent of the choice of the mollifier $\rho$ and $M_\varepsilon\in\mathfrak{R}_0$ such that $M_\varepsilon Z_\varepsilon\rightarrow \hat{Z}$ in probability.

More precisely, for any $\theta<-\frac{5}{2}-\alpha$, any compact set $\mathfrak{R}$ and any $\gamma<r$ we have
$$E\interleave M_\varepsilon Z_\varepsilon;\hat{Z}\interleave_{\gamma;\mathfrak{R}}\lesssim\varepsilon^\theta,$$
uniformly over $\varepsilon\in (0,1]$.

\proof By [16, Theorem 10.7] and [18, Theorem 6.1], it is sufficient to prove that for $\tau\in\mathcal{F}$ with $|\tau|_\mathfrak{s}<0$,  any test function $\varphi\in\mathcal{B}_r$ and every $x\in \mathbb{R}^3$, $t\in[0,T]$, there exist random variables $\hat{\Pi}_x^t\tau(\varphi)$ such that for $\kappa>0$ small enough
\begin{equation}E|(\hat{\Pi}_x^t\tau)(\varphi_x^\lambda)|^2\lesssim \lambda^{2|\tau|_{\mathfrak{s}}+\kappa},\end{equation}
and such that for some $0<\theta<-\frac{5}{2}-\alpha$,
\begin{equation}E|(\hat{\Pi}_x^t\tau-\hat{\Pi}_x^{(\varepsilon,t)}\tau)(\varphi_x^\lambda)|^2\lesssim \varepsilon^{2\theta}\lambda^{2|\tau|_{\mathfrak{s}}+\kappa}.\end{equation}

Since the map $\varphi\mapsto(\hat{\Pi}_x\tau)(\varphi)$ is linear, we can find some functions $\hat{\mathcal{W}}^{(\varepsilon;k)}\tau$ with $(\hat{\mathcal{W}}^{(\varepsilon;k)}\tau)(t,x)\in L^2(\mathbb{R}\times \mathbb{T}^3)^{\otimes k}$, where $(t,x)\in\mathbb{R}^4$ and such that
$$(\hat{\Pi}_0^{(\varepsilon,t)}\tau)(\varphi)=\sum_{k\leq \|\tau\|}I_k\bigg(\int\varphi(y)(\hat{\mathcal{W}}^{(\varepsilon;k)}\tau)(t,y)dy\bigg),$$
where $\|\tau\|$ denotes the number of occurrences of $\Xi$ in the expression $\tau$ and $I_k$ is defined as in [16, Section 10.1].
To obtain (2.10) and (2.11) it is sufficient to find functions $\hat{\mathcal{W}}^{(k)}\tau\in L^2(\mathbb{R}\times \mathbb{T}^3)^{\otimes k}$, define
 $$(\hat{\Pi}_x^t\tau)(\varphi):=\sum_{k\leq \|\tau\|}I_k\bigg(\int\varphi(y)S_x^{\otimes k}(\hat{\mathcal{W}}^{(k)}\tau)(t,y)dy\bigg),$$
 and  estimate the terms $|\langle(\hat{\mathcal{W}}^{(\varepsilon;k)}\tau)(t,y),(\hat{\mathcal{W}}^{(\varepsilon;k)}\tau)(t,\bar{y})\rangle|$ and $|\langle(\delta\hat{\mathcal{W}}^{(\varepsilon;k)}\tau)(t,y),(\delta\hat{\mathcal{W}}^{(\varepsilon;k)}\tau)(t,\bar{y})\rangle|$,
 where $\{S_x\}_{x\in\mathbb{R}^3}$ is the unitary operators associated with translation invariance and $\delta\hat{\mathcal{W}}^{(\varepsilon;k)}\tau=\hat{\mathcal{W}}^{(\varepsilon;k)}\tau-\hat{\mathcal{W}}^{(k)}\tau$.

For $\tau=\mathcal{I}^{ii_1}(\Xi_{i_1}), i,i_1=1,2,3, $ it is easy to conclude that (2.10), (2.11) hold in this case.

For $\tau=\mathcal{I}^{ii_1}(\Xi_{i_1})\mathcal{I}^{jj_1}(\Xi_{j_1})$, $i,i_1,j,j_1=1,2,3$,  we have
$$\hat{\Pi}_x^{(\varepsilon,t)}\tau(y)=\int K^{ii_1}(t-s,y-y_1)\xi^{i_1}_\varepsilon(s,y_1)dsdy_1\int K^{jj_1}(t-s,y-y_1)\xi^{j_1}_\varepsilon(s,y_1)dy_1-C^{1,\varepsilon}_{ii_1jj_1}.$$
If we choose $C^{1,\varepsilon}_{ii_1jj_1}:=\langle K^{ii_1}_\varepsilon,K^{jj_1}_\varepsilon\rangle$ with $K_\varepsilon=\rho_\varepsilon*K$, we have
$$\hat{\Pi}_x^{(\varepsilon,t)}\tau(y)=\int K^{ii_1}(t-s_1,y-y_1)K^{jj_1}(t-s_2,y-y_2)\xi^{i_1}_\varepsilon(s_1,y_1)\diamond \xi^{j_1}_\varepsilon(s_2,y_2)ds_1dy_1ds_2dy_2,$$
so that $\hat{\Pi}_x^{(\varepsilon,t)}\tau(y)$ belongs to the homogeneous chaos of order $2$ with
$$(\hat{\mathcal{W}}^{(\varepsilon;2)}\tau)(t,y;z_1,z_2)=K_\varepsilon^{ii_1}(t-s_1,y-y_1)K_\varepsilon^{jj_1}(t-s_2,y-y_2),$$
for $z_i=(s_i,y_i), i=1,2$.
Since for $i,j=1,2,3,$ $K^{ij}$ is of order $-3$,  applying [16, Lemma 10.14] we deduce that
$$|\langle (\hat{\mathcal{W}}^{(\varepsilon;2)}\tau)(t,y),(\hat{\mathcal{W}}^{(\varepsilon;2)}\tau)(t,\bar{y})\rangle|\lesssim |y-\bar{y}|^{-2}$$
holds uniformly over $\varepsilon\in (0,1]$.
Hence we can choose $$(\hat{\mathcal{W}}^{(2)}\tau)(y;z_1,z_2)=K^{ii_1}(t-s,y-z_1)K^{jj_1}(t-s,y-z_2),$$
and we use it to define $(\hat{\Pi}_x^t\tau)(\psi)$.
In the same way, it is straightforward to obtain an analogous bound on $(\hat{\mathcal{W}}^{(2)})(\tau)$, which implies that (2.10) holds in this case. So it remains to find similar bounds for
$(\delta\hat{\mathcal{W}}^{(\varepsilon;2)}\tau)=(\hat{\mathcal{W}}^{(\varepsilon;2)}\tau)-(\hat{\mathcal{W}}^{(2)}\tau)$.
 Similarly, by [16, Lemma 10.17] we have for $0<\kappa+\theta<-2(2\alpha+5)$
$$|\langle (\delta\hat{\mathcal{W}}^{(\varepsilon;2)}\tau)(t,y),(\delta\hat{\mathcal{W}}^{(\varepsilon;2)}\tau)(t,\bar{y})\rangle|\lesssim \varepsilon^\theta|y-\bar{y}|^{-2-\theta},$$
holds uniformly over $\varepsilon\in (0,1]$.
Then we obtain the bound
$$\aligned&|\int\int\psi^\lambda(y)\psi^\lambda(\bar{y})\langle (\delta\hat{\mathcal{W}}^{(\varepsilon;2)}\tau)(t,y),(\delta\hat{\mathcal{W}}^{(\varepsilon;2)}\tau)(t,\bar{y})\rangle dyd\bar{y}|\lesssim \varepsilon^\theta\lambda^{\kappa+2(2\alpha+4)},\endaligned$$
 which implies (2.11) holds in this case.

For $\tau=\mathcal{I}^{ii_1}_j(\mathcal{I}^{i_1i_2}(\Xi_{i_2})\mathcal{I}^{jj_1}(\Xi_{j_1}))$, $i,i_1,i_2,j,j_1=1,2,3$,  we have
$$(\hat{\mathcal{W}}^{(\varepsilon;2)}\tau)(t,y;z_1,z_2)=\int D_jK^{ii_1}(t-s_0,y-y_0)K_\varepsilon^{i_1i_2}(z_0-z_1)K_\varepsilon^{jj_1}(z_0-z_2)dz_0,$$
for $z_0=(s_0,y_0)$.
Then by [16, Lemma 10.14] we obtain that for any $\delta>0$
$$|\langle (\hat{\mathcal{W}}^{(\varepsilon;2)}\tau)(t,y),(\hat{\mathcal{W}}^{(\varepsilon;2)}\tau)(t,\bar{y})\rangle|\lesssim |y-\bar{y}|^{-\delta},$$
holds uniformly over $\varepsilon\in (0,1]$, which implies the bound
$$\aligned&|\int\int\psi^\lambda(y)\psi^\lambda(\bar{y})\langle (\hat{\mathcal{W}}^{(\varepsilon;2)}\tau)(t,y),(\hat{\mathcal{W}}^{(\varepsilon;2)}\tau)(t,\bar{y})\rangle dyd\bar{y}|\lesssim \lambda^{-6} \int_{|y|\leq\lambda,|\bar{y}|\leq\lambda}|y-\bar{y}|^{-\delta}dyd\bar{y}
 \\\lesssim&\lambda^{-3}\int_{|y|\leq 2\lambda}|y|^{-\delta}dy \lesssim\lambda^{-\delta} \lesssim\lambda^{\kappa+2(2\alpha+5)},\endaligned$$
 for $0<\kappa+\delta<-2(2\alpha+5)$.
 Hence we can choose $$(\hat{\mathcal{W}}^{(2)}\tau)(t,y;z_1,z_2)=\int D_jK^{ii_1}(t-s_0,y-y_0)K^{i_1i_2}(z_0-z_1)K^{jj_1}(z_0-z_2)dz_0,$$
 and deduce easily that (2.10) holds for $\tau=\mathcal{I}^{ii_1}_j(\mathcal{I}^{i_1i_2}(\Xi_{i_2})\mathcal{I}^{jj_1}(\Xi_{j_1}))$.
Similarly for $0<\kappa+\delta+\theta<-2(2\alpha+5)$ we have that the bound
$$\aligned&|\int\int\psi^\lambda(y)\psi^\lambda(\bar{y})\langle (\delta\hat{\mathcal{W}}^{(\varepsilon;2)}\tau)(t,y),(\delta\hat{\mathcal{W}}^{(\varepsilon;2)}\tau)(t,\bar{y})\rangle dyd\bar{y}|\lesssim \varepsilon^\theta\lambda^{\kappa+2(2\alpha+5)},\endaligned$$
holds uniformly over $\varepsilon\in (0,1]$, which also implies that (2.11) holds for $\tau=\mathcal{I}^{ii_1}_j(\mathcal{I}^{i_1i_2}(\Xi_{i_2})\mathcal{I}^{jj_1}(\Xi_{j_1}))$.

For $\tau=\mathcal{I}_j(\mathcal{I}^{i_1i_2}(\Xi_{i_2})\mathcal{I}^{jj_1}(\Xi_{j_1}))$ the same argument also implies (2.10) and (2.11) hold in this case.

In the following we use \includegraphics[height=0.3cm]{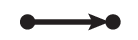} to represent a factor $K$ or $K_\varepsilon$ and \includegraphics[height=0.3cm]{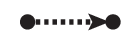} to represent $DK$ or $DK_\varepsilon$, where for simplicity we write $K^{ii_1}=K, D_jK^{ii_1}=DK$ and we do not make a difference between the graphs associated with different $K^{ii_1}$, since they have the same order. In the graphs below we also omit the dependence on $\varepsilon$ if there's no confusion. We also use the convention that if a vertex is drawn in grey, then the corresponding variable is integrated out. We also use $0$ to reprensent $(t,0)$ in the graph for simplicity

For $\tau=\mathcal{I}^{ii_1}_k(\mathcal{I}^{kk_1}(\Xi_{k_1}))\mathcal{I}^{jj_1}(\Xi_{j_1})$, $i,i_1,k,k_1,j,j_1=1,2,3$ we have
$$(\mathcal{W}^{(\varepsilon;2)}\tau)(z)= \includegraphics[height=1cm]{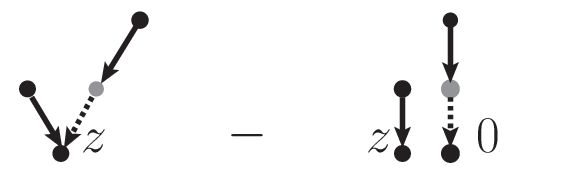}.$$
Defining kernels $Q_\varepsilon^0$, $P_\varepsilon^0$ by
$$P_\varepsilon^0(z-\bar{z})=\includegraphics[height=0.5cm]{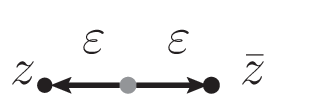},
\quad Q_\varepsilon^0(z,\bar{z})=\includegraphics[height=0.5cm]{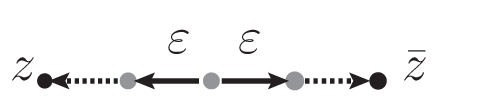},$$
for $z=(t,y)$ and $\bar{z}=(t,\bar{y})$ and we have
$$\aligned\langle\mathcal{W}^{(\varepsilon;2)}\tau(z),\mathcal{W}^{(\varepsilon;2)}\tau(\bar{z})\rangle=&P_\varepsilon^0(z-\bar{z})\delta^{(2)} Q_\varepsilon^0(z,\bar{z}),\endaligned$$
where for any function $Q$ of two variables we have set
$$\delta^{(2)}  Q(z,\bar{z})=Q(z,\bar{z})-Q(z,0)-Q(0,\bar{z})+Q(0,0).$$
It follows from [16, Lemma 10.14, Lemma 10.17] that for every $\delta>0$ we have
$$|Q_\varepsilon^0(z)-Q_\varepsilon^0(0)|\lesssim \|z\|^{1-\delta}_\mathfrak{s},\quad |P_\varepsilon^0(z)|\lesssim \|z\|^{-1}_\mathfrak{s}.$$
As a consequence we have the desired a priori bounds for $\mathcal{W}^{(\varepsilon;2)}\tau$, namely for every $\delta>0$
$$\langle (\hat{\mathcal{W}}^{(\varepsilon;2)}\tau)(z),(\hat{\mathcal{W}}^{(\varepsilon;2)}\tau)(\bar{z})\rangle\lesssim |y-\bar{y}|^{-1}(|y-\bar{y}|^{1-\delta}+|y|^{1-\delta}+|\bar{y}|^{1-\delta})$$
holds uniformly over $\varepsilon\in(0,1]$. As previously, we define $\hat{W}^{(2)}\tau$ like $\hat{\mathcal{W}}^{(\varepsilon;2)}\tau$, but with $K_\varepsilon$ replaced by $K$.
Moreover, we use $\includegraphics[height=0.2cm]{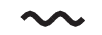}$ to represent the kernel $K-K_\varepsilon$, and we have
$$(\delta\mathcal{W}^{(\varepsilon;2)}\tau)(z)= \includegraphics[height=1cm]{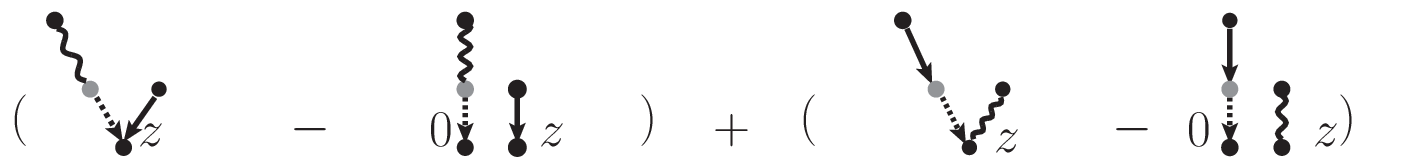}.$$
By a similar calculation as above we obtain the following bounds
$$\aligned\langle (\delta\hat{\mathcal{W}}^{(\varepsilon;2)}\tau)(z),(\delta\hat{\mathcal{W}}^{(\varepsilon;2)}\tau)(\bar{z})\rangle\lesssim& \varepsilon^{2\theta}|y-\bar{y}|^{-1}(|y-\bar{y}|^{1-2\theta-\delta}+
|y|^{1-2\theta-\delta}+|\bar{y}|^{1-2\theta-\delta})
\\&+\varepsilon^{2\theta}|y-\bar{y}|^{-1-2\theta}(|y-\bar{y}|^{1-\delta}+
|y|^{1-\delta}+|\bar{y}|^{1-\delta}),\endaligned$$
which is valid  uniformly over $\varepsilon\in (0,1]$, provided that $\theta<1, \delta>0$. Here we used  [16, Lemma 10.17]. We come to $\hat{\mathcal{W}}^{(\varepsilon;0)}\tau$ and have
$$(\hat{\mathcal{W}}^{(\varepsilon;0)}\tau)(z)=\includegraphics[height=1.1cm]{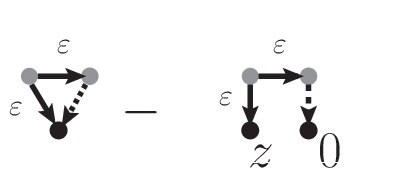}. $$
Since $K$ is symmetric and $DK$ is anti-symmetric with respect to the space variable, we conclude that
$$\includegraphics[height=0.6cm]{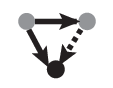}=0,$$
which implies the following
$$(\hat{\mathcal{W}}^{(\varepsilon;0)}\tau)(z)=-\includegraphics[height=1cm]{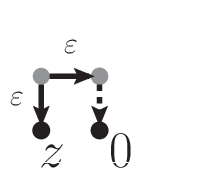}. $$
By [16, Lemma 10.14, Lemma 10.17] we have that for every $\delta>0$
$$|(\hat{\mathcal{W}}^{(\varepsilon;0)}\tau)(z)|\lesssim |y|^{-\delta},$$
holds uniformly over $\varepsilon\in(0,1]$. Similar bounds also hold for $(\delta\hat{\mathcal{W}}^{(\varepsilon;0)}\tau)$. Then we can easily conclude that (2.10) (2.11) hold for $\tau=\mathcal{I}^{ii_1}_k(\mathcal{I}^{kk_1}(\Xi_{k_1}))\mathcal{I}^{jj_1}(\Xi_{j_1})$.

For $\tau=\mathcal{I}^{ii_1}_k(\mathcal{I}^{i_1i_2}(\Xi_{i_2}))\mathcal{I}^{jj_1}(\Xi_{j_1})$, $i,i_1,i_2,k,j,j_1=1,2,3$, we can prove similar bounds as above, since in this case we also have
$$\includegraphics[height=0.6cm]{36.eps}=0.$$

For $\tau=\mathcal{I}^{ii_1}_k(\mathcal{I}^{i_1i_2}(\Xi_{i_2})\mathcal{I}^{kk_1}(\Xi_{k_1}))\mathcal{I}^{jj_1}(\Xi_{j_1})= \includegraphics[height=0.7cm]{06.eps}$, $i,i_1,i_2,k,k_1,j,j_1=1,2,3$,
we have the following identities
$$(\hat{\mathcal{W}}^{(\varepsilon;3)}\tau)(z)= \includegraphics[height=1cm]{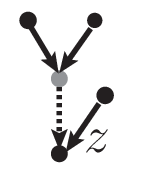},$$
$$(\hat{\mathcal{W}}^{(\varepsilon;1)}_1\tau)(z)= \includegraphics[height=1.2cm]{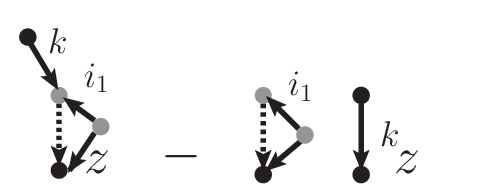},$$
$$(\hat{\mathcal{W}}^{(\varepsilon;1)}_2\tau)(z)= \includegraphics[height=1.2cm]{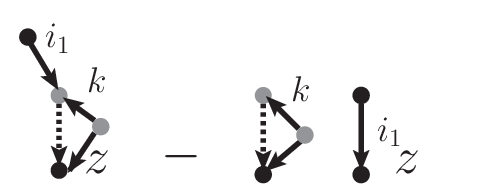}.$$
Then
$$\aligned\langle\hat{\mathcal{W}}^{(\varepsilon;3)}\tau(z),\hat{\mathcal{W}}^{(\varepsilon;3)}\tau(\bar{z})\rangle=&P_\varepsilon^0(z-\bar{z}) Q_\varepsilon(z-\bar{z}),\endaligned$$
for $z=(t,y)$ and $\bar{z}=(t,\bar{y})$, where $$Q_\varepsilon(z-\bar{z})=\includegraphics[height=0.7cm]{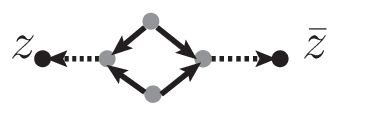}, \quad \includegraphics[height=0.7cm]{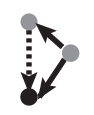}=0.$$
By [16, Lemmas 10.14 and 10.17] for every $\delta>0$ we obtain the bound
$$|Q_\varepsilon(z-\bar{z})|\lesssim |y-\bar{y}|^{-\delta},$$
which implies that
$$|\langle\hat{\mathcal{W}}^{(\varepsilon;3)}\tau(z),\hat{\mathcal{W}}^{(\varepsilon;3)}\tau(\bar{z})\rangle|\lesssim|y-\bar{y}|^{-1-\delta},$$
holds uniformly over $\varepsilon\in(0,1]$.
As previously, we define $\hat{W}^{(3)}\tau$ like $\hat{\mathcal{W}}^{(\varepsilon;3)}\tau$,  but with $K_\varepsilon$ replaced by $K$.
Then $\delta\hat{\mathcal{W}}^{(\varepsilon;3)}\tau$ can be bounded in a manner similar as before. Now for $\hat{\mathcal{W}}^{(\varepsilon;1)}\tau$, we have
$$(\hat{\mathcal{W}}^{(\varepsilon;1)}_1\tau)(z)=((\mathcal{R}_1L_\varepsilon)*K^{kk_1}_\varepsilon)(z),$$
where $L_\varepsilon(z)=\includegraphics[height=0.7cm]{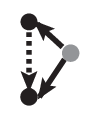}$ and $(\mathcal{R}_1L_\varepsilon)(\psi)=\int L_\varepsilon(x)(\psi(x)-\psi(0))dx$ for $\psi$ smooth with compact support. It follows from [16, Lemma 10.16] that
the bound
$$|\langle(\hat{\mathcal{W}}^{(\varepsilon;1)}_1\tau)(z),(\hat{\mathcal{W}}^{(\varepsilon;1)}_1\tau)(\bar{z})\rangle|\lesssim |y-\bar{y}|^{-1}$$
holds uniformly for $\varepsilon\in(0,1]$. Similarly, this bound also holds for $(\hat{\mathcal{W}}^{(\varepsilon;1)}_2\tau)(z)$. Again,
 $\delta\hat{\mathcal{W}}^{(\varepsilon;1)}_i\tau, i=1,2$ can be bounded in a manner similar as before. Then we can easily conclude that (2.10), (2.11) hold for $\tau=\mathcal{I}^{ii_1}_k(\mathcal{I}^{i_1i_2}(\Xi_{i_2})\mathcal{I}^{kk_1}(\Xi_{k_1}))\mathcal{I}^{jj_1}(\Xi_{j_1})$.

For $\tau=\mathcal{I}^{ii_1}_k(\mathcal{I}^{i_1i_2}(\Xi_{i_2})\mathcal{I}^{kk_1}(\Xi_{k_1}))
\mathcal{I}^{jj_1}_l(\mathcal{I}^{j_1j_2}(\Xi_{j_2})\mathcal{I}^{ll_1}(\Xi_{l_1}))= \includegraphics[height=0.7cm]{07.eps}, i,i_1,i_2,k,k_1,j,j_1,j_2,l,l_1=1,2,3,$ we have the identities
$$(\hat{\mathcal{W}}^{(\varepsilon;4)}\tau)(z)=\includegraphics[height=1cm]{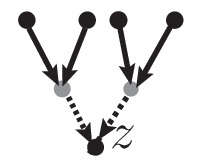}, $$
$$\langle(\hat{\mathcal{W}}^{(\varepsilon;4)}\tau)(z),(\hat{\mathcal{W}}^{(\varepsilon;4)}\tau)(\bar{z})\rangle=\includegraphics[height=1.5cm]{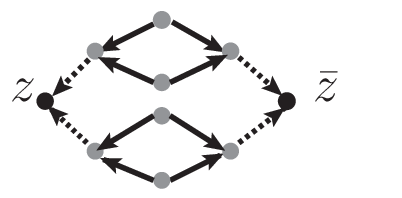},$$
for $z=(t,y), \bar{z}=(t,\bar{y})$.
Then we obtain the bound for every $\delta>0$
$$|\langle(\hat{\mathcal{W}}^{(\varepsilon;4)}\tau)(z),(\hat{\mathcal{W}}^{(\varepsilon;4)}\tau)(\bar{z})\rangle|\lesssim|y-\bar{y}|^{-\delta}.$$
Similarly, we obtain
$$|\langle(\delta\hat{\mathcal{W}}^{(\varepsilon;4)}\tau)(z),(\delta\hat{\mathcal{W}}^{(\varepsilon;4)}\tau)(\bar{z})\rangle|\lesssim
\varepsilon^{2\theta}|y-\bar{y}|^{-2\theta-\delta}$$
holds uniformly for $\varepsilon\in(0,1]$, provided $\theta<1$.

For $(\hat{\mathcal{W}}^{(\varepsilon;2)}\tau)(z), $ we have the identity
$$(\hat{\mathcal{W}}^{(\varepsilon;2)}\tau)(z)=\sum_{i=1}^4(\hat{\mathcal{W}}^{(\varepsilon;2)}_i\tau)(z).$$
$$(\hat{\mathcal{W}}^{(\varepsilon;2)}_1\tau)(z)=\includegraphics[height=1cm]{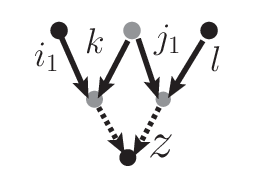}.$$
Other terms can be obtained by changing the position for $i_1,k$ or $j_1,l$. Since the estimates are similar, we omit them here. We also use the notation $\includegraphics[height=0.5cm]{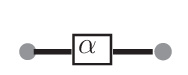}$ for $\|z-\bar{z}\|_\mathfrak{s}^{\alpha}1_{\|z-\bar{z}\|_\mathfrak{s}\leq C}$ for a constant $C$. We obtain that for $\delta>0$, $z=(t,y), \bar{z}=(t,\bar{y})$
$$\aligned&|\langle(\hat{\mathcal{W}}^{(\varepsilon;2)}_1\tau)(z),(\hat{\mathcal{W}}^{(\varepsilon;2)}_1\tau)(\bar{z})\rangle|
=\includegraphics[height=1.2cm]{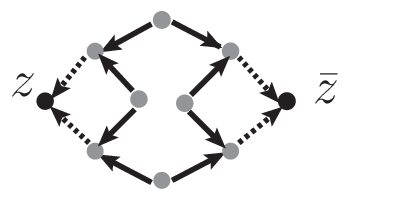}\\&\lesssim\includegraphics[height=1.2cm]{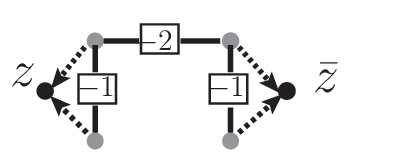}
\lesssim|y-\bar{y}|^{-\delta},\endaligned$$
holds uniformly for $\varepsilon\in(0,1]$, where we used Young's inequality in the first inequality.
Similarly, we have
$$|\langle(\delta\hat{\mathcal{W}}^{(\varepsilon;2)}_1\tau)(z),(\delta\hat{\mathcal{W}}^{(\varepsilon;2)}_1\tau)(\bar{z})\rangle|\lesssim
\varepsilon^{2\theta}|y-\bar{y}|^{-2\theta-\delta},$$
provided $\theta<1$.
Now for $\hat{\mathcal{W}}^{(\varepsilon;0)}\tau$ we have
$$(\hat{\mathcal{W}}^{(\varepsilon;0)}\tau)(z)=\includegraphics[height=1cm]{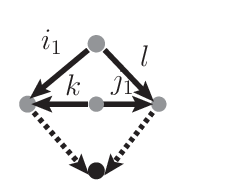}
+\includegraphics[height=1cm]{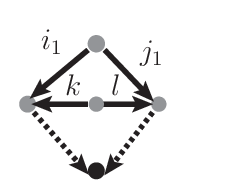}-C^{2,\varepsilon}_{ii_1i_2jj_1j_2kk_1ll_1}. $$
Hence we  choose
$$C^{2,\varepsilon}_{ii_1i_2jj_1j_2kk_1ll_1}=\includegraphics[height=1cm]{2116.eps}+\includegraphics[height=1cm]{2117.eps}$$
and also in this case (2.10), (2.11) follow.

For $\tau=\mathcal{I}^{ii_1}_l(\mathcal{I}^{i_1i_2}_k(\mathcal{I}^{i_2i_3}(\Xi_{i_3})\mathcal{I}^{kk_1}(\Xi_{k_1}))
\mathcal{I}^{ll_1}(\Xi_{l_1}))\mathcal{I}^{jj_1}(\Xi_{j_1})= \includegraphics[height=0.7cm]{08.eps}$, $i,i_1,i_2,i_3,j,j_1,k,k_1,l,l_1=1,2,3$,
we have the following identities:
$$(\hat{\mathcal{W}}^{(\varepsilon;4)}\tau)(z)= \includegraphics[height=1.5cm]{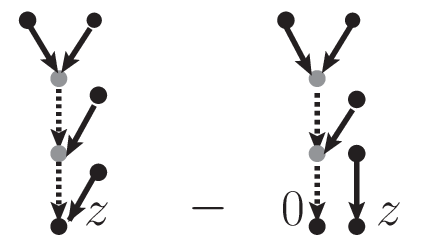}.$$
$$(\hat{\mathcal{W}}^{(\varepsilon;2)}\tau)(z)=\sum_{i=1}^5(\hat{\mathcal{W}}^{(\varepsilon;2)}_i\tau)(z)=\sum_{i=1}^5[(\hat{\mathcal{W}}^{(\varepsilon;2)}_{i1}\tau)(z)
-(\hat{\mathcal{W}}^{(\varepsilon;2)}_{i2}\tau)(z)],$$
where
$$(\hat{\mathcal{W}}^{(\varepsilon;2)}_{11}\tau)(z)= \includegraphics[height=1.5cm]{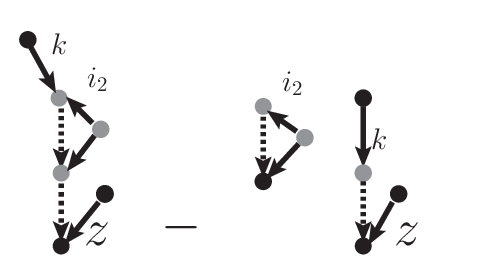}, (\hat{\mathcal{W}}^{(\varepsilon;2)}_{12}\tau)(z)= \includegraphics[height=1.5cm]{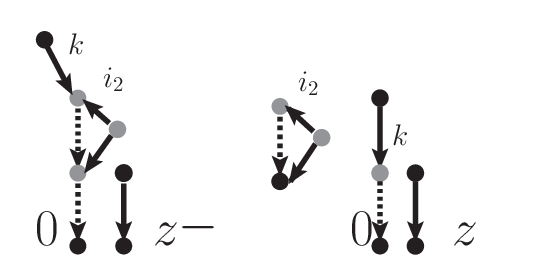},$$
$$(\hat{\mathcal{W}}^{(\varepsilon;2)}_{21}\tau)(z)= \includegraphics[height=1.5cm]{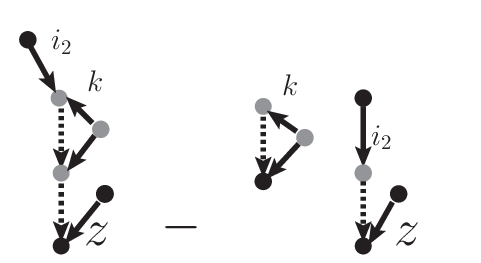}, (\hat{\mathcal{W}}^{(\varepsilon;2)}_{22}\tau)(z)= \includegraphics[height=1.5cm]{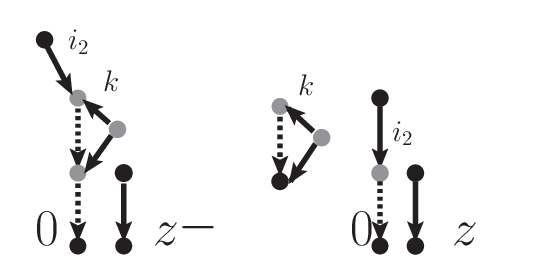},$$
$$(\hat{\mathcal{W}}^{(\varepsilon;2)}_{31}\tau)(z)= \includegraphics[height=1.5cm]{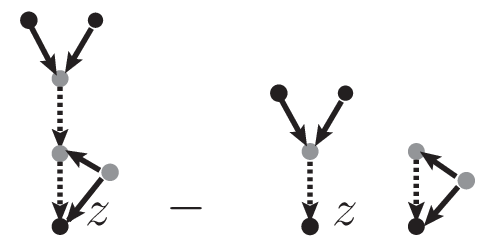}, (\hat{\mathcal{W}}^{(\varepsilon;2)}_{32}\tau)(z)= \includegraphics[height=1.5cm]{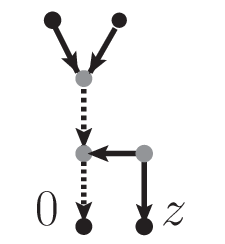},$$
$$(\hat{\mathcal{W}}^{(\varepsilon;2)}_{4}\tau)(z)=(\hat{\mathcal{W}}^{(\varepsilon;2)}_{41}\tau)(z)-(\hat{\mathcal{W}}^{(\varepsilon;2)}_{42}\tau)(z)= \includegraphics[height=1.5cm]{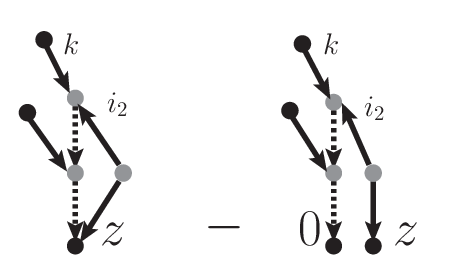},$$
$$ (\hat{\mathcal{W}}^{(\varepsilon;2)}_{5}\tau)(z)=(\hat{\mathcal{W}}^{(\varepsilon;2)}_{51}\tau)(z)-(\hat{\mathcal{W}}^{(\varepsilon;2)}_{52}\tau)(z)= \includegraphics[height=1.5cm]{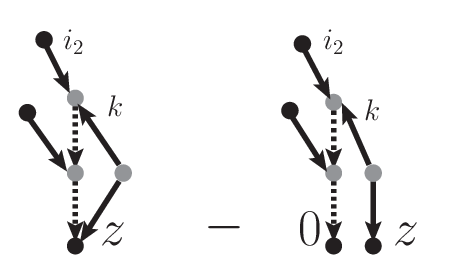}.$$
Now for $\hat{\mathcal{W}}^{(\varepsilon;4)}\tau$ we have
$$\langle\hat{\mathcal{W}}^{(\varepsilon;4)}\tau(z),\hat{\mathcal{W}}^{(\varepsilon;4)}\tau(\bar{z})\rangle=P_\varepsilon^0(z-\bar{z})\delta^{(2)} Q^2_\varepsilon(z,\bar{z}),$$
for $z=(t,y), \bar{z}=(t,\bar{y})$, where
$$Q^2_\varepsilon(z,\bar{z})=\includegraphics[height=1cm]{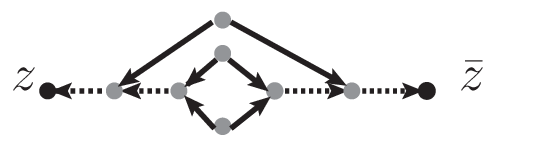},\quad \includegraphics[height=0.7cm]{2148.eps}=0.$$
By [16, Lemmas 10.14, 10.16 and 10.17] for every $\delta>0$ we have that the bound
$$|\langle\hat{\mathcal{W}}^{(\varepsilon;4)}\tau(z),\hat{\mathcal{W}}^{(\varepsilon;4)}\tau(\bar{z})\rangle|\lesssim|y-\bar{y}|^{-1}
(|y-\bar{y}|^{1-\delta}+|y|^{1-\delta}+|\bar{y}|^{1-\delta})$$
holds uniformly for $\varepsilon\in(0,1]$, and that
$$\aligned& |\langle\hat{\mathcal{W}}^{(\varepsilon;2)}_{11}\tau(z)-\hat{\mathcal{W}}^{(\varepsilon;2)}_{12}\tau(z)
,\hat{\mathcal{W}}^{(\varepsilon;2)}_{11}\tau(\bar{z})-\hat{\mathcal{W}}^{(\varepsilon;2)}_{12}\tau(\bar{z})\rangle|
\\\lesssim&|y-\bar{y}|^{-1}|\langle K^{kk_1}*\mathcal{R}_1 L_\varepsilon^1* DK^{ii_1}(z-\cdot)-K^{kk_1}*\mathcal{R}_1 L_\varepsilon^1* DK^{ii_1}((t,0)-\cdot)
,\\&K^{kk_1}*\mathcal{R}_1 L_\varepsilon^1* DK^{ii_1}(\bar{z}-\cdot)-K^{kk_1}*\mathcal{R}_1 L_\varepsilon^1* DK^{ii_1}((t,\bar{y})-\cdot)\rangle|\\\lesssim&|y-\bar{y}|^{-1}
(|y-\bar{y}|^{1-\delta}+|y|^{1-\delta}+|\bar{y}|^{1-\delta})\endaligned$$
holds uniformly for $\varepsilon\in(0,1]$, where $L^1_\varepsilon(z)=\includegraphics[height=0.7cm]{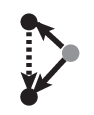}$. Then define $\hat{\mathcal{W}}^{(4)}\tau, \hat{\mathcal{W}}^{(2)}_i\tau, i=1,2,$ in a similar way as before.  Similarly, these bounds also hold for $(\hat{\mathcal{W}}^{(\varepsilon;2)}_2\tau)(z)$. Again,
$\delta\hat{\mathcal{W}}^{(\varepsilon;4)}\tau$, $\delta\hat{\mathcal{W}}^{(\varepsilon;2)}_i\tau, i=1,2$ can be bounded in a manner similar as before. For
$\hat{\mathcal{W}}^{(\varepsilon;2)}_3\tau$ we have
$$(\hat{\mathcal{W}}^{(\varepsilon;2)}_{31}\tau)(z)=((\mathcal{R}_1L^1_\varepsilon)*L^2_\varepsilon)(z),$$
where $L^1_\varepsilon(z)=\includegraphics[height=0.7cm]{2126.eps}$, $L^2_\varepsilon(z)=\includegraphics[height=1cm]{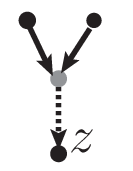}$. It follows from [16, Lemma 10.16] that for every $\delta>0$,
the bound
$$|\langle(\hat{\mathcal{W}}^{(\varepsilon;2)}_{31}\tau)(z),(\hat{\mathcal{W}}^{(\varepsilon;2)}_{31}\tau)(\bar{z})\rangle|\lesssim |y-\bar{y}|^{-\delta}$$
holds uniformly for $\varepsilon\in(0,1]$. Moreover, for $\hat{\mathcal{W}}^{(\varepsilon;2)}_{32}\tau$ we have for every $\delta\in(0,1)$
$$\aligned&|\langle(\hat{\mathcal{W}}^{(\varepsilon;2)}_{32}\tau)(z),(\hat{\mathcal{W}}^{(\varepsilon;2)}_{32}\tau)(\bar{z})\rangle|=
|\includegraphics[height=1cm]{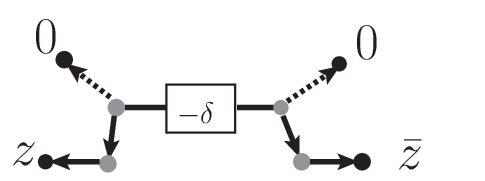}|
\\\lesssim&\includegraphics[height=1cm]{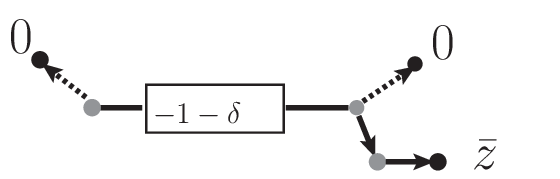}+\includegraphics[height=1cm]{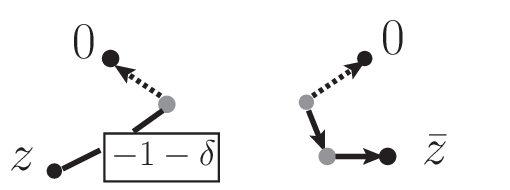}\lesssim |y|^{-\delta}|\bar{y}|^{-\delta}+|\bar{y}|^{-\delta},\endaligned$$
where we used Young's inequality. Again,
 $\delta\hat{\mathcal{W}}^{(\varepsilon;2)}_3\tau, $ can be bounded in a manner similar as before.
 For
$\hat{\mathcal{W}}^{(\varepsilon;2)}_{41}\tau$ we have that for $\delta>0$
 $$\aligned|\langle(\hat{\mathcal{W}}^{(\varepsilon;2)}_{41}\tau)(z),(\hat{\mathcal{W}}^{(\varepsilon;2)}_{41}\tau)(\bar{z})\rangle|
 =&\includegraphics[height=1cm]{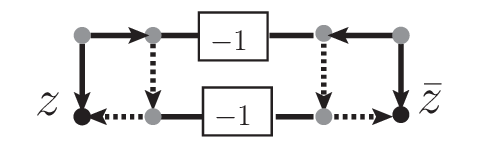}\\
\lesssim&\includegraphics[height=1cm]{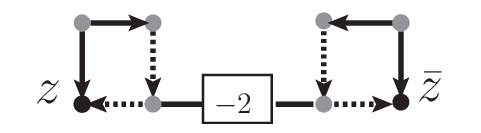}+\includegraphics[height=1cm]{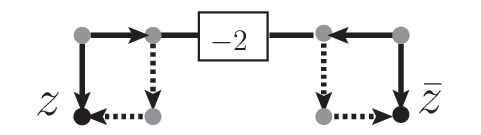}\\\lesssim& |y-\bar{y}|^{-\delta},\endaligned$$
holds uniformly for $\varepsilon\in(0,1]$, where we used Young's inequality. For $\delta\in(0,1)$ we have that
$$\aligned|\langle(\hat{\mathcal{W}}^{(\varepsilon;2)}_{42}\tau)(z),(\hat{\mathcal{W}}^{(\varepsilon;2)}_{42}\tau)(\bar{z})\rangle|
 =&\includegraphics[height=1cm]{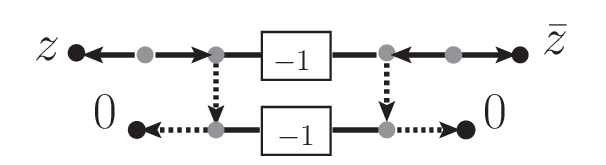}\\
\lesssim&\includegraphics[height=1cm]{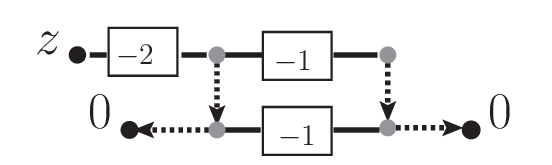}+\includegraphics[height=1cm]{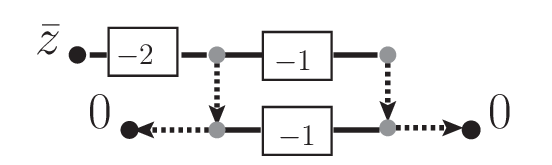}\\
\lesssim&\includegraphics[height=1cm]{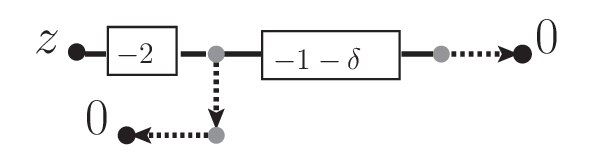}+\includegraphics[height=1cm]{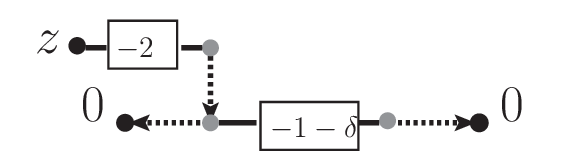}
\\&+\includegraphics[height=1cm]{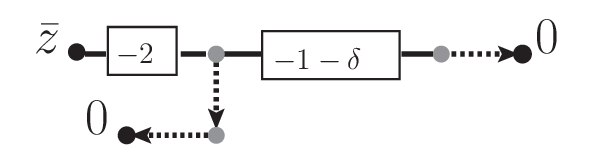}+\includegraphics[height=1cm]{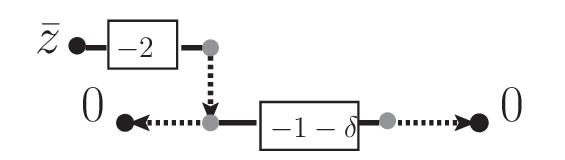}\\
\lesssim&|y|^{-\delta}+|\bar{y}|^{-\delta},\endaligned$$
holds uniformly for $\varepsilon\in(0,1]$, where we used Young's inequality for each   inequality.  Similarly, these bounds also hold for $(\hat{\mathcal{W}}^{(\varepsilon;2)}_5\tau)(z)$. Again, defining $\hat{\mathcal{W}}^{(2)}_i\tau, i=4,5,$ similarly as before and
 $\delta\hat{\mathcal{W}}^{(\varepsilon;2)}_i\tau, i=4,5$ can be bounded in a manner similar as before.

 We now turn to $\hat{\mathcal{W}}^{(\varepsilon;0)}\tau$:
$$(\hat{\mathcal{W}}^{(\varepsilon;0)}\tau)(z)=\sum_{i=1}^2(\hat{\mathcal{W}}^{(\varepsilon;0)}_i\tau)(z)=\sum_{i=1}^2[(\hat{\mathcal{W}}^{(\varepsilon;0)}_{i1}
\tau)(z)
-(\hat{\mathcal{W}}^{(\varepsilon;0)}_{i2}\tau)(z)]-C^{3,\varepsilon}_{ii_1i_2i_3kk_1ll_1jj_1},$$
where
$$(\hat{\mathcal{W}}^{(\varepsilon;0)}_{11}\tau)(z)= \includegraphics[height=1.5cm]{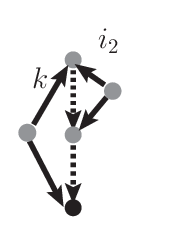}, (\hat{\mathcal{W}}^{(\varepsilon;0)}_{12}\tau)(z)= \includegraphics[height=1.5cm]{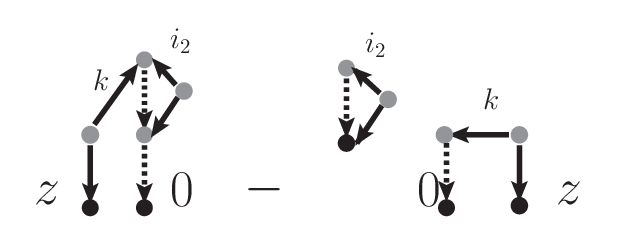},$$
$$(\hat{\mathcal{W}}^{(\varepsilon;0)}_{21}\tau)(z)= \includegraphics[height=1.5cm]{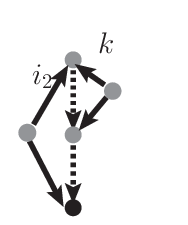}, (\hat{\mathcal{W}}^{(\varepsilon;0)}_{22}\tau)(z)= \includegraphics[height=1.5cm]{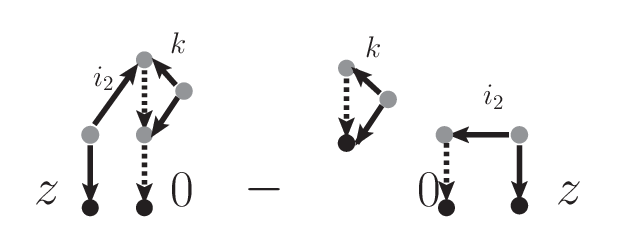},$$
we  choose $C^{3,\varepsilon}_{ii_1i_2i_3kk_1ll_1jj_1}=(\hat{\mathcal{W}}^{(\varepsilon;0)}_{11}\tau)(z)+(\hat{\mathcal{W}}^{(\varepsilon;0)}_{21}\tau)(z).$
By [16, Lemma 10.16] we have that  for every $\delta>0$, $i=1,2$,
$$|(\hat{\mathcal{W}}^{(\varepsilon;0)}_{i2}\tau)(z)|\lesssim|y|^{-\delta}$$
holds uniformly for $\varepsilon\in(0,1]$. Similarly as before, we obtain the bounds for $\delta\hat{\mathcal{W}}^{(\varepsilon;0)}_{i2}\tau$. Then (2.10), (2.11) also follow in this case.

 For $\tau=\mathcal{I}^{ii_1}_l(\mathcal{I}^{ll_1}_k(\mathcal{I}^{l_1l_2}(\Xi_{l_2})\mathcal{I}^{kk_1}(\Xi_{k_1}))
\mathcal{I}^{i_1i_2}(\Xi_{i_2}))\mathcal{I}^{jj_1}(\Xi_{j_1})= \includegraphics[height=0.7cm]{09.eps}$, $i,i_1,i_2,l,l_1,l_2,k,k_1,j,j_1=1,2,3$,
we have similar bounds as above with
 $$C^4_{ii_1i_2kk_1ll_1l_2jj_1}=\includegraphics[height=1.5cm]{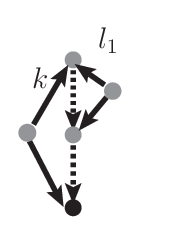}+\includegraphics[height=1.5cm]{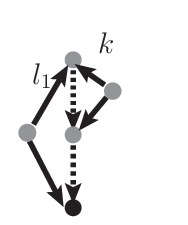}.$$
 $\hfill\Box$

Now combining Theorem 2.17 and Propositions 2.13 and 2.15, we conclude Theorem 1.1 easily.

\section{N-S equation by paracontrolled distributions}

\subsection{Besov spaces and paraproduct}
In the following we recall the definitions and some properties of Besov spaces and paraproducts. For a general introduction to these theories we refer to \cite{BCD11}, \cite{GIP13}. Here the notations are different from the previous section.

First, we introduce the following notations. The space of real valued infinitely differentiable functions of compact support is denoted by $\mathcal{D}(\mathbb{R}^d)$ or $\mathcal{D}$. The space of Schwartz functions is denoted by $\mathcal{S}(\mathbb{R}^d)$. Its dual, the space of tempered distributions is denoted by $\mathcal{S}'(\mathbb{R}^d)$. If $u$ is a vector of $n$ tempered distributions on $\mathbb{R}^d$, then we write $u\in \mathcal{S}'(\mathbb{R}^d,\mathbb{R}^n)$. The Fourier transform and the inverse Fourier transform are denoted by $\mathcal{F}$ and $\mathcal{F}^{-1}$.

 Let $\chi,\theta\in \mathcal{D}$ be nonnegative radial functions on $\mathbb{R}^d$, such that

i. the support of $\chi$ is contained in a ball and the support of $\theta$ is contained in an annulus;

ii. $\chi(z)+\sum_{j\geq0}\theta(2^{-j}z)=1$ for all $z\in \mathbb{R}^d$.

iii. $\textrm{supp}(\chi)\cap \textrm{supp}(\theta(2^{-j}\cdot))=\emptyset$ for $j\geq1$ and $\textrm{supp}(\theta(2^{-i}\cdot))\cap \textrm{supp}(\theta(2^{-j}\cdot))=\emptyset$ for $|i-j|>1$.

We call such a pair $(\chi,\theta)$ a dyadic partition of unity, and for the existence of dyadic partitions of unity  we refer to [1, Proposition 2.10]. The Littlewood-Paley blocks are now defined as
$$\Delta_{-1}u=\mathcal{F}^{-1}(\chi\mathcal{F}u)\quad \Delta_{j}u=\mathcal{F}^{-1}(\theta(2^{-j}\cdot)\mathcal{F}u).$$

For $\alpha\in \mathbb{R}$, the H\"{o}lder-Besov space $\mathcal{C}^\alpha$ is given by $\mathcal{C}^\alpha=B^\alpha_{\infty,\infty}(\mathbb{R}^d,\mathbb{R}^n)$, where for $p,q\in [1,\infty]$ we define
$$B^\alpha_{p,q}(\mathbb{R}^d,\mathbb{R}^n)=\{u=(u^1,...,u^n)\in\mathcal{S}'(\mathbb{R}^d,\mathbb{R}^n):\|u\|_{B^\alpha_{p,q}}=\sum_{i=1}^n(\sum_{j\geq-1}(2^{j\alpha}\|\Delta_ju^i\|_{L^p})^q)^{1/q}<\infty\},$$
with the usual interpretation as the $l^\infty$-norm in case $q=\infty$. We write $\|\cdot\|_{\alpha}$ instead of $\|\cdot\|_{B^\alpha_{\infty,\infty}}$.

We point out that everything above and everything that follows can be applied to distributions on the torus. More precisely, let $\mathcal{D}'(\mathbb{T}^d)$ be the space of distributions on $\mathbb{T}^d$. Therefore, Besov spaces on the torus with general indices $p,q\in[1,\infty]$ are defined as
$$B^\alpha_{p,q}(\mathbb{T}^d,\mathbb{R}^n)=\{u\in\mathcal{S}'(\mathbb{T}^d,\mathbb{R}^n):\|u\|_{B^\alpha_{p,q}}=\sum_{i=1}^n(\sum_{j\geq-1}(2^{j\alpha}\|\Delta_ju^i\|_{L^p(\mathbb{T}^d)})^q)^{1/q}<\infty\}.$$
 We  will need the following Besov embedding theorem on the torus (c.f. [12, Lemma 41]):
\vskip.10in
 \th{Lemma 3.1} Let $1\leq p_1\leq p_2\leq\infty$ and $1\leq q_1\leq q_2\leq\infty$, and let $\alpha\in\mathbb{R}$. Then $B^\alpha_{p_1,q_1}(\mathbb{T}^d)$ is continuously embedded in $B^{\alpha-d(1/p_1-1/p_2)}_{p_2,q_2}(\mathbb{T}^d)$.
\vskip.10in

 Now we recall the following paraproduct introduced by Bony (see \cite{Bon81}). In general, the product $fg$ of two distributions $f\in \mathcal{C}^\alpha, g\in \mathcal{C}^\beta$ is well defined if and only if $\alpha+\beta>0$. In terms of Littlewood-Paley blocks, the product $fg$ can be formally decomposed as
 $$fg=\sum_{j\geq-1}\sum_{i\geq-1}\Delta_if\Delta_jg=\pi_<(f,g)+\pi_0(f,g)+\pi_>(f,g),$$
 with $$\pi_<(f,g)=\pi_>(g,f)=\sum_{j\geq-1}\sum_{i<j-1}\Delta_if\Delta_jg, \quad\pi_0(f,g)=\sum_{|i-j|\leq1}\Delta_if\Delta_jg.$$
We  use the notation
$$S_jf=\sum_{i\leq j-1}\Delta_if.$$
We will use without comment that $\|\cdot\|_\alpha\leq\|\cdot\|_\beta$ for $\alpha\leq\beta$, that $\|\cdot\|_{L^\infty}\lesssim \|\cdot\|_\alpha$ for $\alpha>0$, and that $\|\cdot\|_\alpha\lesssim\|\cdot\|_{L^\infty}$ for $\alpha\leq0$. We will also use that $\|S_ju\|_{L^\infty}\lesssim 2^{-j\alpha}\|u\|_\alpha$ for $\alpha<0$ and $u\in \mathcal{C}^\alpha$.

\vskip.10in
 The basic result about these bilinear operations is given by the following estimates:
\vskip.10in
 \th{Lemma 3.2}(Paraproduct estimates, \cite{Bon81}, [12, Lemma 2]) For any $\beta\in \mathbb{R}$ we have
 $$\|\pi_<(f,g)\|_\beta\lesssim \|f\|_{L^\infty}\|g\|_\beta\quad f\in L^\infty, g\in \mathcal{C}^\beta,$$
 and for $\alpha<0$ furthermore
 $$\|\pi_<(f,g)\|_{\alpha+\beta}\lesssim \|f\|_{\alpha}\|g\|_\beta\quad f\in \mathcal{C}^\alpha, g\in \mathcal{C}^\beta.$$
 For $\alpha+\beta>0$ we have
 $$\|\pi_0(f,g)\|_{\alpha+\beta}\lesssim \|f\|_{\alpha}\|g\|_\beta\quad f\in \mathcal{C}^\alpha, g\in \mathcal{C}^\beta.$$
\vskip.10in
From this lemma we know that $\pi_<(f,g)$ and $\pi_>(f,g)$ are well defined if $f\in L^\infty$. The only term not well defined in defining $fg$ is $\pi_0(f,g).$ Furthermore, if $f$ is smooth, the regularity of $\pi_>(f,g)$ and $\pi_0(f,g)$ will become better than the regularity of $g$. $\pi_<(f,g)$ retains the same regularity as $g$.
\vskip.10in
 The following basic commutator lemma is important for our later use:
\vskip.10in
 \th{Lemma 3.3}([12, Lemma 5]) Assume that $\alpha\in (0,1)$ and $\beta,\gamma\in \mathbb{R}$ are such that $\alpha+\beta+\gamma>0$ and $\beta+\gamma<0$. Then for smooth $f,g,h,$ the trilinear operator
 $$C(f,g,h)=\pi_0(\pi_<(f,g),h)-f\pi_0(g,h)$$ has the bound
 $$\|C(f,g,h)\|_{\alpha+\beta+\gamma}\lesssim\|f\|_\alpha\|g\|_\beta\|h\|_\gamma.$$
 Thus, $C$ can be uniquely extended to a bounded trilinear operator in $L^3(\mathcal{C}^\alpha\times \mathcal{C}^\beta \times \mathcal{C}^\gamma, \mathcal{C}^{\alpha+\beta+\gamma})$.

\vskip.10in
By using this commutator estimate to make sense of the product of $\pi_<(f,g)$ and $h$ for $f\in \mathcal{C}^\alpha,g\in \mathcal{C}^\beta, h\in \mathcal{C}^\gamma$, it is sufficient to define $\pi_0(g,h)$.
\vskip.10in
Now we prove the following commutator estimate for the Leray projection. We follow a similar argument as [4, Lemma A.1]. In the following we use the notation
$f(D)u=\mathcal{F}^{-1}f\mathcal{F}u$.
\vskip.10in
\th{Lemma 3.4} Let $u\in \mathcal{C}^\alpha$ for some $\alpha<1$ and $v\in \mathcal{C}^\beta$ for some $\beta\in \mathbb{R}$. Then for  every $k,l=1,2,3$
$$\|P^{kl}\pi_{<}(u,v)-\pi_<(u,P^{kl}v)\|_{\alpha+\beta}\lesssim \|u\|_\alpha\|v\|_\beta,$$
where  $P$ is the Leray projection.

\proof We have $$P^{kl}\pi_{<}(u,v)-\pi_<(u,P^{kl}v)=\sum_{j=-1}^\infty [P^{kl}(S_{j-1}u\Delta_jv)-S_{j-1}u\Delta _jP^{kl}v]$$
and  every term of this series has a Fourier transform with support in an annulus of the form $2^j\mathcal{A}$ where $\mathcal{A}$ is an annulus. Let $\psi\in \mathcal{D}$ with support in an annulus be such that $\psi =1$ on $\mathcal{A}$. Then
$$P^{kl}(S_{j-1}u\Delta_jv)-S_{j-1}u\Delta _jP^{kl}v=[\hat{P}^{kl}(D),S_{j-1}u]\Delta_jv=[(\psi(2^{-j}\cdot)\hat{P}^{kl})(D),S_{j-1}u]\Delta_jv.$$
Here $\hat{P}^{kl}(x)=\delta_{kl}-\frac{x_k x_l}{|x|^2}$ and  $$[(\psi(2^{-j}\cdot)\hat{P}^{kl})(D),S_{j-1}u]f=(\psi(2^{-j}\cdot)\hat{P}^{kl})(D)(S_{j-1}uf)-S_{j-1}u(\psi(2^{-j}\cdot)\hat{P}^{kl})(D)f$$ denotes the commutator.
By a similar  argument as in the proof of [4, Lemma A.1] we have
$$\|[(\psi(2^{-j}\cdot)\hat{P}^{kl})(D),S_{j-1}u]\Delta_jv\|_{L^\infty}\lesssim\sum_{\eta\in \mathbb{N}^d,|\eta|=1}\|x^\eta\mathcal{F}^{-1}(\psi(2^{-j}\cdot)\hat{P}^{kl})\|_{L^1}\|\partial^\eta S_{j-1}u\|_{L^\infty}\|\Delta_jv\|_{L^\infty}.$$

Moreover, we have the following estimates
$$\aligned &\|x^\eta\mathcal{F}^{-1}(\psi(2^{-j}\cdot)\hat{P}^{kl})\|_{L^1}\\\leq& 2^{-j}\|\mathcal{F}^{-1}(\partial^\eta\psi)(2^{-j}\cdot)\hat{P}^{kl})\|_{L^1}+\|\mathcal{F}^{-1}(\psi(2^{-j}\cdot)\partial^\eta\hat{P}^{kl})\|_{L^1}
\\=&2^{-j}\|\mathcal{F}^{-1}(\partial^\eta\psi(\cdot)\hat{P}^{kl}(2^j\cdot))\|_{L^1}+\|\mathcal{F}^{-1}(\psi(\cdot)\partial^\eta\hat{P}^{kl}(2^j\cdot))\|_{L^1}
\\\lesssim&2^{-j}\|(1+|\cdot|^2)^{d}\mathcal{F}^{-1}(\partial^\eta\psi(\cdot)\hat{P}^{kl}(2^j\cdot))\|_{L^\infty}+ \|(1+|\cdot|^2)^{d}\mathcal{F}^{-1}(\psi(\cdot)\partial^\eta\hat{P}^{kl}( 2^j\cdot))\|_{L^\infty} \endaligned$$
$$\aligned=&2^{-j}\|\mathcal{F}^{-1}((1-\Delta)^{d}(\partial^\eta\psi(\cdot)\hat{P}^{kl}( 2^j\cdot)))\|_{L^\infty}+ \|\mathcal{F}^{-1}((1-\Delta)^{d}(\psi(\cdot)\partial^\eta\hat{P}^{kl}( 2^j\cdot)))\|_{L^\infty}
\\\lesssim&2^{-j}\|(1-\Delta)^{d}(\partial^\eta\psi(\cdot)\hat{P}^{kl}( 2^j\cdot))\|_{L^1}+ \|(1-\Delta)^{d}(\psi(\cdot)\partial^\eta\hat{P}^{kl}( 2^j\cdot))\|_{L^1}
 \\\lesssim&2^{-j}\sum_{0\leq|m|\leq 2d}( 2^j)^{|m|}\frac{1}{( 2^j)^{|m|}}+ \sum_{|m|\leq 2d}( 2^j)^{|m|}\frac{1}{( 2^j)^{|m|+1}}\\\lesssim&2^{-j},\endaligned$$
 where  in the fourth inequality we used $|D^m\hat{P}^{kl}(x)|\lesssim |x|^{-|m|}$ for any multiindices $m$.
Thus we get that
$$\|[\psi(2^{-j}\cdot)\hat{P}^{kl}(D),S_{j-1}u]\Delta_jv\|_{L^\infty}\lesssim 2^{-j(\alpha+\beta)}\|u\|_\alpha\|v\|_\beta,$$which implies the result by a similar argument as in the proof of [4, Lemma A.1].$\hfill\Box$

\vskip.10in
Now we recall the following heat semigroup estimate.
\vskip.10in
\th{Lemma 3.5}([12, Lemma 47]) Let $u\in \mathcal{C}^\alpha$ for some $\alpha\in \mathbb{R}$. Then for every $\delta\geq0$
$$\|P_tu\|_{\alpha+\delta}\lesssim t^{-\delta/2}\|u\|_\alpha,$$
where $P_t$ is the heat semigroup on $\mathbb{T}^d$.

\vskip.10in

For the Leray projection we  have the following estimate on $\mathbb{T}^d$:
\vskip.10in
\th{Lemma 3.6} Let $u\in \mathcal{C}^\alpha$ on $\mathbb{T}^d$ for some $\alpha\in \mathbb{R}$. Then for every $k,l=1,2,3$
$$\|P^{kl}u\|_{\alpha}\lesssim \|u\|_\alpha,$$
where  $P$ is the Leray projection.

\proof Let $\psi\in \mathcal{D}$ with support in an annulus be such that $\psi =1$ on the support of $\theta$. We have that for $j\geq0$
$$\aligned \|\Delta_jP^{kl}u\|_{L^\infty}=&\|\mathcal{F}^{-1}(\hat{P}^{kl}(\cdot)\psi(2^{-j}\cdot))\theta_j\mathcal{F}u\|_{L^\infty}\\\lesssim& \|\mathcal{F}^{-1}(\hat{P}^{kl}(\cdot)\psi(2^{-j}\cdot))\|_{L^1}2^{-j\alpha}\|u\|_\alpha=\|\mathcal{F}^{-1}(\hat{P}^{kl}(2^j\cdot)\psi)\|_{L^1}2^{-j\alpha}\|u\|_\alpha.\endaligned$$
Here $\hat{P}^{kl}(x)=\delta_{kl}-\frac{x^kx^l}{|x|^2}$. By a similar calculaton as in the proof of Lemma 3.4 we obtain that
$$\|\mathcal{F}^{-1}(\hat{P}^{kl}(2^j\cdot)\psi)\|_{L^1}\lesssim \|(1-\Delta)^d(\hat{P}^{kl}(2^j\cdot)\psi)\|_{L^1}\lesssim \sum_{0\leq|m|\leq 2d}(2^j)^{|m|}\frac{1}{(2^j)^{|m|}}\lesssim C.$$
By the theory in \cite{SW71} we know that the above calculations also hold on $\mathbb{T}^d$. Moreover, we have on $\mathbb{T}^d$ for $1<p<\infty$ $$\|\Delta_{-1}P^{kl}u\|_{L^\infty(\mathbb{T}^d)}=\|\mathcal{F}^{-1}\hat{P}^{kl}\chi\mathcal{F}u\|_{L^\infty(\mathbb{T}^d)}\lesssim\|\mathcal{F}^{-1}\hat{P}^{kl}\chi\mathcal{F}u\|_{L^p(\mathbb{T}^d)}\lesssim\|\Delta_{-1}u\|_{L^p(\mathbb{T}^d)}\lesssim\|\Delta_{-1}u\|_{L^\infty(\mathbb{T}^d)},$$
where in the first inequality we used that supp$(\chi\hat{P}\mathcal{F}u)$ is contained in a ball and in the second inequality we used Mihlin's multiplier theorem. Thus the result follows.$\hfill\Box$

\vskip.10in
\subsection{N-S equation}

Let us focus on the equation on  $\mathbb{T}^3$:
$$Lu^i=\sum_{i_1=1}^3P^{ii_1}\xi^{i_1}-\frac{1}{2}\sum_{i_1=1}^3P^{ii_1}(\sum_{j=1}^3D_j(u^iu^j)),\eqno(3.1)$$
$$u(0)=Pu_0\in \mathcal{C}^{-z},$$
where $\xi=(\xi^1,\xi^2,\xi^3)$, $\xi^1, \xi^2, \xi^3$ are the periodic independent space time white noise, $L=\partial_t-\Delta$ and $z\in (1/2,1/2+\delta_0)$ with $0<\delta_0<1/2$. Here without loss of generality we suppose that $\nu=1$.  As we mentioned in the introduction the nonlinear term of this equation is not well defined because of the singularity of $\xi$. In the following we follow the idea of \cite{GIP13} to give the definition of the solution to the equation as a limit of solutions $u^\varepsilon$ to the following equations:
$$Lu^{\varepsilon,i}=\sum_{i_1=1}^3P^{ii_1}\xi^{\varepsilon,i_1}-\frac{1}{2}\sum_{i_1=1}^3P^{ii_1}(\sum_{j=1}^3D_j(u^\varepsilon u^{\varepsilon,j})),$$
$$u(0)=Pu_0\in \mathcal{C}^{-z}.$$
Here $\xi^\varepsilon$ is a family of smooth approximations  of $\xi$ such that $\xi^\varepsilon\rightarrow\xi$ as $\varepsilon\rightarrow0$. Now we prove a uniform estimate for $u^\varepsilon$.

 In the following to avoid heavy notation we omit the dependence on $\varepsilon$ if there's no confusion  and consider (3.1) for smooth $\xi$. We split the equation (3.1) into the following four equations:
$$Lu_1^i=\sum_{i_1=1}^3P^{ii_1}\xi^{i_1},$$
$$Lu_2^i=-\frac{1}{2}\sum_{i_1=1}^3P^{ii_1}(\sum_{j=1}^3D_j(u_1^{i_1}\diamond u_1^j))\quad u_2(0)=0,$$
$$Lu_3^i=-\frac{1}{2}\sum_{i_1=1}^3P^{ii_1}(\sum_{j=1}^3D_j(u_1^{i_1}\diamond u_2^j+u_2^{i_1}\diamond u_1^j)),\quad u_3(0)=0,$$
and $$\aligned Lu_4^i=&-\frac{1}{2}\sum_{i_1,j=1}^3P^{ii_1}D_j[u_1^{i_1}\diamond (u_3^j+u_4^j)+(u_3^{i_1} +u_4^{i_1})\diamond u_1^j+u_2^{i_1}\diamond u_2^j\\&+u_2^{i_1}(u_3^j+u_4^j)+u_2^j(u_3^{i_1}+u_4^{i_1})+(u_3^{i_1}+u_4^{i_1})(u_3^{j}+u_4^{j})],\endaligned\eqno(3.2)$$
$$u_4(0)=Pu_0-u_1(0),$$
where for $i,j=1,2,3$ $$u_1^i\diamond u_3^j=\pi_<(u_3^j,u_1^i)+\pi_>(u_3^j,u_1^i)+\pi_{0,\diamond}(u_3^j,u_1^i)$$
and $$u_1^i\diamond u_4^j=\pi_<(u_4^j,u_1^i)+\pi_>(u_4^j,u_1^i)+\pi_{0,\diamond}(u_4^j,u_1^i).$$
Here for $i=1,2,3$, $u_1^i(t)=\int_{-\infty}^t\sum_{i_1=1}^3P^{ii_1}P_{t-s}\xi^{\varepsilon,i_1} ds$ and we use $\diamond$  to replace the product of some terms, the meaning of which will be given later. In fact, the product of these terms needs to be renormalised such that they converge as $\varepsilon\rightarrow0$. We will discuss this in Section 3.3 below. The results for the renormalised terms not including $u_4$ can be  proved by using a similar idea as in the definition of Wick products. However, $u_4\diamond u_1$ cannot be defined by this trick since $u_4$ is the unknown. To deal with this term we will use the fact that $u_4$ has a specific structure since it satisfies (3.2). Now we do some preparations.  Consider the following equations: $$LK^i=u_1^i,\quad K^i(0)=0.$$Then we obtain that for every $\delta>0$ small enough, if $u_1^i\in C([0,T];\mathcal{C}^{-\frac{1}{2}-\frac{\delta}{2}})$, then $K^i\in C([0,T];\mathcal{C}^{\frac{3}{2}-\delta})$  and by Lemma 3.5
$$\|K^i(t)\|_{\frac{3}{2}-\delta}\lesssim t^{\delta/4}\sup_{s\in[0,t]}\|u_1^i(s)\|_{-1/2-\delta/2}.\eqno(3.3)$$

First we  assume that $u_1^i\in C([0,T];\mathcal{C}^{-\frac{1}{2}-\frac{\delta}{2}})$,  $u_1^i\diamond u_1^j\in C([0,T];\mathcal{C}^{-1-\delta/2})$, $u_1^i\diamond u_2^j=u_2^j\diamond u_1^i\in C([0,T];\mathcal{C}^{-1/2-\delta/2})$, $u_2^i\diamond u_2^j\in C([0,T];\mathcal{C}^{-\delta})$, $\pi_{0,\diamond}(u_3^i,u_1^j)\in C([0,T];\mathcal{C}^{-\delta})$ and $\pi_{0,\diamond}(P^{ii_1}D_jK^j,u_1^{j_1}),$ $\pi_{0,\diamond}(P^{ii_1}D_jK^{i_1},u_1^{j_1})\in C([0,T];\mathcal{C}^{-\delta})$ for $i,j, i_1,j_1=1,2,3,$ and that
$$\aligned C_\xi^\varepsilon:=&\sup_{t\in[0,T]}\bigg[\sum_{i=1}^3\|u_1^{\varepsilon,i}\|_{-1/2-\delta/2}+\sum_{i,j=1}^3\|u_1^{\varepsilon,i}\diamond u_1^{\varepsilon,j}\|_{-1-\delta/2}+\sum_{i,j=1}^3\|u_1^{\varepsilon,i}\diamond u_2^{\varepsilon,j}\|_{-1/2-\delta/2}\\&+\sum_{i,j=1}^3\|u_2^{\varepsilon,i}\diamond u_2^{\varepsilon,j}\|_{-\delta}+\sum_{i,j=1}^3\|\pi_{0,\diamond}( u_3^{\varepsilon,i},u_1^{\varepsilon,j})\|_{-\delta}+\sum_{i,i_1,j,j_1=1}^3\|\pi_{0,\diamond}(P^{ii_1}D_{j}K^{\varepsilon,j},u_1^{\varepsilon,j_1})\|_{-\delta}
\\&+\sum_{i,i_1,j,j_1=1}^3\|\pi_{0,\diamond}
(P^{ii_1}D_{j}K^{\varepsilon,i_1},u_1^{\varepsilon,j_1})\|_{-\delta}\bigg]<\infty.\endaligned$$ By Lemmas 3.5 and 3.6 we easily deduce that
$u_2^i\in C([0,T];\mathcal{C}^{-\delta}),u_3^i\in C([0,T];\mathcal{C}^{1/2-\delta})$  for $i=1,2,3$, and that $$\sup_{t\in[0,T]}(\sum_{i=1}^3\|u_2^i\|_{-\delta}+\sum_{i=1}^3\|u_3^i\|_{1/2-\delta})\lesssim C_\xi.\eqno(3.4)$$

In the following we will fix $\delta>0$ small enough such that
$$\delta<\delta_0\wedge\frac{1-2\delta_0}{3}\wedge\frac{1-z}{4}\wedge(2z-1).$$
 By a fixed point argument it is easy to obtain local existence and uniqueness of solution to equation (3.2): More precisely, for each $\varepsilon\in(0,1)$ there exists a maximal time $T_\varepsilon$ and  $u_4\in C((0,T_\varepsilon);\mathcal{C}^{1/2-\delta_0})$ with respect to the norm
$\sup_{t\in [0,T]}t^{\frac{1/2-\delta_0+z}{2}}\|u_4(t)\|_{1/2-\delta_0}$ such that $u_4$ satisfies equation (3.2) before $T_\varepsilon$ and
$$\sup_{t\in [0,T_\varepsilon)}t^{\frac{1/2-\delta_0+z}{2}}\|u_4(t)\|_{1/2-\delta_0}=\infty.$$
Indeed, since $\xi_\varepsilon$ is smooth, by (3.2) and Lemmas 3.5 and 3.6 we have the following estimate $$\sup_{t\in [0,T]}t^{\frac{1/2-\delta_0+z}{2}}\|u_4(t)\|_{1/2-\delta_0}\lesssim C_\varepsilon(\|u_0\|_{-z},u_1,u_2,u_3)+T^{\frac{1/2+\delta_0-z}{2}}(\sup_{t\in [0,T]}t^{\frac{1/2-\delta_0+z}{2}}\|u_4(t)\|_{1/2-\delta_0})^2,$$
where $C_\varepsilon(\|u_0\|_{-z},u_1,u_2,u_3)$ are constants depending on $\varepsilon$ and we used $z<1/2+\delta_0$.

\textbf{Paracontrolled ansatz}: As we mentioned before, our problem lies in how to define $\pi_{0}(u_4^j,u_1^i)$. Observing that the worst term on the right hand side of (3.2) is $PD\pi_<(u_3+u_4,u_1)$, we write $u_4$ as the following paracontrolled ansatz for $i=1,2,3$:
$$u_4^i=-\frac{1}{2}\sum_{i_1=1}^3P^{ii_1}(\sum_{j=1}^3D_j[\pi_<(u_3^{i_1}+u_4^{i_1},K^j)+\pi_<(u_3^{j}+u_4^{j},K^{i_1})])+u^{\sharp,i}$$ with $u^{\sharp,i}(t)\in \mathcal{C}^{1/2+\beta}$ for some $\delta/2<\beta<(z+2\delta-1/2)<(1/2-2\delta)$ and $t\in(0,T_\varepsilon)$ (which can be done for fixed $\varepsilon>0$ since $\xi_\varepsilon$ is smooth and by (3.2) we note that
$$t^{\frac{1/2+\beta+z}{2}}\|u_4(t)\|_{1/2+\beta}\lesssim C_\varepsilon(\|u_0\|_{-z},u_1,u_2,u_3)+t^{\frac{1/2+\delta_0-z}{2}}(\sup_{s\in [0,t]}s^{\frac{1/2-\delta_0+z}{2}}\|u_4(s)\|_{1/2-\delta_0})^2).$$
From the paracontrolled ansatz and Lemma 3.2 we easily get the following estimate for $i=1,2,3$:
$$\|u_4^i\|_{1/2-\delta}\lesssim \sum_{i_1,j=1}^3\|u_3^{i_1}+u_4^{i_1}\|_{1/2-\delta_0}\|K^j\|_{3/2-\delta}+\|u^{\sharp,i}\|_{1/2+\beta}.\eqno(3.5)$$
Moreover $u_4$ solves (3.2) if and only if $u^\sharp$ solves the following equation:
$$\aligned &Lu^{\sharp,i}= -\frac{1}{2}\sum_{i_1,j=1}^3P^{ii_1}D_j\big{[}u_2^{i_1}\diamond u_2^j+u_2^{i_1}(u_3^j+u_4^j)+u_2^{j}(u_3^{i_1}+u_4^{i_1})+(u_3^{i_1}+u_4^{i_1})(u_3^j+u_4^j)\\&-\pi_<(L(u_3^{i_1}+u_4^{i_1}),K^j)
+2\sum_{l=1}^3\pi_<(D_l(u_3^{i_1}+u_4^{i_1}),D_lK^j)+\pi_>(u_3^{i_1}+u_4^{i_1},u_1^j)+\pi_{0,\diamond}(u_3^{i_1},u_1^j)+\pi_{0,\diamond}(u_4^{i_1},u_1^j)
\\&-\pi_<(L(u_3^{j}+u_4^{j}),K^{i_1})
+2\sum_{l=1}^3\pi_<(D_l(u_3^{j}+u_4^{j}),D_lK^{i_1})+\pi_>(u_3^{j}+u_4^{j},u_1^{i_1})+\pi_{0,\diamond}(u_3^{j},u_1^{i_1})
+\pi_{0,\diamond}(u_4^{j},u_1^{i_1})\big{]}\\&:=\phi^{\sharp,i}.\endaligned\eqno(3.6)$$

\noindent\textbf{Renormalisation of $\pi_{0}(u_4^i,u_1^j)$}: By the paracontrolled ansatz we have for $i,j=1,2,3$,
$$\aligned\pi_{0}(u_4^i,u_1^j)=&-\frac{1}{2}\pi_{0}(\sum_{i_1,j_1=1}^3P^{ii_1}\pi_<(u_3^{i_1}+u_4^{i_1},D_{j_1}K^{j_1}),u_1^j)
-\frac{1}{2}\pi_{0}(\sum_{i_1,j_1=1}^3P^{ii_1}\pi_<(u_3^{j_1}+u_4^{j_1},D_{j_1}K^{i_1}),u_1^j)
\\&-\frac{1}{2}\sum_{i_1,j_1=1}^3\pi_0(P^{ii_1}\pi_<(D_{j_1}(u_3^{i_1}+u_4^{i_1}),K^{j_1}),u_1^j))
-\frac{1}{2}\sum_{i_1,j_1=1}^3\pi_0(P^{ii_1}\pi_<(D_{j_1}(u_3^{j_1}+u_4^{j_1}),K^{i_1}),u_1^j))\\&+\pi_0(u^{\sharp,i},u_1^j).\endaligned$$
The last three terms can be easily controlled by Lemma 3.2, and it is sufficient to consider the first two terms:
For $i,i_1,j,j_1=1,2,3$,
$$\aligned&\pi_{0}(P^{ii_1}\pi_<(u_3^{i_1}+u_4^{i_1},D_{j_1}K^{j_1}),u_1^{j})\\=&\pi_0(P^{ii_1}\pi_<(u_3^{i_1}+u_4^{i_1},D_{j_1}K^{j_1}),u_1^{j})
-\pi_0(\pi_<(u_3^{i_1}+u_4^{i_1},P^{ii_1}
D_{j_1}K^{j_1}),u_1^{j})\\&+\pi_0(\pi_<(u_3^{i_1}+u_4^{i_1},P^{ii_1}D_{j_1}K^{j_1}),u_1^{j})-(u_3^{i_1}+u_4^{i_1})\pi_0(P^{ii_1}D_{j_1}K^{j_1}
,u_1^j)\\&+(u_3^{i_1}+u_4^{i_1})\pi_{0}
(P^{ii_1}D_{j_1}K^{j_1},u_1^{j}).\endaligned$$
Applying Lemmas 3.3 and 3.4 we can control  the first four terms on the right hand side of above equality. As we mentioned above for $\pi_{0}
(P^{ii_1}D_{j_1}K^{j_1},u_1^{j})$ we need to do renormalisation to make it convergent as $\varepsilon\rightarrow0$, which leads to the renormalisation of $\pi_{0}(u_4^i,u_1^j)$. Define
$$\aligned&\pi_{0,\diamond}(u_4^i,u_1^j)\\:=&-\frac{1}{2}(\pi_{0,\diamond}(\sum_{i_1,j_1=1}^3P^{ii_1}\pi_<(u_3^{i_1}+u_4^{i_1},D_{j_1}K^{j_1}),u_1^j)
+\pi_{0,\diamond}(\sum_{i_1,j_1=1}^3P^{ii_1}\pi_<(u_3^{j_1}+u_4^{j_1},D_{j_1}K^{i_1}),u_1^j)
\\&+\sum_{i_1,j_1=1}^3\pi_0(P^{ii_1}\pi_<(D_{j_1}(u_3^{i_1}+u_4^{i_1}),K^{j_1}),u_1^j))
+\sum_{i_1,j_1=1}^3\pi_0(P^{ii_1}\pi_<(D_{j_1}(u_3^{j_1}+u_4^{j_1}),K^{i_1}),u_1^j))\\&+\pi_0(u^{\sharp,i},u_1^j),\endaligned$$
where
$$\aligned&\pi_{0,\diamond}(P^{ii_1}\pi_<(u_3^{i_1}+u_4^{i_1},D_{j_1}K^{j_1}),u_1^{j})\\:=&\pi_0(P^{ii_1}\pi_<(u_3^{i_1}+u_4^{i_1},D_{j_1}K^{j_1}),u_1^{j})
-\pi_0(\pi_<(u_3^{i_1}+u_4^{i_1},P^{ii_1}
D_{j_1}K^{j_1}),u_1^{j})\\&+\pi_0(\pi_<(u_3^{i_1}+u_4^{i_1},P^{ii_1}D_{j_1}K^{j_1}),u_1^{j})-(u_3^{i_1}+u_4^{i_1})\pi_0(P^{ii_1}D_{j_1}K^{j_1}
,u_1^j)\\&+(u_3^{i_1}+u_4^{i_1})\pi_{0,\diamond}
(P^{ii_1}D_{j_1}K^{j_1},u_1^{j}),\endaligned$$
and $\pi_{0,\diamond}(P^{ii_1}\pi_<(u_3^{j_1}+u_4^{j_1},D_{j_1}K^{i_1}),u_1^j)$ can be defined similarly. Using Lemmas 3.2 and 3.3 we get that for $\delta\leq\delta_0<1/2-3\delta/2$
$$\aligned &\|\pi_{0,\diamond}(P^{ii_1}\pi_<(u_3^{i_1}+u_4^{i_1},D_{j_1}K^{j_1}),u_1^{j})\|_{-\delta}\\\lesssim & \|P^{ii_1}\pi_<(u_3^{i_1}+u_4^{i_1},D_{j_1}K^{j_1})-\pi_<(u_3^{i_1}+u_4^{i_1},P^{ii_1}D_{j_1}K^{j_1})\|_{1-\delta-\delta_0}\|u_1^{j}\|_{-1/2-\delta/2}\\&+\|u_3^{i_1}+u_4^{i_1}
\|_{1/2-\delta_0}\|P^{ii_1}
D_{j_1}K^{j_1}\|_{1/2-\delta}\|u_1^{j}\|_{-1/2-\delta/2}\\&+\|u_3^{i_1}+u_4^{i_1}\|_{1/2-\delta_0}\|\pi_{0,\diamond}(P^{ii_1}D_{j_1}K^{j_1},u_1^{j})\|_{-\delta}
\\\lesssim & \|u_3^{i_1}+u_4^{i_1}\|_{1/2-\delta_0}\|K^{j_1}\|_{3/2-\delta}\|u_1^{j}\|_{-1/2-\delta/2}+
\|u_3^{i_1}+u_4^{i_1}\|_{1/2-\delta_0}\|\pi_{0,\diamond}(P^{ii_1}D_{j_1}K^{j_1},u_1^{j})\|_{-\delta}.\endaligned$$
Here in the last inequality we used Lemmas 3.4 and 3.6. Similar estimates can also be deduced for
$\pi_{0,\diamond}(\sum_{i_1,j_1=1}^3P^{ii_1}\pi_<(u_3^{j_1}+u_4^{j_1},D_{j_1}K^{i_1}),u_1^j)$.

Hence we obtain that for $i,j=1,2,3,$
$$\aligned \|\pi_{0,\diamond}(u_4^i,u_1^j)\|_{-\delta}\lesssim & \sum_{i_1=1}^3\|u_3^{i_1}+u_4^{i_1}\|_{1/2-\delta_0}\sum_{j_1=1}^3\|K^{j_1}\|_{3/2-\delta}\|u_1^j\|_{-1/2-\delta/2}\\&
+\sum_{i_1,j_1=1}^3\|u_3^{i_1}+u_4^{i_1}\|_{1/2-\delta_0}
\|\pi_{0,\diamond}(P^{ii_1}D_{j_1}K^{j_1},u_1^j)\|_{-\delta}\endaligned$$
$$\aligned&+\sum_{i_1,j_1=1}^3\|u_3^{j_1}+u_4^{j_1}\|_{1/2-\delta_0}\|\pi_{0,\diamond}(P^{ii_1}D_{j_1}K^{i_1},u_1^j)\|_{-\delta}\\&+\|u^{\sharp,i}\|_{1/2+\beta}\|u_1^{j}\|_{-1/2-\delta/2}
\\\lesssim &C_\xi^3+1+\|u_4\|_{1/2-\delta_0}(C_\xi^2+1)+\|u^\sharp\|C_\xi.\endaligned$$
\vskip.10in
\textbf{Estimate of $\phi^\sharp$}: To obtain a uniform estimate for $u_4^\varepsilon$, we first prove an estimate for $\phi^\sharp$:
\vskip.10in
\th{Lemma 3.7} For $\phi^\sharp$ defined in (3.6), the following estimate holds:
$$\|\phi^{\sharp,i}\|_{-1-2\delta}\lesssim (1+C_\xi^4)\big{[}1+\|u^{\sharp}\|_{1/2+\beta}+\|u_4\|_{1/2-\delta_0}
+\|u_4\|_{\delta}^2\big{]}.\eqno(3.7)$$

\proof First we consider $\pi_<(L(u_3^i+u_4^i),K^j), i,j=1,2,3,$:
Indeed (3.2) implies that for $i=1,2,3,$
$$\aligned L(u_3^i+u_4^i)=&-\frac{1}{2}\sum_{i_1,j=1}^3P^{ii_1}D_j(u_1^{i_1}\diamond u_2^j+u_1^{j}\diamond u_2^{i_1}+u_1^{i_1}\diamond(u_3^j+u_4^j)+u_1^j\diamond (u_3^{i_1} +u_4^{i_1})\\&+u_2^{i_1}\diamond u_2^j+u_2^{i_1}(u_3^j+u_4^j)+u_2^j(u_3^{i_1}+u_4^{i_1})+(u_3^{i_1}+u_4^{i_1})(u_3^{j}+u_4^{j})),\endaligned$$
where for $i,j=1,2,3,$
$$u_1^i\diamond (u_3^j+u_4^j)=\pi_<(u_3^j+u_4^j,u_1^i)+\pi_{0,\diamond}(u_3^j,u_1^i)+\pi_>(u_3^j+u_4^j,u_1^i)+\pi_{0,\diamond}(u_4^j,u_1^i).$$
Using Lemmas 3.6 and 3.2 we obtain that for $i=1,2,3,$ $$\aligned\|L(u_3^i+u_4^i)\|_{-3/2-\delta/2}\lesssim&\sum_{i_1,j_1=1}^3[\|u_1^{i_1}\diamond u_2^{j_1}\|_{-1/2-\delta/2}+\|u_2^{i_1}\diamond u_2^{j_1}\|_{-\delta}+\|u_1^{i_1}\|_{-1/2-\delta/2}\|u_3^{j_1}+u_4^{j_1}\|_{1/2-\delta_0}\\&+\|\pi_{0,\diamond}( u_3^{i_1},u_1^{j_1})\|_{-\delta}+\|u_2^{i_1}\|_{-\delta}\|u_3^{j_1}+u_4^{j_1}\|_{1/2-\delta_0}\\&+\|u_3^{i_1}+u_4^{i_1}\|_{\delta}\|u_3^{j_1}+u_4^{j_1}\|_{\delta}
+\|\pi_{0,\diamond}(u_4^{i_1},u_1^{j_1})\|_{-\delta}]
\\\lesssim&C_\xi^3+1+(1+C_\xi^2)\|u_4\|_{1/2-\delta_0}+C_\xi\|u^\sharp\|_{1/2+\beta}+\|u_4\|_\delta^2,\endaligned$$
where we used $\delta<\delta_0\wedge(\frac{1}{2}-\delta_0)$,
which by Lemma 3.2 yields that
$$\aligned&\|\pi_<(L(u_3^i+u_4^i),K^j)\|_{-3\delta/2}\\\lesssim&\|K^{j}\|_{3/2-\delta}[C_\xi^3+1+(1+C_\xi^2)\|u_4\|_{1/2-\delta_0}+C_\xi\|u^\sharp\|_{1/2+\beta}+\|u_4\|_\delta^2].\endaligned$$
Then we consider $\pi_<(D_l(u_3^{i_1}+u_4^{i_1}),D_lK^j)+\pi_>(u_3^{i_1}+u_4^{i_1},u_1^j)$ for $i_1,l,j=1,2,3$ in (3.6): Indeed  Lemma 3.2 implies that
$$\aligned &\|\pi_<(D_l(u_3^{i_1}+u_4^{i_1}),D_lK^j)+\pi_>(u_3^{i_1}+u_4^{i_1},u_1^j)\|_{-2\delta}\\\lesssim & (\|u_3^{i_1}\|_{1/2-\delta}+\|u_4^{i_1}\|_{1/2-\delta})(\|K^j\|_{3/2-\delta}+\|u_1^j\|_{-1/2-\delta/2})\\\lesssim & (\|u_3^{i_1}\|_{1/2-\delta}+\sum_{i_2,j_1=1}^3\|u_3^{i_2}+u_4^{i_2}\|_{1/2-\delta_0}\|K^{j_1}\|_{3/2-\delta}+\|u^{\sharp,i_1}\|_{1/2+\beta})
C_\xi,\endaligned$$
where in the last inequality we used (3.5).

Combining all these estimates obtained above, by (3.6) we get that
$$\aligned&\|\phi^{\sharp,i}\|_{-1-2\delta}\\\lesssim &\sum_{j=1}^3(\|K^j\|_{3/2-\delta}+1)\sum_{i_1,j_1=1}^3[\|u_2^{i_1}\diamond u_2^{j_1}\|_{-\delta}+\|\pi_{0,\diamond}( u_3^{i_1},u_1^{j_1})\|_{-\delta}+\|u_2^{i_1}\|_{-\delta}\|u_3^{j_1}+u_4^{j_1}\|_{1/2-\delta_0}\\&+\|u_3^{i_1}+u_4^{i_1}\|_{\delta}\|u_3^{j_1}+u_4^{j_1}\|_{\delta}
+C_\xi^3+1+(1+C_\xi^2)\|u_4\|_{1/2-\delta_0}+C_\xi\|u^\sharp\|_{1/2+\beta}+\|u_4\|_\delta^2]\\&+\sum_{i_1,j_1,l=1}^3(\|u_3^{i_1}\|_{1/2-\delta}
+\sum_{i_2,j_1=1}^3\|u_3^{i_2}+u_4^{i_2}\|_{1/2-\delta_0}\|K^{j_1}\|_{3/2-\delta}+\|u^{\sharp,i_1}\|_{1/2+\beta})C_\xi
\\\lesssim &(1+C_\xi^4)\big{[}1+\|u^{\sharp}\|_{1/2+\beta}+\|u_4\|_{1/2-\delta_0}
+\|u_4\|_{\delta}^2\big{]},\endaligned$$
where we used (3.2) (3.3) and $\delta\leq \delta_0$ in the last inequality.$\hfill\Box$
\vskip.10in

\textbf{Construction of the solution:} In the following we will prove a uniform estimate of $u_4^\varepsilon$: By the paracontrolled ansatz (3.3) and Lemma 3.2 we get
$$\|u_4^{i}(t)\|_{1/2-\delta_0}\lesssim t^{\delta/4}C_\xi\sum_{i_1=1}^3\|u_3^{i_1}(t)+u_4^{i_1}(t)\|_{1/2-\delta_0}+\|u^{\sharp,i}(t)\|_{1/2-\delta_0},$$
which shows that for $t\in[0,\bar{T}]$ (with $\bar{T}>0$ only depending on $C_\xi$)
$$ \sum_{i=1}^3\|u_4^{i}(t)\|_{1/2-\delta_0}\lesssim C_\xi^2+ \sum_{i=1}^3\|u^{\sharp,i}(t)\|_{1/2-\delta_0}.\eqno(3.8)$$
Similarly, we have for $t\in[0,\bar{T}]$ (with $\bar{T}>0$ only depending on $C_\xi$)
$$ \sum_{i=1}^3\|u_4^{i}(t)\|_{\delta}\lesssim C_\xi^2+\sum_{i=1}^3\|u^{\sharp,i}(t)\|_{\delta}.\eqno(3.9)$$
Moreover, Lemma 3.5 and (3.6) yield that for $\delta+z<1$
$$\aligned &t^{\delta+z}\| u^{\sharp}(t)\|_{1/2+\beta}\\\lesssim & \|Pu_0-u_1(0)\|_{-z}+t^{\delta+z}\int_0^t(t-s)^{-3/4-\delta-\beta/2}s^{-(\delta+z)}s^{\delta+z}\|\phi^{\sharp}(s)\|_{-1-2\delta}ds,\endaligned\eqno(3.10)$$
where we used the condition on $\beta$ to deduce that $3/4+\beta/2+\delta<1$ and $\frac{1/2+\beta+z}{2}\leq\delta+z$.
Similarly, we deduce that
$$\aligned t^{\delta+z}\|u^{\sharp}(t)\|^2_{\delta}\lesssim & \|Pu_0-u_1(0)\|^2_{-z}+t^{\delta+z}(\int_0^t(t-s)^{-\frac{1+3\delta}{2}}s^{-(\delta+z)}s^{\delta+z}\|\phi^{\sharp}(s)\|_{-1-2\delta}ds)^2 \\\lesssim &\|Pu_0-u_1(0)\|^2_{-z}+t^{(1-3\delta)/2}\int_0^t(t-s)^{-\frac{1+3\delta}{2}}s^{-(\delta+z)} (s^{\delta+z}\|\phi^{\sharp}(s)\|_{-1-2\delta})^2ds.\endaligned\eqno(3.11)$$
Here in the last inequality we used H\"{o}lder's inequality. Thus, by (3.7-3.11) we get that for $t\in[0,\bar{T}]$
$$\aligned &t^{\delta+z}\|\phi^{\sharp}\|_{-1-2\delta}\lesssim (1+C_\xi^4)\big{[}\|Pu_0-u_1(0)\|^2_{-z}+C_\xi^4+1\\&+\int_0^tt^{\delta+z}(t-s)^{-3/4-\delta-\beta/2}s^{-(\delta+z)}(s^{\delta+z}\|\phi^{\sharp}
(s)\|_{-1-2\delta})\\&+ t^{(1-3\delta)/2}(t-s)^{-\frac{1+3\delta}{2}}s^{-(\delta+z)}(s^{\delta+z}\|\phi^{\sharp}(s)\|_{-1-2\delta})^2ds\big{]}.\endaligned$$
Then Bihari's inequality implies that for $\delta<\frac{1-z}{4}$ there exists some $0<T_0\leq \bar{T}$ such that
$$\sup_{t\in[0,T_0]}t^{\delta+z}\|\phi^{\sharp}\|_{-1-2\delta}\lesssim C(T_0,C_\xi,\|u_0\|_{-z}),\eqno(3.12)$$
where $C(T_0,C_\xi,\|u_0\|_{-z})$ depends on $T_0, \|u_0\|_{-z}$ and $C_\xi$. Here $T_0$ can be chosen independent of $\varepsilon$ such that (3.12) holds for all $\varepsilon\in(0,1)$, if $C_\xi^\varepsilon$ and $\|u_0\|_{-z}$ is uniformly bounded over $\varepsilon\in(0,1)$. Similarly as (3.10) we have
$$\aligned &t^{(1/2-\delta_0+z)/2}\| u^{\sharp}(t)\|_{1/2-\delta_0}\\\lesssim & \|Pu_0-u_1(0)\|_{-z}+t^{(1/2-\delta_0+z)/2}\int_0^t(t-s)^{-3/4-\delta+\delta_0/2}s^{-(\delta+z)}s^{\delta+z}\|\phi^{\sharp}(s)\|_{-1-2\delta}ds
\\\lesssim &\|Pu_0-u_1(0)\|_{-z}+t^{(1-4\delta-z)/2}\sup_{s\in[0,t]}s^{\delta+z}\|\phi^{\sharp}(s)\|_{-1-2\delta}.\endaligned\eqno(3.13)$$
Then by (3.8) (3.13) we obtain that
$$\sup_{t\in [0,T_0]}t^{\frac{1/2-\delta_0+z}{2}}\|u_4(t)\|_{1/2-\delta_0}\lesssim C_\xi^2+\|u_0\|_{-z}+C(T_0,C_\xi,\|u_0\|_{-z}),$$
which implies that $T_\varepsilon\geq T_0$. Here we used $z\geq1/2+\delta/2$.
Moreover, similarly as for (3.8) one also gets that
for $t\in [0,T_0]$
$$\aligned \|u_4(t)\|_{-z}&\lesssim C_\xi^2+\|u^\sharp(t)\|_{-z}\\&\lesssim C_\xi^2+\|u_0\|_{-z}+\int_0^t(t-s)^{\frac{-1-2\delta+z}{2}}s^{-(\delta+z)}
s^{\delta+z}\|\phi^{\sharp,\lambda}\|_{-1-2\delta}ds,\endaligned$$
where in the last inequality we used Lemma 3.5. This gives us our final estimate for $u^4$:
$$\sup_{t\in[0,T_0]}\|u_4(t)\|_{-z}\lesssim C_\xi^2+\|u_0\|_{-z}+C(T_0,C_\xi,\|u_0\|_{-z}).$$
We define $\mathbb{Z}(\xi^\varepsilon):=(u_1^\varepsilon, u_1^\varepsilon\diamond u_1^\varepsilon,u_1^\varepsilon\diamond u_2^\varepsilon,u_2^\varepsilon\diamond u_2^\varepsilon,\pi_{0,\diamond}( u_3^\varepsilon,u_1^\varepsilon),\pi_{0,\diamond}(PDK^\varepsilon,u_1^\varepsilon))\in \mathbb{X}:=C([0,T];\mathcal{C}^{-1/2-\delta/2})\times  C([0,T];\mathcal{C}^{-1-\delta/2})\times C([0,T];\mathcal{C}^{-1/2-\delta/2})\times C([0,T];\mathcal{C}^{-\delta})\times C([0,T];\mathcal{C}^{-\delta})\times C([0,T];\mathcal{C}^{-\delta})$. Here $\mathbb{X}$ is equipped with the product topology.

Similar arguments show that for every $a>0$ there exists a sufficiently small $T_0>0$ such that the map $( u_0,\mathbb{Z}(\xi_\varepsilon))\mapsto u_4$ is Lipschitz continuous on the set $$\max\{\|u_0\|_{-z},C_\xi\}\leq a.$$Here we consider $u_4$ with respect to the norm given by $$\sup_{t\in [0,T_0]}\|u_4(t)\|_{-z}.$$
Hence we obtain that there exists a local solution $u$ to (3.1) with initial condition $u_0$, which is the limit of the solutions $u^\varepsilon,\varepsilon>0,$ to the following equation $$Lu^{\varepsilon,i}=\sum_{i_1=1}^3P^{ii_1}\xi^{\varepsilon,i_1}-\frac{1}{2}\sum_{i_1=1}^3P^{ii_1}(\sum_{j=1}^3D_j(u^{\varepsilon,i_1} u^{\varepsilon,j}))\quad u^\varepsilon(0)=u_0,$$
provided that $\mathbb{Z}(\xi^\varepsilon)$ converges in $\mathbb{X}$, i.e.
for $i,i_1,j,j_2=1,2,3$, there exist $v_1^{i},v_2^{ij},v_3^{ij},v_4^{ij},v_5^{ij},v_6^{ii_1jj_2},v_7^{ii_1jj_2}$ such that for any $\delta>0$, $u_1^{\varepsilon,i}\rightarrow v_1^i$ in $C([0,T];\mathcal{C}^{-1/2-\delta/2})$, $u_1^{\varepsilon,i}\diamond u_1^{\varepsilon,j}\rightarrow v_2^{ij}$ in $C([0,T];\mathcal{C}^{-1-\delta/2})$, $u_1^{\varepsilon,i}\diamond u_2^{\varepsilon,j}\rightarrow v_3^{ij}$ in $C([0,T];\mathcal{C}^{-1/2-\delta/2})$, $u_2^{\varepsilon,i}\diamond u_2^{\varepsilon,j}\rightarrow v_4^{ij}$ in $C([0,T];\mathcal{C}^{-\delta})$, $\pi_{0,\diamond}(u_3^{\varepsilon,i}, u_1^{\varepsilon,j})\rightarrow v_5^{ij}$ in $C([0,T];\mathcal{C}^{-\delta})$, $\pi_{0,\diamond}(P^{ii_1}D_jK^{\varepsilon,j}, u_1^{\varepsilon,j_2})\rightarrow v_6^{ii_1jj_2}$ in $C([0,T];\mathcal{C}^{-\delta})$ and $\pi_{0,\diamond}(P^{ii_1}D_jK^{\varepsilon,i_1}, u_1^{\varepsilon,j_2})\rightarrow v_7^{ii_1jj_2}$ in $C([0,T];\mathcal{C}^{-\delta})$.
Here $$u^{\varepsilon,i}_1\diamond u^{\varepsilon,j}_1:=u^{\varepsilon,i}_1u^{\varepsilon,j}_1-C^{\varepsilon,ij}_0,$$
 $$u_1^{\varepsilon,i}\diamond u_2^{\varepsilon,j}:=u_1^{\varepsilon,i} u_2^{\varepsilon,j},$$
 $$u_2^{\varepsilon,i}\diamond u_2^{\varepsilon,j}:=u_2^{\varepsilon,i} u_2^{\varepsilon,j}-C_2^{\varepsilon,ij},$$
 $$\pi_{0,\diamond}(u_3^{\varepsilon,i}, u_1^{\varepsilon,j}):=\pi_{0}(u_3^{\varepsilon,i}, u_1^{\varepsilon,j})-C_1^{\varepsilon,ij},$$
 $$\pi_{0,\diamond}(P^{ii_1}D_jK^{\varepsilon,j}, u_1^{\varepsilon,j_2}):=\pi_{0}(P^{ii_1}D_jK^{\varepsilon,j}, u_1^{\varepsilon,j_2}),$$
  $$\pi_{0,\diamond}(P^{ii_1}D_jK^{\varepsilon,i_1}, u_1^{\varepsilon,j_2}):=\pi_{0}(P^{ii_1}D_jK^{\varepsilon,i_1}, u_1^{\varepsilon,j_2}),$$
 and $C^\varepsilon_0\in\mathbb{R}$ is defined in Section 3.3,
 $C_1^\varepsilon$ is defined in Section 3.3.1 and  $C_2^\varepsilon$  is defined in Appendix 4.2.
 Hence we obtain the following theorem:

\vskip.10in
\th{Theorem 3.8} Let $z\in (1/2,1/2+\delta_0)$ with $0<\delta_0<1/2$ and assume that $(\xi^\varepsilon)_{\varepsilon>0}$ is a family of smooth functions converging to $\xi$ as $\varepsilon\rightarrow0$. Let for $\varepsilon>0$ the function $u^\varepsilon$ be the unique maximal solution to the Cauchy problem
$$Lu^{\varepsilon,i}=\sum_{i_1=1}^3P^{ii_1}\xi^{\varepsilon,i_1}-\frac{1}{2}\sum_{i_1=1}^3P^{ii_1}(\sum_{j=1}^3D_j(u^{\varepsilon,i_1} u^{\varepsilon,j}))\quad u^\varepsilon(0)=Pu_0,$$
such that $u_4^\varepsilon$ defined as above belongs to $C((0,T_\varepsilon);\mathcal{C}^{1/2-\delta_0})$, where $u_0\in \mathcal{C}^{-z}$. Suppose that $\mathbb{Z}(\xi^\varepsilon)$ converges to $(v_1,v_2,v_3,v_4,v_5, v_6,v_7)$ in $\mathbb{X}$.  Then there exist $\tau=\tau(u_0,v_1,v_2,v_3,v_4,v_5,v_6,v_7)>0$ and $u\in C([0,\tau];\mathcal{C}^{-z})$ such that $$\sup_{t\in[0,\tau]}\|u^\varepsilon-u\|_{-z}\rightarrow0.$$
The limit $u$ depends only on $(u_0,v_i),i=1....,7$, and not on the approximating family.

\vskip.10in
\th{Remark 3.9} Indeed we can define the solution space as follows: $u-u_1\in \mathcal{D}_X^L$ if
$$u-u_1=u_2+u_3-\frac{1}{2}\int_0^tP_{t-s}P\sum_{j=1}^3D_j[\pi_<(\Phi',u_1^j)+\pi_<(\Phi'^{j},u_1)]ds+\Phi^{\sharp}$$
such that
$$\|\Phi^\sharp\|_{\star,1,L,T}:=\sup_{t\in[0,T]}t^{\frac{1-\eta+z}{2}}\|\Phi_t^\sharp\|_{1-\eta}+\sup_{t\in[0,T]}t^{\frac{\gamma+z}{2}}\|\Phi_t^\sharp\|_{\gamma}+\sup_{s,t\in[0,T]}s^{\frac{z+a}{2}}\frac{\|\Phi_t^\sharp-\Phi_s^\sharp\|_{a-2b}}{|t-s|^b}<\infty,$$
and
$$\|\Phi'\|_{\star,2,L,T}:=\sup_{t\in[0,T]}t^{\frac{2\gamma+z}{2}}\|\Phi_t'\|_{1/2-\kappa}+\sup_{s,t\in[0,T]}s^{\frac{z+a}{2}}\frac{\|\Phi'_t-\Phi'_s\|_{c-2d}}{|t-s|^d}<\infty.$$
Here $\eta,\gamma\in(0,1),a\geq2b,0<\kappa<1/2,c\geq 2d$. By a similar argument as in \cite{CC13},  if $u-u_1\in \mathcal{D}_X^L$ then the equation
 $$u-u_1=P_t(u_0-u_1(0))-\frac{1}{2}\int_0^tP_{t-s}P\sum_{j=1}^3D_j(u_1\diamond u_1^j+(u-u_1)\diamond u_1^j+u_1\diamond (u-u_1)^j+(u-u_1)\diamond(u-u_1)^jds$$
 can be well defined and by a fixed point argument we also obtain local existence and uniqueness of solutions. The calculations for this method are more complicated and we will not go into details here.

\subsection{Renormalisation}
In the following we use the notation $X$ to represent $u_1$,  $k_{1,...,n}:=\sum_{i=1}^nk_i$ and $$\hat{f}(k):= (2\pi)^{-\frac{3}{2}}\int_{\mathbb{T}^3} f(x)e^{\imath x\cdot k}dx$$ for $k\in\mathbb{Z}^3$. To simplify the arguments below, we assume that $\hat{\xi}(0)=0$ and restrict ourselves to the flow of $\int_{\mathbb{T}^3} u(x)dx=0$. Then we know that $X_t=\sum _{k\in\mathbb{Z}^3\backslash\{0\}}\hat{X}_t(k)e_k$ is a centered Gaussian process with covariance function given by
$$E[\hat{X}_t^i(k)\hat{X}_s^j(k')]=1_{k+k'=0}\sum_{i_1=1}^3\frac{e^{-|k|^2|t-s|}}{2|k|^2}\hat{P}^{ii_1}(k)\hat{P}^{ji_1}(k),$$
and $\hat{X}_t(0)=0$, where $e_k(x)=(2\pi)^{-3/2}e^{\imath x\cdot k},x\in\mathbb{T}^3$ and $\hat{P}^{ii_1}(k)=\delta_{ii_1}-\frac{k_ik_{i_1}}{|k|^2}$ for $k\in\mathbb{Z}^3\backslash\{0\}$. Let us take a smooth radial function $f$ with compact support such that $f(0)=1$. We regularize $X$ in the following way
$$X_t^{\varepsilon,i}=\int_{-\infty}^t\sum_{i_1=1}^3P^{ii_1}P_{t-s}\xi^{\varepsilon,i_1} ds$$ with $\xi^\varepsilon=\sum_{k\in\mathbb{Z}^3\backslash\{0\}}f(\varepsilon k)\hat{\xi}(k)$.
In this subsection we will prove that
there exist $v_1,v_2,v_3,v_4,v_5,v_6,v_7$ such that $\mathbb{Z}(\xi^\varepsilon)$ converges to $(v_1,v_2,...,v_7)$ in $\mathbb{X}$.

\vskip.10in

It is easy to obtain that there exists $v_1$ such that $u_1^\varepsilon\rightarrow v_1$ in $L^p(\Omega,P,C([0,T];\mathcal{C}^{-1/2-\delta/2}))$ for every $p\geq1$. The renormalisation of $u_1^{\varepsilon,i}\diamond u_1^{\varepsilon,j},i,j=1,2,3$ and the fact that there exists $v_2\in C([0,T];\mathcal{C}^{-1-\delta})$ such that $u_1^{\varepsilon,i}\diamond u_1^{\varepsilon,j}\rightarrow v_2^{ij}$ in $L^p(\Omega,P,C([0,T];\mathcal{C}^{-1-\delta}))$ for every $p\geq1$ can be easily obtained by using the Wick product (c.f.\cite{CC13}), where  $$C_0^{\varepsilon,ij}=(2\pi)^{-3}\sum_{i_1=1}^3\sum_{k\in \mathbb{Z}^3\backslash\{0\}}\frac{f(\varepsilon k)^2}{2|k|^2}\hat{P}^{ii_1}(k)\hat{P}^{ji_1}(k).$$ It is obvious that $C_0^{\varepsilon,ij}\rightarrow\infty$ as $\varepsilon\rightarrow0$. Here $u_1^\varepsilon$ and $u_1^{\varepsilon,i}\diamond u_1^{\varepsilon,j}$ correspond to $\includegraphics[height=0.5cm]{01.eps} $ and $\includegraphics[height=0.5cm]{02.eps} $ in Section 2 respectively.
By a similar argument as in the proof of Theorem 2.17 we could  conclude that  $u_1^{\varepsilon,i}\diamond u_2^{\varepsilon,j}\rightarrow v_3^{ij}$ in $C([0,T];\mathcal{C}^{-1/2-\delta})$, $u_2^{\varepsilon,i}\diamond u_2^{\varepsilon,j}\rightarrow v_4^{ij}$ in $C([0,T];\mathcal{C}^{-\delta})$. We could also use Fourier analysis to obtain it. Here for completeness of this method we calculate it in the appendix. For the terms including $\pi_0$ we cannot use a similar argument as in the proof of Theorem 2.17 to obtain the results since the definition of $\pi_0$ depends on the Fourier analysis. That is one of difference between these two approaches (see Remark 3.13).

 We first prove the following two lemmas for  later use, the first of which is inspired by [16, Lemma 10.14].

\vskip.10in
\th{Lemma 3.10} Let $0<l,m<d,l+m-d>0$. Then
$$\sum_{k_1,k_2\in \mathbb{Z}^d\backslash\{0\},k_1+k_2=k}\frac{1}{|k_1|^{l}|k_2|^{m}}\lesssim \frac{1}{|k|^{l+m-d}}.$$
\proof We have the following estimate:$$\aligned\sum_{k_1,k_2\in \mathbb{Z}^d\backslash\{0\},k_1+k_2=k}\frac{1}{|k_1|^{l}|k_2|^{m}}\lesssim&\sum_{k_1,k_2\in \mathbb{Z}^d\backslash\{0\},k_1+k_2=k,|k_1|\leq \frac{|k|}{2}}\frac{1}{|k_1|^{l}|k_2|^{m}}\\&+\sum_{k_1,k_2\in \mathbb{Z}^d\backslash\{0\},k_1+k_2=k,|k_2|\leq \frac{|k|}{2} }\frac{1}{|k_1|^{l}|k_2|^{m}}\\&+\sum_{k_1,k_2\in \mathbb{Z}^d\backslash\{0\},k_1+k_2=k,|k_1|>\frac{|k|}{2},|k_2|>\frac{|k|}{2}}\frac{1}{|k_1|^{l}|k_2|^{m}}.\endaligned$$
Since  $|k_1|\leq |k|/2$ implies that $|k_2|\geq |k|-|k_1|\geq |k|/2$, we obtain
$$\sum_{k_1,k_2\in \mathbb{Z}^d\backslash\{0\},k_1+k_2=k,|k_1|\leq \frac{|k|}{2}}\frac{1}{|k_1|^{l}|k_2|^{m}}\lesssim \sum_{k_1\in \mathbb{Z}^d\backslash\{0\},|k_1|\leq \frac{|k|}{2}}\frac{1}{|k_1|^{l}|k|^{m}}\lesssim  |k|^{-l-m+d}.$$
For the second term  a similar argument also yields the desired  estimate.
For the third term:  by $|k_2|\geq |k_1|-|k|$ and the triangle inequality, one has
$$|k_2|\geq \frac{1}{4} (|k_1|-|k|)+\frac{1-1/4}{2}|k|\geq\frac{1}{4}|k_1|,$$
which implies that $$\sum_{k_1,k_2\in \mathbb{Z}^d\backslash\{0\},k_1+k_2=k,|k_1|>\frac{|k|}{2},|k_2|>\frac{|k|}{2}}\frac{1}{|k_1|^{l}|k_2|^{m}}\lesssim |k|^{-l-m+d}.$$
Hence the result follows.$\hfill\Box$

\vskip.10in
  \th{Lemma 3.11} For any $0<\eta<1$, $i,j,l=1,2,3$ and for $t>0$ the following estimate holds:
  $$|e^{-|k_{12}|^2t}k_{12}^i\hat{P}^{jl}(k_{12})-e^{-|k_{2}|^2t}k_{2}^i\hat{P}^{jl}(k_{2})|\lesssim |k_1|^\eta|t|^{-(1-\eta)/2}.$$
  Here $\hat{P}^{ij}(x)=\delta_{ij}-\frac{x^i x^j}{|x|^2}$.

  \proof First we have the following bound:
  $$|e^{-|k_{12}|^2t}k_{12}\hat{P}(k_{12})-e^{-|k_{2}|^2t}k_{2}\hat{P}(k_{2})|\lesssim |t|^{-1/2}.$$
  Consider the function $F(x)=e^{-|x|^2t}x\hat{P}(x)$. Then it is easy to check that $|DF|$ is bounded, which implies that
   $$|e^{-|k_{12}|^2t}k_{12}\hat{P}(k_{12})-e^{-|k_{2}|^2t}k_{2}\hat{P}(k_{2})|\lesssim |k_1|.$$
   Thus, the result follows by the interpolation.$\hfill\Box$

\subsubsection{Renormalisation for $\pi_0(u_3^{\varepsilon,i_0},u_1^{\varepsilon,j_0})$}
Now we consider $\pi_0(u_{31}^{\varepsilon, i_0},u_1^{\varepsilon, j_0})$. The estimates for $\pi_0(u_3^{\varepsilon, i_0}-u_{31}^{\varepsilon, i_0},u_1^{\varepsilon, j_0})$ can be obtained similarly, where $Lu_{31}^{i_0}=-\frac{1}{2}\sum_{i_1=1}^3P^{i_0i_1}\sum_{j=1}^3D_j(u_2^{i_1}\diamond u_1^j).$ We have the following identity:
$$\pi_0(u_{31}^{\varepsilon, i_0i_1},u_1^{\varepsilon, j_0})(t)=\frac{1}{4}\sum_{i=1}^7I_t^i,$$
where
$$\aligned I_t^1&=(2\pi)^{-9/2}\sum_{k\in \mathbb{Z}^3\backslash\{0\}}\sum_{|i-j|\leq1}\sum_{k_{1234}=k}\sum_{i_1,i_2,i_3,j_1=1}^3\theta(2^{-i} k_{123})\theta(2^{-j}k_4)\int_0^tdse^{-|k_{123}|^2(t-s)}\int_0^s:\hat{X}_\sigma^{\varepsilon ,i_2}(k_1)\\&\hat{X}^{\varepsilon, i_3}_\sigma(k_2)\hat{X}_s^{\varepsilon,j_1}(k_3)\hat{X}^{\varepsilon, j_0}_t(k_4):
e^{-|k_{12}|^2(s-\sigma)}d\sigma \imath k_{12}^{i_3}\imath k_{123}^{j_1}\hat{P}^{i_1i_2}(k_{12})\hat{P}^{i_0i_1}(k_{123})e_k,\\ I_t^2+I_t^3&=(2\pi)^{-9/2}\sum_{k\in \mathbb{Z}^3\backslash\{0\}}\sum_{|i-j|\leq1}\sum_{k_{23}=k,k_1}\sum_{i_1,i_2,i_3,j_1=1}^3\theta(2^{-i} k_{123})\theta(2^{-j}k_1)\int_0^tdse^{-|k_{123}|^2(t-s)}\\&\int_0^s:\hat{X}^{\varepsilon,i_5}_\sigma(k_2)\hat{X}^{\varepsilon,j_1}_s(k_3)
:\frac{e^{-|k_1|^2(t-\sigma)}f(\varepsilon k_1)^2}{2|k_1|^2}\sum_{i_4=1}^3\hat{P}^{i_6i_4}(k_1)\hat{P}^{j_0i_4}(k_1)e^{-|k_{12}|^2(s-\sigma)}d\sigma\\& \imath k_{12}^{i_3}\imath k_{123}^{j_1}\hat{P}^{i_1i_2}(k_{12})\hat{P}^{i_0i_1}(k_{123})(1_{i_5=i_3,i_6=i_2}+1_{i_5=i_2,i_6=i_3})e_k,\\I_t^4&=(2\pi)^{-9/2}\sum_{k\in \mathbb{Z}^3\backslash\{0\}}\sum_{|i-j|\leq1}\sum_{k_{12}=k,k_3}\sum_{i_1,i_2,i_3,j_1=1}^3\theta(2^{-i} k_{123})\theta(2^{-j}k_3)\int_0^tdse^{-|k_{123}|^2(t-s)}\int_0^s:\hat{X}_\sigma^{\varepsilon,i_2}(k_1)
 \\&\hat{X}_\sigma^{\varepsilon,i_3}(k_2):\frac{e^{-|k_3|^2(t-s)}f(\varepsilon k_3)^2}{2|k_3|^2}\sum_{i_4=1}^3\hat{P}^{j_1i_4}(k_3)\hat{P}^{j_0i_4}(k_3)e^{-|k_{12}|^2(s-\sigma)}d\sigma \imath k_{12}^{i_3}\imath k_{123}^{j_1}\hat{P}^{i_1i_2}(k_{12})\hat{P}^{i_0i_1}(k_{123})e_k,\\ I_t^5+I_t^6&=(2\pi)^{-9/2}\sum_{k\in \mathbb{Z}^3\backslash\{0\}}\sum_{|i-j|\leq1}\sum_{k_{14}=k,k_2}\sum_{i_1,i_2,i_3,j_1=1}^3\theta(2^{-i} k_{1})\theta(2^{-j}k_4)\int_0^tdse^{-|k_{1}|^2(t-s)}\\&\int_0^s:\hat{X}^{\varepsilon,i_5}_\sigma(k_1)\hat{X}^{\varepsilon,j_0}_t(k_4):\frac{e^{-|k_2|^2(s-\sigma)}f(\varepsilon k_2)^2}{2|k_2|^2}\sum_{i_4=1}^3\hat{P}^{i_6i_4}(k_2)\hat{P}^{j_1i_4}(k_2)\\&e^{-|k_{12}|^2(s-\sigma)}d\sigma \imath k_{12}^{i_3}\imath k_{1}^{j_1}\hat{P}^{i_1i_2}(k_{12})\hat{P}^{i_0i_1}(k_{1})(1_{i_5=i_2,i_6=i_3}+1_{i_5=i_3,i_6=i_2})e_k,\endaligned$$ $$\aligned I_t^7&=(2\pi)^{-6}\sum_{|i-j|\leq1}\sum_{k_1,k_2}\sum_{i_1,i_2,i_3,j_1=1}^3\theta(2^{-i} k_2)\theta(2^{-j}k_2)\int_0^tdse^{-|k_2|^2(t-s)}\int_0^s\frac{f(\varepsilon k_1)^2f(\varepsilon k_2)^2}{4|k_1|^2|k_2|^2}\\&\sum_{i_4,i_5=1}^3\big(\hat{P}^{i_3i_4}(k_1)\hat{P}^{j_1i_4}(k_1)\hat{P}^{i_2i_5}(k_2)\hat{P}^{j_0i_5}(k_2)
+\hat{P}^{i_2i_4}(k_1)\hat{P}^{j_1i_4}(k_1)\hat{P}^{i_3i_5}(k_2)\hat{P}^{j_0i_5}(k_2)\big)\\&e^{-|k_{12}|^2(s-\sigma)-|k_1|^2(s-\sigma)-|k_2|^2(t-\sigma)}d\sigma \imath k_{12}^{i_3}\imath k_{2}^{j_1}\hat{P}^{i_1i_2}(k_{12})\hat{P}^{i_0i_1}(k_{2})].\endaligned$$
Here $I_t^2, I_t^3$ and $I_t^5, I_t^6$ correspond to the terms associated with each indicator function respectively.
To make it more readable we write each term corresponding to the tree notation in Section 2. $ \pi_0(u_{31}^{\varepsilon, i_0},u_1^{\varepsilon, j_0})$ corresponds to $\includegraphics[height=0.7cm]{08.eps} $ and $I_t^1, I_t^2, I_t^3, I_t^4, I_t^5, I_t^6, I_t^7$ correspond to the associated $\hat{\mathcal{W}}^{(\varepsilon,4)}, \hat{\mathcal{W}}^{(\varepsilon,2)}_4, \hat{\mathcal{W}}^{(\varepsilon,2)}_5, \hat{\mathcal{W}}^{(\varepsilon,2)}_3, \hat{\mathcal{W}}^{(\varepsilon,2)}_1, \hat{\mathcal{W}}^{(\varepsilon,2)}_2, \hat{\mathcal{W}}^{(\varepsilon,0)}$ in the proof of Theorem 2.17 respectively.

First we consider $I_t^7$: by simple calculations we have
$$\aligned I_t^7=&(2\pi)^{-6}\sum_{k_1,k_2}\sum_{i_1,i_2,i_3,j_1=1}^3\imath k_{12}^{i_3}\imath k_{2}^{j_1}\hat{P}^{i_1i_2}(k_{12})\hat{P}^{i_0i_1}(k_{2})\frac{f(\varepsilon k_1)^2f(\varepsilon k_2)^2}{4|k_1|^2|k_2|^2(|k_1|^2+|k_2|^2+|k_{12}|^2)}\\&\sum_{i_4,i_5=1}^3\big(\hat{P}^{i_3i_4}(k_1)\hat{P}^{j_1i_4}(k_1)\hat{P}^{i_2i_5}(k_2)\hat{P}^{j_0i_5}(k_2)
+\hat{P}^{i_2i_4}(k_1)\hat{P}^{j_1i_4}(k_1)\hat{P}^{i_3i_5}(k_2)\hat{P}^{j_0i_5}(k_2)\big)
\\&\bigg[\frac{1-e^{-2|k_2|^2t}}{2|k_2|^2}-\int_0^tdse^{-2|k_2|^2(t-s)}e^{-(|k_{12}|^2+|k_1|^2+|k_2|^2)s}\bigg]. \endaligned$$
Let $$C_{11}^{\varepsilon,i_0j_0}(t)=I_t^7$$
We could easily conclude that $C_{11}^{\varepsilon,i_0j_0}(t)\rightarrow\infty$, as $\varepsilon\rightarrow0$.

Similarly, we can also find  $C_{12}^\varepsilon$ for $u_3-u_{31}$. Define $C_1^\varepsilon=C_{11}^\varepsilon+C_{12}^\varepsilon$.

\textbf{Terms in the second chaos}: We come to $I_t^2$ and have the following calculations:
$$\aligned &E|\Delta_qI_t^2|^2\\\lesssim&\sum_{k\in \mathbb{Z}^3\backslash\{0\}}\sum_{|i-j|\leq1,|i'-j'|\leq1}\sum_{k_{23}=k,k_1,k_4}
\theta(2^{-i} k_{123})\theta(2^{-j}k_1)
\theta(2^{-i'} k_{234})\theta(2^{-j'}k_4)\theta(2^{-q}k)^2\\&\Pi_{i=1}^4\frac{f(\varepsilon k_i)^2}{|k_i|^2}\int_0^t\int_0^tdsd\bar{s}e^{-|k_{123}|^2(t-s)-|k_{234}|^2(t-\bar{s})}\int_0^s\int_0^{\bar{s}}d\sigma d\bar{\sigma}e^{-|k_1|^2(t-\sigma)-|k_4|^2(t-\bar{\sigma})}\\&e^{-(|k_{12}|^2(s-\sigma)+|k_{24}|^2(s-\bar{\sigma}))} |
k_{12}k_{123}k_{24}k_{234}|
\\\lesssim& \sum_{k\in \mathbb{Z}^3\backslash\{0\}}\sum_{|i-j|\leq1,|i'-j'|\leq1}\sum_{k_{23}=k,k_1,k_4}\theta(2^{-i} k_{123})\theta(2^{-j}k_1)\theta(2^{-i'} k_{234})\theta(2^{-j'}k_4)\theta(2^{-q}k)^2\\&\frac{t^\eta}{|k_2|^2|k_3|^2|k_1|^{4-\eta}|k_4|^{4-\eta}}\\\lesssim& \sum_{k\in \mathbb{Z}^3\backslash\{0\}}\sum_{q\lesssim i}2^{-(1-\eta-\epsilon)i}\sum_{q\lesssim i'}2^{-(1-\eta-\epsilon)i'}\sum_{k_{23}=k}\theta(2^{-q}k)^2\frac{t^\eta}{|k_2|^2|k_3|^{2}}\lesssim t^\eta2^{2q(\eta+2\epsilon)},\endaligned$$
where $\eta,\epsilon>0$ are small enough, we used  $\sup_{a\in\mathbb{R}} |a|^r \exp(-a^2)\leq C$ for $r\geq0$ in the second inequality and  Lemma 3.10 in the last inequality. Furthermore,  $q\lesssim i$ follows from $|k|\leq |k_{123}|+|k_1|\lesssim 2^{i}$ and similarly one gets $q\lesssim i'$. Also for $I^3_t$ we have a similar estimate.

Now we deal with $I_t^4=I_t^4-\tilde{I}_t^4+\tilde{I}_t^4-\sum_{i_1=1}^3u^{i_1}_2(t)C_3^{\varepsilon,i_1}(t)$ where
$$\aligned \tilde{I}_t^4=&(2\pi)^{-\frac{9}{2}}\sum_{k\in \mathbb{Z}^3\backslash\{0\}}\sum_{|i-j|\leq1}\sum_{k_{12}=k,k_3}\sum_{i_1,i_2,i_3,j_1=1}^3\theta(2^{-i} k_{123})\theta(2^{-j}k_3)\int_0^t:\hat{X}^{\varepsilon,i_2}_\sigma(k_1)
\hat{X}^{\varepsilon,i_3}_\sigma(k_2):e^{-|k_{12}|^2(t-\sigma)}\imath k_{12}^{i_3}\\&\hat{P}^{i_1i_2}(k_{12})e_k d\sigma\int_0^tdse^{-|k_{123}|^2(t-s)}\frac{e^{-|k_3|^2(t-s)}f(\varepsilon k_3)^2}{|k_3|^2}\sum_{i_4}\hat{P}^{j_1i_4}(k_3)\hat{P}^{j_0i_4}(k_3) \imath k_{123}^{j_1}\hat{P}^{i_0i_1}(k_{123}),\endaligned$$
and $$\aligned C_3^{\varepsilon,i_1}(t)=&(2\pi)^{-\frac{9}{2}}\sum_{|i-j|\leq1}\sum_{k_3}\sum_{j_1=1}^3\theta(2^{-i} k_{3})\theta(2^{-j}k_3)\int_0^tds\frac{e^{-2|k_3|^2(t-s)}f(\varepsilon k_3)^2}{|k_3|^2}\\& \sum_{i_4}\hat{P}^{j_1i_4}(k_3)\hat{P}^{j_0i_4}(k_3)\imath k_{3}^{j_1}\hat{P}^{i_0i_1}(k_{3})=0.\endaligned$$
Let $c_{k_{123},k_3}^{j_1}(t-s)=\sum_{i_1=1}^3e^{-|k_{123}|^2(t-s)}\frac{e^{-|k_3|^2(t-s)}f(\varepsilon k_3)^2}{|k_3|^2} |k_{123}^{j_1}\hat{P}^{i_0i_1}(k_{123})|$. Then we have for $\epsilon>0$ small enough,
$$\aligned& E|\Delta_q(I_t^4-\tilde{I}_t^4)|^2\\\lesssim &\sum_{k\in \mathbb{Z}^3\backslash\{0\}}\sum_{|i-j|\leq1,|i'-j'|\leq1}\sum_{k_{12}=k,k_3,k_4}\theta(2^{-q}k)^2\theta(2^{-i} k_{123})\theta(2^{-j}k_3)\theta(2^{-i'} k_{124})\theta(2^{-j'}k_4)\\& \int_0^tds\int_0^td\bar{s}\frac{1}{|k_1|^2|k_2|^2}\sum_{j_1,j_1'=1}^3
c^{j_1}_{k_{123},k_3}(t-s)c^{j_1'}_{k_{124},k_4}(t-\bar{s})\bigg[
\int_0^sd\sigma\int_0^{\bar{s}}d\bar{\sigma}(e^{-|k_{12}|^2(s-\sigma)}
-e^{-|k_{12}|^2(t-\sigma)})\\&(e^{-|k_{12}|^2(\bar{s}-\bar{\sigma})}-e^{-|k_{12}|^2(t-\bar{\sigma})})|k_{12}|^2
+\int_s^td\sigma\int^t_{\bar{s}}d\bar{\sigma}e^{-|k_{12}|^2(t-\sigma)-|k_{12}|^2(t-\bar{\sigma})}|k_{12}|^2
\bigg]
\\\lesssim &\sum_{k\in \mathbb{Z}^3\backslash\{0\}}\sum_{|i-j|\leq1,|i'-j'|\leq1}\sum_{k_{12}=k,k_3,k_4}\theta(2^{-q}k)^2\theta(2^{-i} k_{123})\theta(2^{-j}k_3)\theta(2^{-i'} k_{124})\theta(2^{-j'}k_4)\\
&\int_0^tds\int_0^td\bar{s}\frac{1}{|k_{12}||k_1|^2|k_2|^2}(t-s)^{1/4}(t-\bar{s})^{1/4}\sum_{j_1,j_1'=1}^3c_{k_{123},k_3}^{j_1}(t-s)c^{j_1'}_{k_{124},k_4}(t-\bar{s})\\\lesssim &\sum_{k\in \mathbb{Z}^3\backslash\{0\}}\sum_{|i-j|\leq1,|i'-j'|\leq1}\sum_{k_{12}=k,k_3,k_4}\theta(2^{-q}k)^2\theta(2^{-i} k_{123})\theta(2^{-j}k_3)\theta(2^{-i'} k_{124})\theta(2^{-j'}k_4)\\&\frac{t^{2\epsilon}}{|k_{12}||k_1|^2|k_2|^2|k_3|^2|k_4|^2(|k_{123}|^2+|k_3|^2)^{3/4-\epsilon}(|k_{124}|^2+|k_4|^2)^{3/4-\epsilon}}
\\\lesssim&t^{2\epsilon}\sum_{q\lesssim i}\sum_{q\lesssim i'}2^{-(i+i')(1/2-3\epsilon)}\sum_k\sum_{k_{12}=k}\theta(2^{-q}k)\frac{1}{|k_{12}||k_1|^2|k_2|^2}\\\lesssim&t^{2\epsilon}2^{-2q(1/2-3\epsilon)}\sum_k\sum_{k_{12}=k}\theta(2^{-q}k)\frac{1}{|k_{12}||k_1|^2|k_2|^2}
\lesssim t^{2\epsilon}2^{2q(3\epsilon)},\endaligned$$
where in the last inequality we used Lemma 3.10 and $q\lesssim i$ follows $|k|\leq |k_{123}|+|k_3|\lesssim 2^i$ and similarly one gets $q\lesssim i'$.
Moreover, by a similar argument as in the proof of Lemma 3.11 we obtain that for $\eta>\epsilon>0$ small enough
$$\aligned& E[|\Delta_q(\tilde{I}_t^4-\sum_{i_1=1}^3u^{\varepsilon,i_1}_2(t) C^{\varepsilon,i_1}_3(t))|^2]\\\lesssim&\sum_{k}\sum_{k_{12}=k}\frac{1}{|k_1|^2|k_2|^2|k_{12}|^2}\theta(2^{-q}k)^2
\bigg[\sum_{i_1,j_1=1}^3\sum_{|i-j|\leq1}\sum_{k_3}\theta(2^{-j}k_3)\int_0^t\frac{e^{-|k_3|^2(t-s)}f(\varepsilon k_3)^2}{|k_3|^2}\\&(\theta(2^{-i}k_{123})e^{-|k_{123}|^2(t-s)}k_{123}^{j_1}\hat{P}^{i_0i_1}(k_{123})-\theta(2^{-i}k_3)
e^{-|k_{3}|^2(t-s)}k^{j_1}_{3}\hat{P}^{i_0i_1}(k_{3}))ds\bigg]^2\\\lesssim & \sum_{k}\sum_{k_{12}=k}\frac{1}{|k_1|^2|k_2|^2|k_{12}|^{2-2\eta}}\theta(2^{-q}k)^2\bigg[\sum_{j=0}^\infty\sum_{k_3}
\theta(2^{-j}k_3)\int_0^t\frac{e^{-|k_3|^2(t-s)}}{|k_3|^2}(t-s)^{-(1-\eta)/2}ds\bigg]^2 \\\lesssim& t^{\eta-\epsilon}2^{q(2\eta)},\endaligned$$
where in the last inequality we used Lemma 3.10.

Now we consider $I_t^5=I_t^5-\tilde{I}_t^5+\tilde{I}_t^5-\bar{I}_t^5$, where
$$\aligned \tilde{I}_t^5=&(2\pi)^{-9/2}\sum_{k\in \mathbb{Z}^3\backslash\{0\}}\sum_{|i-j|\leq1}\sum_{k_{14}=k,k_2}\sum_{i_1,i_2,i_3,j_1=1}^3\theta(2^{-i} k_{1})\theta(2^{-j}k_4)\int_0^t:\hat{X}^{\varepsilon,i_2}_s(k_1)\hat{X}^{\varepsilon,j_0}_t(k_4):e^{-|k_{1}|^2(t-s)}\\&\imath k_{1}^{j_1}\hat{P}^{i_0i_1}(k_{1})e_k ds\int_0^sd\sigma e^{-|k_{12}|^2(s-\sigma)}\frac{e^{-|k_2|^2(s-\sigma)}f(\varepsilon k_2)^2}{|k_2|^2} \imath k_{12}^{i_3}\hat{P}^{i_1i_2}(k_{12})\sum_{i_4=1}^3\hat{P}^{i_3i_4}(k_2)\hat{P}^{j_1i_4}(k_2),\endaligned$$
and $$\aligned \bar{I}_t^5=&(2\pi)^{-9/2}\sum_{k\in \mathbb{Z}^3\backslash\{0\}}\sum_{|i-j|\leq1}\sum_{k_{14}=k,k_2}\sum_{i_1,i_2,i_3,j_1=1}^3\theta(2^{-i} k_{1})\theta(2^{-j}k_4)\int_0^t:\hat{X}^{\varepsilon,i_2}_s(k_1)\hat{X}^{\varepsilon,j_0}_t(k_4):e^{-|k_{1}|^2(t-s)}\endaligned$$
$$\aligned&\imath k_{1}^{j_1}\hat{P}^{i_0i_1}(k_{1})e_k ds\int_0^sd\sigma e^{-|k_{2}|^2(s-\sigma)}\frac{e^{-|k_2|^2(s-\sigma)}f(\varepsilon k_2)^2}{|k_2|^2} \imath k_{2}^{i_3}\hat{P}^{i_1i_2}(k_{2})\sum_{i_4=1}^3\hat{P}^{i_3i_4}(k_2)\hat{P}^{j_1i_4}(k_2)=0.\endaligned$$
Let $d_{k_{12},k_2}(s-\sigma)=\sum_{i_2,i_3=1}^3e^{-|k_{12}|^2(s-\sigma)}\frac{e^{-|k_2|^2(s-\sigma)}f(\varepsilon k_2)^2}{|k_2|^2} |k_{12}^{i_3}\hat{P}^{i_1i_2}(k_{12})|$.
Since by H\"{o}lder's inequality we obtain
$$\aligned &E(:\hat{X}^{\varepsilon,i_2}_s(k_1)\hat{X}^{\varepsilon,j_0}_t(k_4):-:\hat{X}^{\varepsilon,i_2}_\sigma(k_1)\hat{X}^{\varepsilon,j_0}_t(k_4):)
(\overline{:\hat{X}^{\varepsilon,i_2}_{\bar{s}}(k'_1)\hat{X}^{\varepsilon,j_0}_t(k'_4):
-:\hat{X}^{\varepsilon,i_2}_{\bar{\sigma}}(k'_1)\hat{X}^{\varepsilon,j_0}_t(k'_4):})
\\\lesssim&(1_{k_1=k'_1}1_{k_4=k'_4}+1_{k_1=k'_4}1_{k_4=k'_1})(\frac{1-e^{-|k_1|^2|s-\sigma|}}{|k_1|^2|k_4|^2})^{1/2}(\frac{1-e^{-|k'_1|^2|\bar{s}-\bar{\sigma}|}}{|k'_1|^2|k'_4|^2})^{1/2}
\\\lesssim&(1_{k_1=k'_1}1_{k_4=k'_4}+1_{k_1=k'_4}1_{k_4=k'_1})\frac{|k_1|^\eta|k'_1|^\eta}{|k_1||k'_1||k_4||k'_4|}|s-\sigma|^{\eta/2}|\bar{s}-\bar{\sigma}|^{\eta/2},\endaligned$$
it follows that for $\eta,\epsilon>0$ small enough
$$\aligned E|\Delta_q(I_t^5-\tilde{I}_t^5)|^2\lesssim &\sum_{k\in \mathbb{Z}^3\backslash\{0\}}\sum_{|i-j|\leq1,|i'-j'|\leq1}\sum_{k_{14}=k,k_3,k_2}\theta(2^{-q}k)^2\theta(2^{-i} k_{1})\theta(2^{-j}k_4)\theta(2^{-i'} k_{1})\theta(2^{-j'}k_4)\endaligned$$
$$\aligned& \int_0^tds\int_0^td\bar{s}\int_0^sd\sigma\int_0^{\bar{s}}d\bar{\sigma}e^{-|k_{1}|^2(t-s)}e^{-|k_{1}|^2(t-\bar{s})}|k_{1}|^2\frac{1}{|k_1|^{2-2\eta}|k_4|^2}\\&|s-\sigma|^{\eta/2}|\bar{s}-\bar{\sigma}|^{\eta/2}d_{k_{12},k_2}(s-\sigma)d_{k_{13},k_3}(\bar{s}-\bar{\sigma})
\\&+\sum_{k\in \mathbb{Z}^3\backslash\{0\}}\sum_{|i-j|\leq1,|i'-j'|\leq1}\sum_{k_{14}=k,k_3,k_2}\theta(2^{-q}k)^2\theta(2^{-i} k_{1})\theta(2^{-j}k_4)\theta(2^{-i'} k_{4})\theta(2^{-j'}k_1)\\& \int_0^tds\int_0^td\bar{s}\int_0^sd\sigma\int_0^{\bar{s}}d\bar{\sigma}e^{-|k_{1}|^2(t-s)}e^{-|k_{4}|^2(t-\bar{s})}|k_{1}||k_4|\frac{1}{|k_1|^{2-\eta}|k_4|^{2-\eta}}\\&|s-\sigma|^{\eta/2}|\bar{s}-\bar{\sigma}|^{\eta/2}d_{k_{12},k_2}(s-\sigma)d_{k_{34},k_3}(\bar{s}-\bar{\sigma})
\\\lesssim &\sum_{k\in \mathbb{Z}^3\backslash\{0\}}\sum_{|i-j|\leq1,|i'-j'|\leq1}\sum_{k_{14}=k}\theta(2^{-q}k)^2\theta(2^{-i} k_{1})\theta(2^{-j}k_4)\theta(2^{-i'} k_{1})\theta(2^{-j'}k_4)\\&(\frac{t^{\epsilon}}{|k_1|^{4-2\eta-2\epsilon}|k_4|^2}+\frac{t^{\epsilon}}{|k_1|^{3-\eta-\epsilon}|k_4|^{3-\eta-\epsilon}})
\\\lesssim&t^{\epsilon}\sum_k\sum_{k_{14}=k}\theta(2^{-q}k)\sum_{q\lesssim i}2^{-i}\frac{1}{|k_1|^{3-2\eta-2\epsilon}|k_4|^2}\\&+t^{\epsilon}\sum_k\sum_{k_{14}=k}\theta(2^{-q}k)\sum_{q\lesssim j}2^{-j\epsilon}\frac{1}{|k_1|^{3-\eta-2\epsilon}|k_4|^{3-\eta-\epsilon}}
\\\lesssim& t^{\epsilon}2^{q(2\epsilon+2\eta)},\endaligned$$where in the last inequality we used Lemma 3.10 and $q\lesssim i$ follows from $|k|\leq |k_1|+|k_4|\lesssim 2^i$.

Moreover, it follows by Lemma 3.11 that for $\eta,\epsilon>0$ small enough
$$\aligned E[|\Delta_q(\tilde{I}_t^5-\bar{I}_t^5)|^2]\lesssim&\sum_{k\in \mathbb{Z}^3\backslash\{0\}}\sum_{|i-j|\leq1,|i'-j'|\leq1}\sum_{k_{14}=k,k_3,k_2}\theta(2^{-q}k)^2\theta(2^{-i} k_{1})\theta(2^{-j}k_4)\theta(2^{-i'} k_{1})\theta(2^{-j'}k_4)
\\&\int_0^t\int_0^t|k_1|^{2+2\eta}e^{-|k_1|^2(t-s+t-\bar{s}+|s-\bar{s}|)}\frac{1}{|k_1|^2|k_4|^2}\int_0^s\frac{e^{-|k_2|^2(s-\sigma)}}{|k_2|^2}(s-
\sigma)^{-(1-\eta)/2}\\&\int_0^{\bar{s}}\frac{e^{-|k_3|^2(\bar{s}-\bar{\sigma})}}{|k_3|^2}(\bar{s}-\bar{\sigma})^{-(1-\eta)/2}dsd\bar{s}d\sigma d\bar{\sigma}\\&+\sum_{k\in \mathbb{Z}^3\backslash\{0\}}
\sum_{|i-j|\leq1,|i'-j'|\leq1}\sum_{k_{14}=k,k_3,k_2}\theta(2^{-q}k)^2\theta(2^{-i} k_{1})\theta(2^{-j}k_4)\theta(2^{-i'} k_{4})\theta(2^{-j'}k_1)\\&\int_0^t\int_0^t|k_1|^{1+2\eta}|k_4|e^{-2|k_1|^2(t-s)-2|k_4|^2(t-\bar{s})}\frac{1}{|k_1|^2|k_4|^2}\int_0^s\frac{e^{-|k_2|^2(s-\sigma)}}{|k_2|^2}(s-\sigma)^{-(1-\eta)/2}\\&\int_0^{\bar{s}}\frac{e^{-|k_3|^2(\bar{s}-\bar{\sigma})}}{|k_3|^2}(\bar{s}-\bar{\sigma})^{-(1-\eta)/2}dsd\bar{s}d\sigma d\bar{\sigma}\\\lesssim&t^{\epsilon}\sum_k\sum_{k_{14}=k}\theta(2^{-q}k)\sum_{q\lesssim i}2^{-i}\frac{1}{|k_1|^{3-2\eta-2\epsilon}|k_4|^2}\\&+t^{\epsilon}\sum_k\sum_{k_{14}=k}\theta(2^{-q}k)\sum_{q\lesssim j}2^{-j\epsilon}\frac{1}{|k_1|^{3-2\eta-2\epsilon}|k_4|^{3-\epsilon}}
\\\lesssim& t^{\epsilon}2^{q(2\epsilon+2\eta)},\endaligned$$
where in the last inequality we used Lemma 3.10 and $q\lesssim i$ follows from $|k|\leq |k_1|+|k_4|\lesssim 2^i$. Similar estimates can also be obtained for $I_t^6$.

\textbf{Terms in the fourth chaos}: Now for $I_t^1$ we have the following calculations:
$$\aligned &E[|\Delta_qI_t^1|^2]\\\lesssim&\sum_{k\in \mathbb{Z}^3\backslash\{0\}}\sum_{|i-j|\leq1,|i'-j'|\leq1}\sum_{k_{1234}=k,k'_{1234}=k}\theta(2^{-q}k)^2\theta(2^{-i} k_{123})\theta(2^{-j}k_4)\theta(2^{-i'} k'_{123})\theta(2^{-j'}k'_4)\\&(1_{k_1=k'_1,k_2=k'_2,k_3=k'_3,k_4=k'_4}+1_{k_1=k'_4,k_2=k'_2,k_3=k'_3,k_4=k'_1}
+1_{k_1=k'_1,k_2=k'_2,k_3=k'_4,k_4=k'_3}+1_{k_1=k'_3,k_2=k'_4,k_3=k'_1,k_4=k'_2}\\&+1_{k_1=k'_1,k_2=k'_3,k_3=k'_2,k_4=k'_4}+1_{k_1=k'_3,k_2=k'_2,k_3=k'_4,k_4=k'_1}
+1_{k_1=k'_4,k_2=k'_2,k_3=k'_1,k_4=k'_3})\\&
\int_0^tds\int_0^td\bar{s}e^{-|k_{123}|^2(t-s)-|k_{123}'|^2(t-\bar{s})}\int_0^s\int_0^{\bar{s}}\frac{1}{|k_1|^2|k_2|^2|k_3|^2|k_4|^2}e^{-|k_{12}|^2(s-\sigma)-|k_{12}'|^2
(\bar{s}-\bar{\sigma})}d\sigma d\bar{\sigma} |k_{12}k_{123}k_{12}'k_{123}'|\\=&E_t^1+E_t^2+E_t^3+E_t^4+E_t^5+E_t^6+E_t^7. \endaligned.$$
Here each $E_t^i$ corresponds to the term associated with each indicator function.

For $\epsilon,\eta>0$ small enough by Lemma 3.10 we have
$$\aligned E_t^1\lesssim&\sum_{k\in \mathbb{Z}^3\backslash\{0\}}\sum_{|i-j|\leq1,|i'-j'|\leq1}\sum_{k_{1234}=k}\frac{\theta(2^{-q}k)^2\theta(2^{-i} k_{123})\theta(2^{-j}k_4)\theta(2^{-i'} k_{123})\theta(2^{-j'}k_4)t^\eta}{|k_1|^2|k_2|^2|k_3|^2|k_4|^2|k_{12}|^2|k_{123}|^{2-2\eta}}
\\\lesssim& \sum_{k\in \mathbb{Z}^3\backslash\{0\}}\sum_{|i-j|\leq1,|i'-j'|\leq1}\sum_{k_{1234}=k}\theta(2^{-q}k)^2\theta(2^{-i} k_{123})\theta(2^{-j}k_4)\theta(2^{-i'} k_{123})\theta(2^{-j'}k_4)\frac{t^\eta}{|k_4|^2|k_{123}|^{4-2\eta-\epsilon}}
\\\lesssim& \sum_{k\in \mathbb{Z}^3\backslash\{0\}}\sum_{q\lesssim i}2^{-(2-2\eta-\epsilon)i}\theta(2^{-q}k)^2\frac{t^\eta}{|k|}
\lesssim2^{q(2\eta+\epsilon)}t^\eta,\endaligned$$
and $$\aligned E_t^2\lesssim&\sum_{k\in \mathbb{Z}^3\backslash\{0\}}\sum_{|i-j|\leq1,|i'-j'|\leq1}\sum_{k_{1234}=k}\frac{\theta(2^{-q}k)^2\theta(2^{-i} k_{123})\theta(2^{-j}k_4)\theta(2^{-i'} k_{234})\theta(2^{-j'}k_1)t^\eta}{|k_1|^2|k_2|^2|k_3|^2|k_4|^2|k_{12}||k_{24}||k_{123}|^{1-\eta}|k_{234}|^{1-\eta}}
\\\lesssim& \sum_{k\in \mathbb{Z}^3\backslash\{0\}}\sum_{k_{1234}=k}\frac{\theta(2^{-q}k)^2t^\eta2^{-q(2-2\eta)}}{|k_1|^{1+\eta}|k_2|^2|k_3|^2|k_4|^{1+\eta}|k_{12}||k_{24}||k_{123}|^{1-\eta}|k_{234}|^{1-\eta}}
\\\lesssim& \sum_{k\in \mathbb{Z}^3\backslash\{0\}}(\sum_{k_{1234}=k}\frac{\theta(2^{-q}k)^2t^\eta2^{-q(2-2\eta)}}{|k_1|^{1+\eta}|k_2|^2|k_3|^2|k_4|^{1+\eta}|k_{12}|^2|k_{123}|^{2-2\eta}})^{1/2}\\&(\sum_{k_{1234}=k}\frac{\theta(2^{-q}k)^2t^\eta2^{-q(2-2\eta)}}{|k_1|^{1+\eta}|k_2|^2|k_3|^2|k_4|^{1+\eta}|k_{24}|^2|k_{234}|^{2-2\eta}})^{1/2}
\\\lesssim& \sum_{k\in \mathbb{Z}^3\backslash\{0\}}2^{-(2-2\eta)q}\frac{t^\eta}{|k|}
\lesssim2^{q(2\eta)}t^\eta.\endaligned$$
By a similar argument we can also obtain the same bounds for $E_t^3, E_t^4$, $E_t^5$, $E_t^6$ and $E_t^7$, which implies that for $\epsilon,\eta>0$ small enough
$$E[|\Delta_qI_t^1|^2]\lesssim2^{q(2\eta+\epsilon)}t^\eta.$$
By a similar calculation as above we get that for $\eta,\epsilon, \gamma>0$ small enough
$$\aligned &E[|\Delta_q(\pi_{0,\diamond}(u_3^{\varepsilon_1,i_0},u_1^{\varepsilon_1,j_0})(t_1)-\pi_{0,\diamond}(u_3^{\varepsilon_1,i_0},u_1^{\varepsilon_1,j_0})(t_2)
-\pi_{0,\diamond}(u_3^{\varepsilon_2,i_0},u_1^{\varepsilon_2,j_0})(t_1)\\&+\pi_{0,\diamond}(u_3^{\varepsilon_2,i_0},u_1^{\varepsilon_2,j_0})(t_2))|^2]\\\lesssim & (\varepsilon_1^{2\gamma}+\varepsilon_2^{2\gamma})|t_1-t_2|^\eta2^{q(\epsilon+2\eta)},\endaligned$$
 which by Gaussian hypercontractivity and Lemma 3.1 implies that
 $$\aligned &E[\|\pi_{0,\diamond}(u_3^{\varepsilon_1,i_0},u_1^{\varepsilon_1,j_0})(t_1)-\pi_{0,\diamond}(u_3^{\varepsilon_1,i_0},u_1^{\varepsilon_1,j_0})(t_2)
-\pi_{0,\diamond}(u_3^{\varepsilon_2,i_0},u_1^{\varepsilon_2,j_0})(t_1)\\&+\pi_{0,\diamond}(u_3^{\varepsilon_2,i_0},u_1^{\varepsilon_2,j_0})(t_2)\|^p_{\mathcal{C}^{-\eta-\epsilon-3/p}}]\\\lesssim &E[\|\pi_{0,\diamond}(u_3^{\varepsilon_1,i_0},u_1^{\varepsilon_1,j_0})(t_1)-\pi_{0,\diamond}(u_3^{\varepsilon_1,i_0},u_1^{\varepsilon_1,j_0})(t_2)
-\pi_{0,\diamond}(u_3^{\varepsilon_2,i_0},u_1^{\varepsilon_2,j_0})(t_1)\\&+\pi_{0,\diamond}(u_3^{\varepsilon_2,i_0},u_1^{\varepsilon_2,j_0})(t_2)\|^p_{B^{-\eta-\epsilon}_{p,p}}]\\\lesssim& (\varepsilon_1^{p\gamma}+\varepsilon_2^{p\gamma})|t_1-t_2|^{p(\eta-\epsilon)/2},\endaligned\eqno(3.14)$$
 (see the proof of (4.2), (4.3)).  Thus, for every $i_0,j_0=1,2,3$ we choose $p$ large enough and deduce that  there exist $v_5^{i_0j_0}\in C([0,T],\mathcal{C}^{-\delta}), i_0,j_0=1,2,3,$ such that for $p>1$
$$\pi_{0,\diamond}(u_3^{\varepsilon,i_0},u_1^{\varepsilon,j_0})\rightarrow v_5^{i_0j_0}\textrm{ in } L^p(\Omega,P,C([0,T],\mathcal{C}^{-\delta})).$$
Here $\delta>0$ depending on $\eta,\epsilon, p$ can be chosen small enough.

\subsubsection{Renormalisation for $\pi_0(P^{i_1i_2}D_{j_0}K^{\varepsilon,j_0},u_1^{\varepsilon,j_1})$ and $\pi_0(P^{i_1i_2}D_{j_0}K^{\varepsilon,i_2},u_1^{\varepsilon,j_1})$}
In this subsection we consider  $\pi_0(P^{i_1i_2}D_{j_0}K^{\varepsilon,j_0},u_1^{\varepsilon,j_1})$ and $\pi_0(P^{i_1i_2}D_{j_0}K^{\varepsilon,i_2},u_1^{\varepsilon,j_1})$ for $i_1,i_2,j_0,j_1=1,2,3$ and have the following identity:
$$\aligned &\pi_0(P^{i_1i_2}D_{j_0}K^{\varepsilon,j_0},u_1^{\varepsilon,j_1})(t)\\=&(2\pi)^{-\frac{3}{2}}\sum_{k\in\mathbb{Z}^3\backslash\{0\}}\sum_{|i-j|\leq1}
\sum_{k_{12}=k}\theta(2^{-i}k_1)\theta(2^{-j}k_2)\int_0^te^{-(t-s)|k_1|^2}
\imath k_1^{j_0}:\hat{X}^{\varepsilon,j_0}_s(k_1)\hat{X}^{\varepsilon,j_1}_t(k_2):dse_k
\hat{P}^{i_1i_2}(k_1)\\&+(2\pi)^{-3}\sum_{|i-j|\leq1}\sum_{k_{1}}\theta(2^{-i}k_1)\theta(2^{-j}k_1)\int_0^te^{-2(t-s)|k_1|^2}\imath k_1^{j_0}\frac{f(\varepsilon k_1)^2}{2|k_1|^2}ds\hat{P}^{i_1i_2}(k_1)\sum_{j_2=1}^3\hat{P}^{j_0j_2}(k_1)\hat{P}^{j_1j_2}(k_1).\endaligned$$
Here $\pi_0(P^{i_1i_2}D_{j_0}K^{\varepsilon,j_0},u_1^{\varepsilon,j_1})$ corresponds to $\includegraphics[height=0.7cm]{05.eps} $ and the first term and the second term on the right hand side of the above equality correspond to the associated $\hat{\mathcal{W}}^{(\varepsilon,2)}, \hat{\mathcal{W}}^{(\varepsilon,0)}$ in the proof of Theorem 2.17 respectively.
It is easy to get that the second term on the right hand side of the above equality equals zero. It is straightforward to calculate for $\epsilon>0$ small enough:

$$\aligned &E|\Delta_q\pi_0(P^{i_1i_2}D_{j_0}K^{\varepsilon,j_0},u_1^{\varepsilon,j_1})|^2\\\lesssim&\sum_{k\in\mathbb{Z}^3
\backslash\{0\}}\sum_{|i-j|\leq1,|i'-j'|\leq1}\sum_{k_{12}=k}\theta(2^{-q}k)^2\theta(2^{-i}k_1)\theta(2^{-j}k_2)\theta(2^{-i'}k_1)\theta(2^{-j'}k_2)\\&
\bigg[\int_0^t\int_0^te^{-(t-s+t-\bar{s})|k_1|^2}|k_1|^2\frac{e^{-|k_1|^2|s-\bar{s}|}}{|k_1|^2|k_2|^2}dsd\bar{s}
\\&+
\int_0^t\int_0^te^{-2(t-s)|k_1|^2-2(t-\bar{s})|k_2|^2}|k_1||k_2|\frac{1}{|k_1|^2|k_2|^2}dsd\bar{s}\bigg]\endaligned$$
$$\aligned\lesssim &t^{\epsilon}\sum_k\sum_{q\lesssim i}\sum_{k_{12}=k}\theta(2^{-q}k)\theta(2^{-i}k_1)\frac{1}{|k_1|^{4-2\epsilon}|k_2|^2}\\&+t^\epsilon\sum_k\sum_{q\lesssim i}\sum_{k_{12}=k}\theta(2^{-q}k)\theta(2^{-j}k_2)\frac{1}{|k_1|^{3-2\epsilon}|k_2|^3}
\\\lesssim& t^{\epsilon}2^{2q\epsilon},\endaligned$$
where in the last inequality we used Lemma 3.10.
By a similar calculation we also get that for $\epsilon,\eta>0$, $\gamma>0$ small enough
$$\aligned &E[|\Delta_q(\pi_{0,\diamond}(P^{i_1i_2}D_{j_0}K^{\varepsilon,j_0},u_1^{\varepsilon,j_1})(t_1)
-\pi_{0,\diamond}(P^{i_1i_2}D_{j_0}K^{\varepsilon,j_0},u_1^{\varepsilon,j_1})(t_2)\\&-\pi_{0,\diamond}(P^{i_1i_2}D_{j_0}K^{\varepsilon,j_0},u_1^{\varepsilon,j_1}
)(t_1)+\pi_{0,\diamond}(P^{i_1i_2}D_{j_0}K^{\varepsilon,j_0},u_1^{\varepsilon,j_1})(t_2))|^2]\\\lesssim & (\varepsilon_1^{2\gamma}+\varepsilon_2^{2\gamma})|t_1-t_2|^\eta2^{q(\epsilon+2\eta)},\endaligned$$
 which by Gaussian hypercontractivity, Lemma 3.1 and similar arguments as for (3.14) implies that there exists $v_6^{i_1i_2j_0j_1}\in C([0,T];\mathcal{C}^{-\delta})$ for $i_1,i_2,j_0,j_1=1,2,3$ such that for $p>1$
$$\pi_{0,\diamond}(P^{i_1i_2}D_{j_0}K^{\varepsilon,j_0},u_1^{\varepsilon,j_1})\rightarrow v_6^{i_1i_2j_0j_1}\textrm{ in }L^p(\Omega, P, C([0,T];\mathcal{C}^{-\delta})).$$
Here $\delta>0$ depending on $\eta,\epsilon, p$ can be chosen small enough.
By a similar argument we also obtain that there exists $v_7^{i_1i_2j_0j_1}\in C([0,T];\mathcal{C}^{-\delta})$ for $i_1,i_2,j_0,j_1=1,2,3$ such that
$$\pi_{0,\diamond}(P^{i_1i_2}D_{j_0}K^{\varepsilon,i_2},u_1^{\varepsilon,j_1})\rightarrow v_7^{i_1i_2j_0j_1}\textrm{ in }L^p(\Omega, P,C([0,T];\mathcal{C}^{-\delta})).$$

Combining all the convergence results we obtained above and Theorem 3.8 we obtain local existence and uniqueness of the solutions to the 3D Navier-Stokes equation driven by space-time white noise.
\vskip.10in

\th{Theorem 3.12} Let $z\in (1/2,1/2+\delta_0)$ with $0<\delta_0<1/2$ and $u_0\in \mathcal{C}^{-z}$. Then there exists a unique local solution to
$$Lu^i=\sum_{i_1=1}^3P^{ii_1}\xi-\frac{1}{2}\sum_{i_1=1}^3P^{ii_1}(\sum_{j=1}^3D_j(u^{i_1}u^{j}))\quad u(0)=Pu_0,$$
in the following sense: For $\xi^\varepsilon=\sum_{k}f(\varepsilon k)\hat{\xi}(k)e_k$ with $f$ a smooth radial function with compact support satisfying $f(0)=1$ and for $\varepsilon>0$ consider the maximal unique solution $u^\varepsilon$ to the following equation, such that $u_4^\varepsilon$ defined above belongs to  $C((0,T^\varepsilon);\mathcal{C}^{1/2-\delta_0})$, $$Lu^{\varepsilon,i}=\sum_{i_1=1}^3P^{ii_1}\xi^\varepsilon-\frac{1}{2}\sum_{i_1=1}^3P^{ii_1}(\sum_{j=1}^3D_j(u^{\varepsilon,i_1}u^{\varepsilon,j})),\quad u^\varepsilon(0)=Pu_0.$$
Then there exists $u\in C([0,\tau);\mathcal{C}^{-z})$ and a sequence of random time $\tau_L$ converging to the explosion time $\tau$ of $u$ such that $$\sup_{t\in[0,\tau_L]}\|u^\varepsilon-u\|_{-z}\rightarrow^P0,.$$

\proof  By a similar argument as above we have that there exists some $\gamma>0$ and $u_1\in C([0,T];\mathcal{C}^{-1/2-\delta/2})$, $u_2\in C([0,T];\mathcal{C}^{-\delta})$, $u_3\in C([0,T];\mathcal{C}^{\frac{1}{2}-\delta})$  such that for every $p>0$ $$E\|u_1^{\varepsilon}-u_1\|_{C([0,T];\mathcal{C}^{-1/2-\delta/2})}^p\lesssim \varepsilon^{\gamma p},$$
$$E\|u_2^{\varepsilon}-u_2\|_{C([0,T];\mathcal{C}^{-\delta})}^p\lesssim \varepsilon^{\gamma p}.$$
$$E\|u_3^{\varepsilon}-u_3\|_{C([0,T];\mathcal{C}^{1/2-\delta})}^p\lesssim \varepsilon^{\gamma p}.$$
Then for  $\varepsilon_k=2^{-k}\rightarrow0$ and $\epsilon>0$
$$\sum_{k=1}^\infty P(\|u_1^{\varepsilon_k}-u_1\|_{C([0,T];\mathcal{C}^{-1/2-\delta/2})}>\epsilon)\lesssim \sum_{k=1}^\infty 2^{-k\gamma}/\epsilon<\infty,$$
which by the Borel-Cantelli Lemma implies that
$u_1^{\varepsilon_k,i}-u_1^i\rightarrow0$ in $C([0,T];\mathcal{C}^{-1/2-\delta/2})$ a.s., as $k\rightarrow\infty$.  The results for the other terms are similar. Thus we obtain that $\sup_{\varepsilon_k=2^{-k},k\in\mathbb{N}}\bar{C}_\xi^{\varepsilon_k}<\infty $ a.s.,  $T_0$ independent of $\varepsilon$, $u_4:=\lim_{k\rightarrow\infty}u_4^{\varepsilon_k}$ on $[0,T_0]$, $u=u_1+u_2+u_3+u_4$ as the solution to (3.1) on $[0,T_0]$ and
$$\sup_{t\in[0,T_0]}\|u^{\varepsilon_k}-u\|_{-z}\rightarrow0\quad a.s..$$
Now we can extend the solution to the maximal solution such that
$$\sup_{t\in[0,\tau)}\|{u}\|_{-z}=\infty.$$
Indeed, a similar argument as in the proof in Section 3.2 implies that there exists some $T_1(C(T_0))$ (for simplicity we assume $T_1\leq T_0$) such that for every $t^*\in[0,T_0]$
$$\sup_{t\in[t^*,t^*+T_1]}\big{[}(t-t^*)^{\delta+z+\kappa}\|\bar{u}^{ \varepsilon ,\sharp}\|_{1/2+\beta}+(t-t^*)^{\frac{\delta+z+\kappa}{2}}\|\bar{u}^{ \varepsilon,\sharp}(t)\|_{\delta}\big{]}\lesssim C(T_1,C^\varepsilon_\xi,C(T_0),\|u(t^*)\|_{-z}),$$
where $\bar{u}^{ \varepsilon }$ denotes the solution starting at $t^*$ with initial condition $\bar{u}^{ \varepsilon }(t^*)=u(t^*)$ and we can also define $\bar{u}^{\varepsilon,\sharp}$.  Here the only difference is that  $\bar{K}^{\varepsilon,i}$ satisfies the following equation
$$d\bar{K}^{\varepsilon,i}=(\Delta\bar{K}^{\varepsilon,i}+u_1^{\varepsilon,i})dt,\quad \bar{K}^{\varepsilon,i}(t^*)=0,$$
and by a similar argument as above we obtain that there exists some $\gamma>0$ such that for every $p>1$
$$E\sup_{r\in[0,T]}\|\pi_0(PD\int_r^{\cdot}P_{\cdot-s}u^\varepsilon_1ds,u^\varepsilon_1(\cdot))
-\pi_0(PD\int_r^{\cdot}P_{\cdot-s}u_1ds,u_1(\cdot))\|_{C([0,T];\mathcal{C}^{-\delta})}^p\lesssim\varepsilon^{p\gamma},$$
 which implies that a similar convergence also holds for  $\pi_0(PD\bar{K}^\varepsilon,u^\varepsilon_1)$ in this case. Here we omit superscripts for simplicity.

Therefore for $t^*=T_0-\frac{T_1(C(T_0))}{2}$ we obtain the following estimate
$$\aligned&\sup_{t\in[T_0,T_0+\frac{T_1}{2}]}(t^{\delta+z+\kappa}\|\bar{u}^{ \varepsilon ,\sharp}\|_{1/2+\beta}+t^{\frac{\delta+z+\kappa}{2}}\|\bar{u}^{ \varepsilon,\sharp}(t)\|_{\delta})\\\lesssim& \sup_{t\in[T_0,T_0+\frac{T_1}{2}]}((t-t^*)^{\delta+z+\kappa}\|\bar{u}^{ \varepsilon ,\sharp}(t)\|_{1/2+\beta}+(t-t^*)^{\frac{\delta+z+\kappa}{2}}\|\bar{u}^{ \varepsilon,\sharp}(t)\|_{\delta})\\\lesssim &C(T_1,C^\varepsilon_\xi,C(T_0),\|u_0\|_{-z}).\endaligned$$
Hence  by a similar argument as above we obtain the solution $u=\lim_{k\rightarrow\infty}\bar{u}^{\varepsilon_k}$ on $[T_0,T_0+\frac{T_1}{2}]$. Iterating  the above arguments  we get that there exist the explosion time $\tau>0$ and the maximal solution $u$ on $[0,\tau)$ such that
$$\sup_{t\in[0,\tau)}\|u(t)\|_{-z}=\infty.$$
In the following we prove $u^\varepsilon$ converges to $u$ before some random time. For $L\geq0$ define $\tau_L:=\inf\{ t:\|u(t)\|_{-z}\geq L\}\wedge L$. Then $\tau_L$ increases to $\tau$.  Also define $\tau^\varepsilon_L:=\inf\{ t:\|u^\varepsilon(t)\|_{-z}\geq L\}\wedge L$  and    $\rho_L^\varepsilon:=\inf\{ t:C^\varepsilon_\xi(t)\geq L\}.$ Then by the proof in Section 3.2 we obtain for any $L, L_1,L_2>0,$
$$\sup_{t\in[0, \rho^{\varepsilon}_{L_1}\wedge \tau_L\wedge \tau_{L_2}^\varepsilon]}\|u^{\varepsilon}-u\|_{-z}\rightarrow0\quad a.s..$$
Now we have for any $\epsilon>0$
$$P(\sup_{t\in[0,\tau_L]}\|u^{\varepsilon}-u\|_{-z}>\epsilon)\leq P(\sup_{t\in[0,\tau_L\wedge \rho_{L_1}^\varepsilon\wedge \tau_{L_2}^\varepsilon]}\|u^{\varepsilon}-u\|_{-z}>\epsilon)+P(\tau_L>\rho_{L_1}^\varepsilon)+P(\tau_L\wedge\rho_{L_1}^\varepsilon>\tau_{L_2}^\varepsilon).$$
Here the first term goes to zero by the above result, the second term goes to zero as $L_1$ goes to infinity and for $L_2>L+\epsilon$
$$P(\tau_L\wedge\rho_{L_1}^\varepsilon>\tau_{L_2}^\varepsilon)\leq P(\sup_{t\in[0,\tau_L\wedge \rho_{L_1}^\varepsilon\wedge \tau_{L_2}^\varepsilon]}\|u^{\varepsilon}-u\|_{-z}>\epsilon),$$
which goes to zero as $\varepsilon\rightarrow0$ by the above result. Thus the result follows.$\hfill\Box$
\vskip.10in

\th{Remark 3.13}  We used two different approaches and obtained the same results in Theorem 1.1 and Theorem 3.12. As we mentioned in the introduction from a “philosophical” perspective, the theory of regularity structures and  the paracontrolled distribution  are inspired by the theory
of controlled rough paths \cite{Lyo98}, \cite{Gub04}. The main difficulty for this problem lies in how to define multiplication for the unknowns. In the regularity structure theory we used an extension of the Taylor expansion and split the unknown into elements of different orders of homogeneity (i.e. regularity structure). Then it suffices to define the multiplications for these elements of different orders of homogeneity. In the paracontrolled distribution method  using Bony's paraproduct we split the unknown into good terms and bad terms ($\pi_<(\cdot,\cdot)$), where the singularity of the bad term is the same as the singularity of some functional of the Gaussian field. Then by using the commutator estimate it suffices to define the multiplication of some functionals of the Gaussian field.

 From the proof we see that the terms required to be renormalized in the two methods are similar: The terms not including the terms with $|\cdot|_{\mathfrak{s}}>0$ in the  theory of the regularity sturctures are the same as the associated terms in the paracontrolled distribution, while the terms including the terms with $|\cdot|_{\mathfrak{s}}>0$ ( like $\mathcal{I}_l(\mathcal{I}_k(\mathcal{I}(\Xi)\mathcal{I}(\Xi))
\mathcal{I}(\Xi))\mathcal{I}(\Xi)$ and $\mathcal{I}_k(\mathcal{I}(\Xi))\mathcal{I}(\Xi)$) are different from the terms in the paracontrolled distributions ($\pi_0(u_3,u_1)$ and $\pi_0(PDK,u_1)$).  In the theory of regularity structures  a distribution is  divided into  the elements of different orders of homogeneity. For example, the  terms of good regularity ( e.g. $u_3$) are split  into constants, polynomials and some other terms with  positive order ( e.g. $\mathcal{I}_l(\mathcal{I}_k(\mathcal{I}(\Xi)\mathcal{I}(\Xi))
\mathcal{I}(\Xi))$). In the paracontrolled distribution method  using Bony's paraproduct for these terms it is sufficient to define $\pi_0(\cdot,\cdot)$, which plays a similar role as the term of positive order  in the regularity structure theory.

\vskip.10in
\th{Acknowledgement.} We are very grateful to Professor Martin Hairer for giving us some hints to complete the paper and for pointing out an mistake in an early version of our paper, which helped us to improve the results of this paper. We would also like to thank Professor Michael R\"{o}ckner for his encouragement and suggestions for this work.
We are also grateful to the referees for their  comments that lead us to improve
the exposition and revise several details for the present version of our work.

\section{Appendix}

\subsection{Renormalisation for $u_1^\varepsilon u_2^\varepsilon$}
In this subsection we focus on $u_1^\varepsilon u_2^\varepsilon$ and prove that $u_1^{\varepsilon,i}\diamond u_2^{\varepsilon,j}\rightarrow v_3^{ij}$ in $C([0,T];\mathcal{C}^{-1/2-\delta})$ for $i,j=1,2,3$. Now we have the following identity: for $t\in[0,T]$, $i,j=1,2,3$
 $$\aligned u_1^{\varepsilon,j} u_2^{\varepsilon,i}(t)=&\frac{(2\pi)^{-3}}{2}\sum_{i_1,i_2=1}^3\sum_{k\in \mathbb{Z}^3\backslash\{0\}}\sum_{k_{123}=k}\int_0^te^{-|k_{12}|^2(t-s)} \imath k_{12}^{i_2}:\hat{X}^{\varepsilon,i_1}_s(k_1)\hat{X}^{\varepsilon,i_2}_s(k_2)\hat{X}^{\varepsilon,j}_t(k_3):ds \hat{ P}^{ii_1}(k_{12})e_k\\&+\frac{(2\pi)^{-3}}{2}\sum_{i_1,i_2,i_3=1}^3\sum_{k_1,k_2\in\mathbb{Z}^3\backslash\{0\}}\int_0^te^{-|k_{12}|^2(t-s)}\imath k_{12}^{i_2}\hat{X}^{\varepsilon,i_1}_s(k_1)\frac{e^{-|k_2|^2(t-s)}f(\varepsilon k_2)^2}{2|k_2|^2}ds\\&\hat{P}^{ii_1}(k_{12})\hat{P}^{i_2i_3}(k_{2})\hat{P}^{ji_3}(k_{2})e_{k_1}
 \\&+\frac{(2\pi)^{-3}}{2}\sum_{i_1,i_2,i_3=1}^3\sum_{k_1,k_2\in\mathbb{Z}^3\backslash\{0\}}\int_0^te^{-|k_{12}|^2(t-s)}\imath k_{12}^{i_2}\hat{X}^{\varepsilon,i_2}_s(k_2)\frac{e^{-|k_1|^2(t-s)}f(\varepsilon k_1)^2}{2|k_1|^2}ds\\&\hat{P}^{ii_1}(k_{12})\hat{P}^{i_2i_3}(k_{1})\hat{P}^{ji_3}(k_{2})e_{k_2}\\=&I_t^1+I_t^2+I_t^3.\endaligned$$
To make it more readable we write each term corresponding to the tree notation in Section 2: $ u_1^{\varepsilon,j} u_2^{\varepsilon,i}$ corresponds to $\includegraphics[height=0.7cm]{06.eps} $ and $I_t^1, I_t^2, I_t^3$ correspond to the associated $\hat{\mathcal{W}}^{(\varepsilon,3)}, \hat{\mathcal{W}}^{(\varepsilon,1)}_2, \hat{\mathcal{W}}^{(\varepsilon,1)}_1$ in the proof of Theorem 2.16 respectively.

\no\textbf{Term in the first chaos:}
First, we consider $I_t^2$. We have $$I_t^2=I_t^2-\tilde{I}_t^2+\tilde{I}_t^2-\sum_{i_1=1}^3X_t^{\varepsilon,i_1} C^{\varepsilon,i_1}_t,$$
 where $$\aligned\tilde{I}_t^2=&\frac{(2\pi)^{-3}}{2}\sum_{i_1,i_2,i_3=1}^3\sum_{k_1,k_2\in\mathbb{Z}^3\backslash\{0\}}\hat{X}^{\varepsilon,i_1}_t(k_1)e_{k_1}
 \int_0^te^{-|k_{12}|^2(t-s)}\imath k_{12}^{i_2}\frac{e^{-|k_2|^2(t-s)}f(\varepsilon k_2)^2}{2|k_2|^2}ds\\&\hat{P}^{ii_1}(k_{12})\hat{P}^{i_2i_3}(k_{2})\hat{P}^{ji_3}(k_{2}),\endaligned$$
 and $$C^{\varepsilon,i_1}_t=\frac{(2\pi)^{-3}}{2}\sum_{i_2,i_3=1}^3\sum_{k_2\in\mathbb{Z}^3\backslash\{0\}}\int_0^te^{-2|k_2|^2(t-s)}\imath k_2^{i_2}\frac{f(\varepsilon k_2)^2}{2|k_2|^2}\hat{P}^{ii_1}(k_2)\hat{P}^{i_2i_3}(k_{2})\hat{P}^{ji_3}(k_{2})ds=0.$$
A straightforward calculation yields that for $\eta>0$ small enough
$$\aligned E[|\Delta_q(I_t^2-\tilde{I}_t^2)|^2]\lesssim& E\bigg[\bigg|\sum_{i_1,i_2,i_3=1}^3\int_0^t\sum_{k_1}\theta(2^{-q}k_1)e_{k_1}a_{k_1}^{i_1i_2i_3}(t-s)(\hat{X}_s^{\varepsilon,i_1}(k_1)
-\hat{X}_t^{\varepsilon,i_1}(k_1))ds\bigg|^2\bigg]
\\\lesssim &\sum_{i_1,i_2,i_3=1}^3\sum_{i_1',i_2',i_3'=1}^3\int_0^t\int_0^tdsd\bar{s}\sum_{k_1,k'_1}\theta(2^{-q}k_1)\theta(2^{-q}k'_1)
|a_{k_1}^{i_1i_2i_3}(t-s)a_{k'_1}^{i_1'i_2'i_3'}(t-\bar{s})|\\&E\big|(\hat{X}_s^{\varepsilon,i_1}
(k_1)-\hat{X}_t^{\varepsilon,i_1}(k_1))(\overline{\hat{X}_{\bar{s}}^{\varepsilon,i_1'}(k'_1)-\hat{X}_t^{\varepsilon,i_1'}(k'_1)})\big|\endaligned$$
$$\aligned\lesssim &\sum_{k_1}\theta(2^{-q}k_1)^2\frac{f(\varepsilon k_1)^2}{|k_1|^{2(1-\eta)}}\bigg(\int_0^t|t-s|^{\eta/2}|a_{k_1}^{i_1i_2i_3}(t-s)|ds\bigg)^2.\endaligned$$
Here $$a_{k_1}^{i_1i_2i_3}(t-s)=\sum_{k_2}e^{-|k_{12}|^2(t-s)}k_{12}^{i_2}\frac{e^{-|k_2|^2(t-s)}f(\varepsilon k_2)^2}{|k_2|^2}\hat{P}^{ii_1}(k_{12})\hat{P}^{i_2i_3}(k_{12})\hat{P}^{ji_3}(k_{12}),$$
and in the third inequality we used that for $\eta>0$ small enough
$$\aligned &E|(\hat{X}_s^{\varepsilon,i_1}
(k_1)-\hat{X}_t^{\varepsilon,i_1}(k_1))(\overline{\hat{X}_{\bar{s}}^{\varepsilon,i_1'}(k'_1)-\hat{X}_t^{\varepsilon,i_1'}(k'_1)})|
\\\leq& 1_{k_1=k_1'}(E|(\hat{X}_s^{\varepsilon,i_1}
(k_1)-\hat{X}_t^{\varepsilon,i_1}(k_1))|^2)^{1/2}(E|(\overline{\hat{X}_{\bar{s}}^{\varepsilon,i_1'}(k'_1)-\hat{X}_t^{\varepsilon,i_1'}(k'_1)})|^2)^{1/2}
\\\lesssim& 1_{k_1=k_1'}(\frac{f(\varepsilon k_1)^2}{|k_1|^2}(1-e^{-|k_1|^2(t-s)}))^{1/2}(\frac{f(\varepsilon k'_1)^2}{|k'_1|^2}(1-e^{-|k'_1|^2(t-\bar{s})}))^{1/2}\\\lesssim&
\frac{f(\varepsilon k_1)^2}{|k_1|^2}|k_1|^{2\eta}|t-s|^{\eta/2}|t-\bar{s}|^{\eta/2}.\endaligned$$
Since  $\sup_{a\in\mathbb{R}} |a|^r \exp(-a^2)\leq C$ for $r\geq0$ implies that for $\eta>\epsilon>0$, $\epsilon$ small enough $|a_{k_1}^{i_1i_2i_3}(t-s)|\lesssim |t-s|^{-1-\epsilon/2}\sum_{k_2}\frac{1}{|k_2|^{3+\epsilon}}$, it follows that
$$\int_0^t|t-s|^{\eta/2}|a_{k_1}^{i_1i_2i_3}(t-s)|ds\lesssim \int_0^t|t-s|^{\eta/2-1-\epsilon/2}ds\sum_{k_2}\frac{1}{|k_2|^{3+\epsilon}}\lesssim t^{(\eta-\epsilon)/2},$$
which implies that
$$E[|\Delta_q(I_t^2-\tilde{I}_t^2)|^2]\lesssim 2^{q(1+2\eta)}t^{\eta-\epsilon}.$$

Moreover, by Lemma 3.11 we deduce that for $\epsilon>0$ small enough
$$\aligned &E[|\Delta_q(\tilde{I}_t^2-\sum_{i_1=1}^3X_t^{\varepsilon,i_1} C^{\varepsilon,i_1}_t)|^2]\\\lesssim &\sum_{k_1}\frac{f(\varepsilon k_1)^2}{2|k_1|^2}\theta(2^{-q}k_1)^2\bigg[\sum_{i_1,i_2,i_3=1}^3\sum_{k_2}\int_0^t\frac{e^{-|k_2|^2(t-s)}f(\varepsilon k_2)^2}{|k_2|^2}\\&\big(e^{-|k_{12}|^2(t-s)}k_{12}^{i_2}\hat{P}^{ii_1}(k_{12})\hat{P}^{i_2i_3}(k_{2})\hat{P}^{ji_3}(k_{2})-e^{-|k_{2}|^2(t-s)}k_{2}^{i_2}
\hat{P}^{ii_1}(k_{2})\hat{P}^{i_2i_3}(k_{2})\hat{P}^{ji_3}(k_{2})\big)ds\bigg]^2\endaligned$$
$$\aligned\lesssim & \sum_{k_1}\frac{f(\varepsilon k_1)^2}{|k_1|^{2-2\eta}}\theta(2^{-q}k_1)^2\bigg(\sum_{k_2}\int_0^t\frac{e^{-|k_2|^2(t-s)}f(\varepsilon k_2)^2}{|k_2|^2}(t-s)^{-(1-\eta)/2}ds\bigg)^2 \\\lesssim& t^{\eta-\epsilon}2^{q(1+2\eta)},\endaligned\eqno(4.1)$$
holds uniformly over $\varepsilon\in(0,1)$, which is the desired bound for $I_t^2$. Here in the third inequality we also used  $\sup_{a\in\mathbb{R}} |a|^r \exp(-a^2)\leq C$ for $r\geq0$.

Similarly, we obtain that
$$\aligned &E[|\Delta_qI_t^3|^2]\lesssim t^{\eta-\epsilon}2^{q(1+2\eta)}.\endaligned$$

\textbf{Term in the third chaos:}
Now we focus on the bounds for $I_t^1$.
Let $b_{k_{12}}^{i_1,i_2}(t-s)=e^{-|k_{12}|^2(t-s)} k_{12}^{i_2}\hat{P}^{ii_1}(k_{12}).$ We obtain the following inequalities:
$$\aligned &E|\Delta_q I_t^1|^2\\\lesssim&2\sum_{i_1,i_2=1}^3\sum_{i_1',i_2'=1}^3\sum_k\theta(2^{-q}k)\sum_{k_{123}=k}\Pi_{i=1}^3\frac{f(\varepsilon k_i)^2}{|k_i|^2}\int_0^t\int_0^te^{-(|k_1|^2+|k_2|^2)|s-\bar{s}|}|b_{k_{12}}^{i_1,i_2}(t-s)b_{k_{12}}^{i_1',i_2'}(t-\bar{s})|dsd\bar{s}
\\&+2\sum_{i_1,i_2=1}^3\sum_{i_1',i_2'=1}^3\sum_k\theta(2^{-q}k)\sum_{k_{123}=k}\Pi_{i=1}^3\frac{f(\varepsilon k_i)^2}{|k_i|^2}\int_0^t\int_0^te^{-|k_2|^2|s-\bar{s}|-|k_1|^2(t-s)-|k_3|^2(t-\bar{s})}\\&|b_{k_{12}}^{i_1,i_2}(t-s)b_{k_{32}}^{i_1',i_2'}(t-\bar{s})|dsd\bar{s}
\\:=&J_t^1+J_t^2.\endaligned$$
Since $|b_{k_{12}}^{i_1,i_2}(t-s)|\lesssim \frac{1}{|k_{12}|^{1-\eta}(t-s)^{1-\eta/2}}$ it follows by Lemma 3.10 that for $\eta>0$ small enough
$$\aligned J_t^1\lesssim &\sum_k\theta(2^{-q}k)\sum_{k_{123}=k}\Pi_{i=1}^3\frac{1}{|k_i|^2}\frac{t^\eta}{|k_{12}|^{2-2\eta}}\\\lesssim &\sum_k\theta(2^{-q}k)\sum_{k_{123}=k}\frac{t^\eta}{|k_3|^{2}|k_{12}|^{3-2\eta}}
\\\lesssim &t^\eta2^{q(1+2\eta)},\endaligned$$
and
$$\aligned J_t^2\lesssim &\sum_k\theta(2^{-q}k)\sum_{k_{123}=k}\frac{t^\eta}{|k_1|^{2}|k_{2}|^{2}|k_3|^{2}|k_{12}|^{1-\eta}|k_{32}|^{1-\eta}}\\\lesssim & \sum_k\theta(2^{-q}k)(\sum_{k_{123}=k}\frac{t^\eta}{|k_1|^{2}|k_2|^{2}|k_3|^{2}|k_{12}|^{2-2\eta}})^{1/2}(\sum_{k_{123}=k}\frac{t^\eta}{|k_1|^{2}|k_2|^{2}|k_3|^{2}|k_{32}|^{2-2\eta}})^{1/2}
\\\lesssim& t^\eta2^{q(1+2\eta)},\endaligned$$
which yield the desired estimate for $I_t^1$.
By a similar calculation we also obtain that for $\eta>\epsilon>0$, $\gamma>0$ small enough,
$$E[|\Delta_q(u_2^{\varepsilon_1,i}u_1^{\varepsilon_1,j}(t_1)-u_2^{\varepsilon_1,i}u_1^{\varepsilon_1,j}(t_2)-u_2^{\varepsilon_2,i}u_1^{\varepsilon_1,j}
(t_1)+u_2^{\varepsilon_2,i}u_1^{\varepsilon_1,j}(t_2))|^2]\lesssim (\varepsilon_1^{2\gamma}+\varepsilon_2^{2\gamma})|t_1-t_2|^{\eta-\epsilon}2^{q(1+2\eta)},\eqno(4.2)$$
 which by Gaussian hypercontractivity and Lemma 3.1 implies that
 $$\aligned &E[\|(u_2^{\varepsilon_1,i}u_1^{\varepsilon_1,j}(t_1)-u_2^{\varepsilon_1,i}u_1^{\varepsilon_1,j}(t_2)-u_2^{\varepsilon_2,i}u_1^{\varepsilon_1,j}
(t_1)+u_2^{\varepsilon_2,i}u_1^{\varepsilon_1,j}(t_2))\|^p_{\mathcal{C}^{-1/2-\eta-\epsilon-3/p}}]\\\lesssim &E[\|(u_2^{\varepsilon_1,i}u_1^{\varepsilon_1,j}(t_1)-u_2^{\varepsilon_1,i}u_1^{\varepsilon_1,j}(t_2)-u_2^{\varepsilon_2,i}u_1^{\varepsilon_1,j}
(t_1)+u_2^{\varepsilon_2,i}u_1^{\varepsilon_1,j}(t_2))\|^p_{B^{-1/2-\eta-\epsilon}_{p,p}}]\\\lesssim& (\varepsilon_1^{p\gamma}+\varepsilon_2^{p\gamma})|t_1-t_2|^{p(\eta-\epsilon)/2}.\endaligned\eqno(4.3)$$
  Thus, for every $i,j=1,2,3$ we choose $p$ large enough and deduce that there exists $v_3^{ij}\in C([0,T];\mathcal{C}^{-1/2-\delta/2})$ such that
$$u_2^{\varepsilon,i} \diamond u_1^{\varepsilon,j}\rightarrow v_3^{ij}\textrm{ in } L^p(\Omega,P,C([0,T];\mathcal{C}^{-1/2-\delta/2})).$$
Here $\delta>0$ depending on $\eta,\epsilon,p$ can be chosen small enough. For the proof of  (4.2)  we only calculate the corresponding  term as in (4.1) and the other terms can be obtained similarly. It is straightforward to calculate that for $0\leq t_1<t_2\leq T$
$$\aligned &E[|\Delta_q(\tilde{I}_{t_1}^2-\sum_{i_1=1}^3X_{t_1}^{\varepsilon,i_1} C^{\varepsilon,i_1}_{t_1}-\tilde{I}_{t_2}^2+\sum_{i_1=1}^3X_{t_2}^{\varepsilon,i_1} C^{\varepsilon,i_1}_{t_2})|^2]\\\lesssim&E\bigg|\sum_{i_1,i_2,i_3=1}^3\sum_{k_1}\hat{X}_{t_1}^{\varepsilon,i_1}(k_1)\theta(2^{-q}k_1)e_{k_1}
\bigg[\sum_{k_2}\int_0^{t_1}\frac{e^{-|k_2|^2(t_1-s)}f(\varepsilon k_2)^2}{|k_2|^2}\bigg(e^{-|k_{12}|^2(t_1-s)}k_{12}^{i_2}\hat{P}^{ii_1}(k_{12})\\&\hat{P}^{i_2i_3}(k_{2})\hat{P}^{ji_3}(k_{2})
-e^{-|k_{2}|^2(t_1-s)}k_{2}^{i_2}\hat{P}^{ii_1}(k_{2})\hat{P}^{i_2i_3}(k_{2})\hat{P}^{ji_3}(k_{2})\bigg)ds
-\sum_{k_2}\int_0^{t_2}\frac{e^{-|k_2|^2(t_2-s)}f(\varepsilon k_2)^2}{|k_2|^2}\endaligned$$
$$\aligned&\bigg(e^{-|k_{12}|^2(t_2-s)}k_{12}^{i_2}\hat{P}^{ii_1}(k_{12})\hat{P}^{i_2i_3}(k_{2})\hat{P}^{ji_3}(k_{2})
-e^{-|k_{2}|^2(t_2-s)}k_{2}^{i_2}\hat{P}^{ii_1}(k_{2})\hat{P}^{i_2i_3}(k_{2})\hat{P}^{ji_3}(k_{2})\bigg)ds\bigg]\bigg|^2
\\&+E\bigg|\sum_{i_1,i_2,i_3=1}^3\sum_{k_1}(\hat{X}_{t_1}^{\varepsilon,i_1}(k_1)-\hat{X}_{t_2}^{\varepsilon,i_1}(k_1))
\theta(2^{-q}k_1)e_{k_1}\int_0^{t_2}\frac{e^{-|k_2|^2(t_2-s)}f(\varepsilon k_2)^2}{|k_2|^2}\\&\bigg(e^{-|k_{12}|^2(t_2-s)}k_{12}^{i_2}\hat{P}^{ii_1}(k_{12})\hat{P}^{i_2i_3}(k_{2})\hat{P}^{ji_3}(k_{2})
-e^{-|k_{2}|^2(t_2-s)}k_{2}^{i_2}\hat{P}^{ii_1}(k_{2})\hat{P}^{i_2i_3}(k_{2})\hat{P}^{ji_3}(k_{2})\bigg)ds\bigg|^2
\\\lesssim&L_t^1+L_t^2+L_t^3+L_t^4,\endaligned$$
where
$$\aligned L_t^1=&\sum_{k_1}\sum_{i_1,i_2=1}^3\frac{1}{|k_1|^2}\theta(2^{-q}k_1)^2\bigg[\sum_{k_2}\int_0^{t_1}\frac{e^{-|k_2|^2(t_1-s)}(1-e^{-|k_2|^2(t_2-t_1)})f(\varepsilon k_2)^2}{|k_2|^2}\\&\big(e^{-|k_{12}|^2(t_1-s)}k_{12}^{i_2}\hat{P}^{ii_1}(k_{12})-e^{-|k_{2}|^2(t_1-s)}k_{2}^{i_2}\hat{P}^{ii_1}(k_{2})\big)ds\bigg]^2
\\L_t^2=&\sum_{k_1}\sum_{i_1,i_2=1}^3\frac{1}{|k_1|^2}\theta(2^{-q}k_1)^2\bigg[\sum_{k_2}\int_0^{t_1}\frac{e^{-|k_2|^2(t_2-s)}f(\varepsilon k_2)^2}{|k_2|^2}\bigg(e^{-|k_{12}|^2(t_1-s)}k_{12}^{i_2}\hat{P}^{ii_1}(k_{12})\\&-e^{-|k_{2}|^2(t_1-s)}k_{2}^{i_2}\hat{P}^{ii_1}(k_{2})
-e^{-|k_{12}|^2(t_2-s)}k_{12}^{i_2}
\hat{P}^{ii_1}(k_{12})+e^{-|k_{2}|^2(t_2-s)}k_{2}^{i_2}\hat{P}^{ii_1}(k_{2})\bigg)ds\bigg]^2\endaligned$$
$$\aligned L_t^3=&\sum_{k_1}\sum_{i_1,i_2=1}^3\frac{1}{|k_1|^2}\theta(2^{-q}k_1)^2
\bigg[\sum_{k_2}\int_{t_1}^{t_2}\frac{e^{-|k_2|^2(t_2-s)}f(\varepsilon k_2)^2}{|k_2|^2}\bigg(e^{-|k_{12}|^2(t_2-s)}k_{12}^{i_2}\hat{P}^{ii_1}(k_{12})\\&-e^{-|k_{2}|^2(t_2-s)}k_{2}^{i_2}\hat{P}^{ii_1}(k_{2})\bigg)ds\bigg]^2
\\L_t^4=&E\bigg|\sum_{k_1}\sum_{i_1,i_2=1}^3(
\hat{X}^{\varepsilon,i_1}_{t_1}(k_1)-\hat{X}^{\varepsilon,i_1}_{t_2}(k_1))\sum_{k_2}\theta(2^{-q}k_1)e_{k_1}\int_0^{t_2}\frac{e^{-|k_2|^2(t_2-s)}f(\varepsilon k_2)^2}{|k_2|^2}\\&\big(e^{-|k_{12}|^2(t_2-s)}k_{12}^{i_2}\hat{P}^{ii_1}(k_{12})-e^{-|k_{2}|^2(t_2-s)}k_{2}^{i_2}\hat{P}^{ii_1}(k_{2})\big)ds\bigg|^2.\endaligned$$
It is easy to deduce the desired estimates for $L_t^1,L_t^3,L_t^4$ as for (4.1) and it is sufficient to consider $L_t^2$: for some $0<\beta_0<1/2,\eta>0$ small enough, by Lemma 3.11 and  interpolation we have
$$\aligned L_t^2 \lesssim& \sum_{k_1}\frac{1}{|k_1|^2}\theta(2^{-q}k_1)^2(\sum_{k_2}\int_0^{t_1}\frac{e^{-|k_2|^2(t_1-s)}}{|k_2|^2}[|k_1|^{\eta} \wedge|t_2-t_1|^{\frac{\eta}{2}}(|k_{12}|^{2\eta}+|k_2|^{2\eta})](t_1-s)^{-\frac{1-\eta}{2}}ds)^2\\\lesssim& \sum_{k_1}\frac{1}{|k_1|^2}\theta(2^{-q}k_1)^2(\sum_{k_2}\int_0^{t_1}\frac{e^{-|k_2|^2(t_1-s)}}{|k_2|^2}|k_1|^{\eta(1-\beta_0)} |t_2-t_1|^{\frac{\eta\beta_0}{2}}(|k_{12}|^{2\eta\beta_0}+|k_2|^{2\eta\beta_0})(t_1-s)^{-\frac{1-\eta}{2}}ds)^2  \\\lesssim& |t_1-t_2|^{\eta\beta_0/2}2^{q(1+2\eta(1+\beta_0))},\endaligned$$
which is the required estimate for $L_t^2$.

\subsection{Renormalisation for $u_2^{\varepsilon,i}u_2^{\varepsilon,j}$}
In this subsection we deal with $u_2^{\varepsilon,i}u_2^{\varepsilon,j}, i,j=1,2,3,$ and prove that $u_2^{\varepsilon,i}\diamond u_2^{\varepsilon,j}\rightarrow v_4^{ij}$ in $C([0,T];\mathcal{C}^{-\delta})$.
Recall that for $i,j=1,2,3$
$${u}_2^{\varepsilon,i}\diamond {u}_2^{\varepsilon,j}:={u}_2^{\varepsilon,i} {u}_2^{\varepsilon,j}-{C}_2^{\varepsilon,ij},$$
We have the following identities:
$$\aligned &u_2^{\varepsilon,i}u_2^{\varepsilon,j}:=L^1+L^2+L^3,\endaligned$$
where
$$\aligned L_t^1=&(2\pi)^{-\frac{9}{2}}\sum_{i_1,i_2,j_1,j_2=1}^3\sum_{k_{1234}=k}\int_0^t\int_0^te^{-|k_{12}|^2(t-s)-|k_{34}|^2(t-\bar{s})}
:\hat{X}_s^{\varepsilon,i_1}(k_1)\hat{X}^{\varepsilon,i_2}_s(k_2)\hat{X}^{\varepsilon,j_1}_{\bar{s}}(k_3)\hat{X}^{\varepsilon,j_2}_{\bar{s}}(k_4):dsd\bar{s}e_k
\\&\hat{P}^{ii_1}(k_{12})\imath k^{i_2}_{12}\hat{P}^{jj_1}(k_{34})
\imath k^{j_2}_{34}\endaligned$$
$$\aligned L^2_t=&\sum_{i=1}^4{I^i_t}\\=& (2\pi)^{-\frac{9}{2}}\sum_{i_1,i_2,j_1,j_2=1}^3\sum_{k_{24}=k,k_1}\int_0^t\int_0^te^{-|k_{12}|^2(t-s)-|k_{4}-k_1|^2(t-\bar{s})}\frac{f(\varepsilon k_1)^2e^{-|k_1|^2|s-\bar{s}|}}{2|k_1|^2}\\&:\hat{X}^{\varepsilon,i_3}_s(k_2)\hat{X}^{\varepsilon,j_3}_{\bar{s}}(k_4):dsd\bar{s}e_k\hat{P}^{ii_1}
(k_{12})\imath k^{i_2}_{12}\hat{P}^{jj_1}(k_{4}-k_1)\imath(k^{j_2}
_{4}-k_1^{j_2})\sum_{j_5=1}^3\hat{P}^{i_4j_5}(k_1)\hat{P}^{j_4j_5}(k_1)\\&(1_{i_3=i_2,i_4=i_1,j_3=j_2,j_4=j_1}+1_{i_3=i_2,i_4=i_1,j_3=j_1,j_4=j_2}
+1_{i_3=i_1,i_4=i_2,j_3=j_2,j_4=j_1}+1_{i_3=i_1,i_4=i_2,j_3=j_1,j_4=j_2}),\endaligned$$
and
$$\aligned L_t^3=&(2\pi)^{-6}\sum_{i_1,i_2,j_1,j_2=1}^3\sum_{k_1,k_2}\int_0^t\int_0^t
e^{-|k_{12}|^2(t-s+t-\bar{s})}\frac{f(\varepsilon k_1)^2f(\varepsilon k_2)^2e^{-(|k_1|^2+|k_2|^2)|s-\bar{s}|}}{4|k_1|^2|k_2|^2}dsd\bar{s}\hat{P}^{ii_1}(k_{12})\hat{P}^{jj_1}(k_{12})\\&\imath k_{12}^{i_2}(-\imath k_{12}^{j_2})\sum_{j_3,j_4=1}^3
(\hat{P}^{i_1j_3}(k_1)\hat{P}^{j_1j_3}(k_1)\hat{P}^{i_2j_4}(k_2)\hat{P}^{j_2j_4}(k_2)+\hat{P}^{i_1j_3}(k_1)\hat{P}^{j_2j_3}(k_1)\hat{P}^{i_2j_4}(k_2)
\hat{P}^{j_1j_4}(k_2)).\endaligned$$
Here each $I_t^i$ corresponds to the term associated with each indicator function respectively.
To make it more readable we write each term corresponding to the tree notation in Section 2. $u_2^{\varepsilon,i}u_2^{\varepsilon,j}$ corresponds to $\includegraphics[height=0.7cm]{07.eps} $ and $L_t^1, L_t^2, I_t^1, L_t^3$ correspond to the associated $\hat{\mathcal{W}}^{(\varepsilon,4)}, \hat{\mathcal{W}}^{(\varepsilon,2)}, \hat{\mathcal{W}}^{(\varepsilon,2)}_1, \hat{\mathcal{W}}^{(\varepsilon,0)}$ in the proof of Theorem 2.16 respectively.

By an easy computation we obtain that
$$\aligned L_t^3=&(2\pi)^{-6}\sum_{i_1,i_2,j_1,j_2=1}^3\sum_{k_1,k_2}f(\varepsilon k_1)^2f(\varepsilon k_2)^2\hat{P}^{ii_1}(k_{12})\hat{P}^{jj_1}(k_{12})k_{12}^{i_2}k_{12}^{j_2}
\sum_{j_3,j_4=1}^3(\hat{P}^{i_1j_3}(k_1)\hat{P}^{j_1j_3}(k_1)\hat{P}^{i_2j_4}(k_2)\\&\hat{P}^{j_2j_4}(k_2)+\hat{P}^{i_1j_3}(k_1)\hat{P}^{j_2j_3}(k_1)\hat{P}^{i_2j_4}(k_2)\hat{P}^{j_1j_4}(k_2))
\frac{1}{2|k_1|^2|k_2|^2(|k_1|^2+|k_2|^2+|k_{12}|^2)}\\&\bigg[\frac{1-e^{-2|k_{12}|^2t}}{2|k_{12}|^2}-
\int_0^te^{-2|k_{12}|^2(t-s)-(|k_1|^2+|k_2|^2+|k_{12}|^2)s}ds\bigg].\endaligned$$
Let $$\aligned C_2^{\varepsilon,ij}(t)=&L_t^3.\endaligned$$

\textbf{Terms in the second chaos}:
Now we come to $L_t^2$: it is sufficient to consider $I^1_t$
 and the
desired estimates for the other terms can be obtained similarly. For $\epsilon>0$ small enough we have the following inequality
$$\aligned &E|\Delta_qI_t^1|^2\\\lesssim &\sum_k\sum_{k_{24}=k,k_{24}'=k,k_1,k_1'}\theta(2^{-q}k)^2\int_0^t\int_0^t\int_0^t\int_0^t
e^{-|k_{12}|^2(t-\sigma)-|k_{4}-k_1|^2(t-\bar{\sigma})}|k_{12}(k_4-k_1)|\\&e^{-|k_{12}'|^2(t-s)-|k_{4}'-k_1'|^2(t-\bar{s})} |k_{12}'(k_4'-k_1')|\frac{1}{|k_1|^2|k_1'|^2|k_2|^2|k_4|^2}1_{\{k_2=k_2',k_4=k_4'\}}+1_{\{k_2=k_4',k_4=k_2'\}}dsd\bar{s}d\sigma d\bar{\sigma}
,\endaligned$$
Now in the following we only estimate the  term corresponding to the first characteristic function on the right hand side of the inequality. The second term can be estimated similarly:
$$\aligned E|\Delta_qI_t^1|^2\lesssim &\sum_k\sum_{k_{24}=k,k_1,k_3}\theta(2^{-q}k)^2\int_0^t\int_0^t\int_0^t\int_0^te^{-|k_{12}|^2(t-s)-|k_{4}-k_1|^2(t-\bar{s})-|k_{23}|^2(t-\sigma)-|k_{4}-k_3|^2(t-\bar{\sigma})}
\\&\frac{1}{|k_1|^2|k_2|^2|k_3|^2|k_4|^2}dsd\bar{s}|k_{12}(k_{4}-k_1)k_{23}(k_4-k_3)|
\\\lesssim & t^\epsilon\sum_k \sum_{k_{24}=k,k_1,k_3}\frac{\theta(2^{-q}k)^2}{|k_1|^2|k_2|^2|k_3|^2|k_4|^2|k_1-k_4|^{1-\epsilon}|k_4-k_3||k_{12}|^{1-\epsilon}|k_{23}|}
\endaligned$$
$$\aligned\lesssim & t^\epsilon\sum_k\sum_{k_{24}=k}\frac{\theta(2^{-q}k)^2}{|k_2|^2|k_4|^2}\sum_{k_1}\frac{1}{|k_1-k_4||k_1|^2|k_{12}|}\sum_{k_3}\frac{1}{|k_3-k_4||k_3|^2|k_{23}|}\\\lesssim & t^\epsilon\sum_k\sum_{k_{24}=k}\frac{\theta(2^{-q}k)^2}{|k_2|^2|k_4|^2}(\sum_{k_1}\frac{1}{|k_1-k_4|^{2-2\epsilon}|k_1|^2})^{1/2}(\sum_{k_1}\frac{1}{|k_{12}|^{2-2\epsilon}|k_1|^2})^{1/2}\\&(\sum_{k_3}\frac{1}{|k_3-k_4|^2|k_3|^2})^{1/2}(\sum_{k_3}\frac{1}{|k_{23}|^2|k_3|^2})^{1/2}\\\lesssim &t^\epsilon\sum_k\sum_{k_{24}=k}\frac{\theta(2^{-q}k)^2}{|k_2|^{3-\epsilon}|k_4|^{3-\epsilon}}\lesssim t^\epsilon2^{2q\epsilon},\endaligned$$
where in the last two inequalities we used Lemma 3.10.

\textbf{Terms in the fourth chaos}:

Now we consider $L_t^1$. For $\epsilon>0$ small enough we have the following calculations:
$$\aligned &E|\Delta_qL_t^1|^2\\\lesssim&\sum_k\sum_{k_{1234}=k}\theta(2^{-q}k)^2\bigg[\int_0^t\int_0^t\int_0^t\int_0^te^{-|k_{12}|^2(t-s+t-\sigma)-|k_{34}|^2(t-\bar{s}+t-\bar{\sigma})}\frac{e^{-(|k_1|^2+|k_2|^2)|s-\sigma|-(|k_3|^2+|k_4|^2)|\bar{s}-\bar{\sigma}|}}{|k_1|^2|k_2|^2|k_3|^2|k_4|^2}\\&dsd\bar{s}d\sigma d\bar{\sigma}|k_{12}k_{34}|^2+\int_0^t\int_0^t\int_0^t\int_0^te^{-|k_{12}|^2(t-s)-|k_{23}|^2(t-\sigma)-|k_{34}|^2(t-\bar{s})-|k_{14}|^2(t-\bar{\sigma})}\frac{1}{|k_1|^2|k_2|^2|k_3|^2|k_4|^2}\\&dsd\bar{s}d\sigma d\bar{\sigma}|k_{12}k_{34}k_{14}k_{23}|\bigg]\\\lesssim&t^\epsilon\sum_k\sum_{k_{1234}=k}\theta(2^{-q}k)^2
\bigg(\frac{1}{|k_1|^2|k_2|^2|k_3|^2|k_4|^2|k_{12}|^{2-\epsilon}|k_{34}|^{2-\epsilon}}\endaligned$$
$$\aligned&+\frac{1}{|k_1|^2|k_2|^2|k_3|^2|k_4|^2|k_{12}|^{1-\epsilon/2}
|k_{34}|^{1-\epsilon/2}|k_{14}|^{1-\epsilon/2}|k_{23}|^{1-\epsilon/2}}\bigg)
\\\lesssim&t^\epsilon\bigg[2^{2q\epsilon}+\bigg(\sum_{k_{1234}=k}\frac{\theta(2^{-q}k)^2 }{|k_1|^2|k_2|^2|k_3|^2|k_4|^2|k_{12}|^{2-\epsilon}|k_{34}|^{2-\epsilon}}\bigg)^{1/2}\bigg(\sum_{k_{1234}=k}\frac{\theta(2^{-q}k)^2 }{|k_1|^2|k_2|^2|k_3|^2|k_4|^2|k_{14}|^{2-\epsilon}|k_{23}|^{2-\epsilon}}\bigg)^{1/2}\bigg]\\\lesssim&t^\epsilon2^{2q\epsilon},\endaligned$$
where we used Lemma 3.10 in the last inequality.
By a similar calculation we also get that for $\epsilon,\eta>0$, $\gamma>0$ small enough
$$\aligned &E[|\Delta_q(u_2^{\varepsilon_1,i}\diamond u_2^{\varepsilon_1,j}(t_1)-u_2^{\varepsilon_1,i}\diamond u_2^{\varepsilon_1,j}(t_2)-u_2^{\varepsilon_2,i}\diamond u_2^{\varepsilon_2,j}(t_1)+u_2^{\varepsilon_2,i}\diamond u_2^{\varepsilon_2,j})(t_2))|^2]\\\lesssim & (\varepsilon_1^{2\gamma}+\varepsilon_2^{2\gamma})|t_1-t_2|^\eta2^{q(\epsilon+2\eta)},\endaligned$$
which together with Gaussian hypercontractivity, Lemma 3.1 and similar arguments as for (4.3) implies that there exist $v_4^{ij}\in C([0,T];\mathcal{C}^{-\delta}),i,j=1,2,3$ such that for $p>1$
$$u_2^{\varepsilon,i}\diamond u_2^{\varepsilon,j}\rightarrow v_4^{ij}\textrm{ in } L^p(\Omega,P,C([0,T];\mathcal{C}^{-\delta})).$$
 Here $\delta>0$ depending on $\eta, \epsilon, p$ can be chosen small enough.

\end{document}